%     Ce ficher est du plain TeX.  
%     An   introduction to finite volumes
%     Fran\c{c}ois Dubois, 15 octobre 2000,  26, 27 octobre 2001. 
%     asci le 26, 27, 28, 29    aout 2003.
%     Versailles le 27 janvier, 25 fevrier  2005.
%     Edition 'hal', 13, 15, 21 janvier 2011 
%%%%%%%%%%%%%%%%%%%%%%%%%%%%%%%%%%%%%%%%%%%%%%%%%%%%%%%%%%%%%%%%%%%
%                       thanks to Gabriel Turinici
\input epsfx.tex
%                       thanks to Gabriel Turinici
%%%%%%%%%%%%%%%%%%%%%%%%%%%%%%%%%%%%%%%%%%%%%%%%%%%%%%%%%%%%%%%%%%%
\overfullrule=0pt

% \nopagenumbers

%%%%%%%%%%%%%%%%%%%%%%%%%%%%%%%%%%%%%%%%%   reglages 2005 
%  taille d'agrandissement  :
%  \magnification=1400
%  \hsize=12.5cm
%  \vsize=16.5cm  
%    marges : 
%   \hoffset=-.6cm
%   \voffset=.5cm
%%%%%%%%%%%%%%%%%%%%%%%%%%%%%%%%%%%%%%%%%   fin reglages 2005 

%taille d'agrandissement  :
\magnification=1400
\hsize=12.2cm
\vsize=16.0cm
% marge gauche : \hoffset=-1cm
\hoffset=-0.4cm
\voffset=.8cm
\baselineskip 16 true pt

\font \smcaps=cmbx10 at 12 pt
\font \pecaps=cmcsc10 at 9 pt
\font \fauteur=cmbx10 scaled 1167 
\font \gcaps=cmbx10 scaled 1600

%  Pagination d'apres Raymond S\'eroul page 231 et 68 
\newtoks \hautpagegauche  \hautpagegauche={\hfil}
\newtoks \hautpagedroite  \hautpagedroite={\hfil}
\newtoks \titregauche     \titregauche={\hfil}
\newtoks \titredroite     \titredroite={\hfil}
\newif \iftoppage        \toppagefalse   
\newif \ifbotpage         \botpagefalse    
\titregauche={\pecaps    Fran\c{c}ois Dubois }
\titredroite={\pecaps   An introduction to finite volumes for gas dynamics} 
\hautpagegauche = { \hfill \the \titregauche  \hfill  }
\hautpagedroite = { \hfill \the \titredroite  \hfill  }
\headline={ \vbox  { \line {  
\iftoppage    \ifodd  \pageno \the \hautpagedroite  \else \the
\hautpagegauche \fi \fi }     \bigskip  \bigskip  }}
\footline={ \vbox  {   \bigskip  \bigskip \line {  \ifbotpage  
\hfil {\oldstyle \folio} \hfil  \fi }}}

%petite boullette

%carr? de fin de d?monstration
\def\sqr#1#2{{\vcenter{\vbox{\hrule height.#2pt \hbox{\vrule width .#2pt height#1pt 
\kern#1pt \vrule width.#2pt} \hrule height.#2pt}}}}
\def\square{\mathchoice\sqr64\sqr64\sqr{4.2}3\sqr33} 

% Les nombres math?matiques
\def\N{{\rm I}\! {\rm N}}
\def\R{{\rm I}\! {\rm R}}
\def\Z{{\rm Z}\!\! {\rm Z}}

% indice bas

% les valeurs absolues, les modules et les nornes
\def\abs#1{\mid \! #1 \! \mid }

\def\mod#1{\setbox1=\hbox{\kern 3pt{#1}\kern 3pt}%
\dimen1=\ht1 \advance\dimen1 by 0.1pt \dimen2=\dp1 \advance\dimen2 by 0.1pt
\setbox1=\hbox{\vrule height\dimen1 depth\dimen2\box1\vrule}%
\advance\dimen1 by .1pt \ht1=\dimen1
\advance \dimen2 by .01pt \dp1=\dimen2 \box1 \relax}

\def\nor#1{\setbox1=\hbox{\kern 3pt{#1}\kern 3pt}%
\dimen1=\ht1 \advance\dimen1 by 0.1pt \dimen2=\dp1 \advance\dimen2 by 0.1pt
\setbox1=\hbox{\kern 1pt  \vrule \kern 2pt \vrule height\dimen1 depth\dimen2\box1
\vrule
\kern 2pt \vrule \kern 1pt  }%
\advance\dimen1 by .1pt \ht1=\dimen1
\advance \dimen2 by .01pt \dp1=\dimen2 \box1 \relax}

\rm
\hyphenation {per-tur-ba-tion com-pres-sible com-pres-si-bles  me-cha-nics 
sprin-ger con-ve-na-ble
si-mu-la-tion si-tua-tion uni-di-men-sio-nal com-pres-sible com-pres-si-bles  
me-cha-nics 
sprin-ger con-ve-na-ble si-mu-la-tion si-tua-tion com-pres-sible 
com-pres-si-bles  me-cha-nics 
sprin-ger con-ve-na-ble si-mu-la-tion si-tua-tion scien-ti-fi-ques }
%%%%%%%%%%%%%%%%%%%%%%%%%%%%%%%%%%%%%%%%%%%%%%%%%%%%%%%%%%%%%%%%%%%%%%%%%%%%%%%%%%%%%%%%%%%
% \nopagenumbers
%  \null\vskip 0.1 cm 

$\,$

\bigskip \bigskip \bigskip

\centerline{\gcaps     An  introduction  to finite volumes  } 
\bigskip 
\centerline {\gcaps     for gas dynamics  }

\bigskip \bigskip \bigskip

\centerline { \fauteur  Fran\c{c}ois Dubois }

\smallskip
\smallskip
\centerline { Conservatoire National des Arts et M\'etiers }
\centerline {  15 rue Marat, F-78 210 $\,$ Saint Cyr
l'Ecole,  France.  } 

\smallskip
\centerline {  and }
\centerline { Centre National de la Recherche Scientifique}
\centerline {   Laboratoire ASCI,  b\^at. 506, BP
167,  F-91~403 Orsay  Cedex.  }

\bigskip
\bigskip
\centerline {  October   2000  
\footnote{$ ^{^{\scriptstyle  \square}}$}  {\rm 
Research report CNAM-IAT n$^{\rm o}$342-2000.  
Published with the title  ``An  Introduction  to Finite Volumes 
Methods'' in  {\it Encyclopedia Of Life Support Systems} 
(EOLSS, Unesco),  
Mathematical Sciences,   Computational Methods and Algorithms, 
Vladimir V. Shaidurov and  Olivier Pironneau Editors, volume 2, 
p.~36-105, 2009.  Present edition   21 January 2011.  }}

\bigskip
\bigskip  \noindent {\bf  Summary }

\medskip

We propose an elementary introduction to the finite volume method in the context of 
gas dynamics conservation laws. Our approach is founded on the advection equation, the exact
integration of the associated Cauchy problem, and the so-called upwind scheme in one
space dimension. It is then extended in three directions~: hyperbolic linear systems
and particularily the system of acoustics, gas dynamics with the help of the Roe
matrix and two space dimensions by following 
 the approach proposed by Van Leer. A special
emphasis on boundary conditions is proposed all along the text.

\medskip \noindent 
{\bf AMS Subject Classification: } 35L40, 35L60, 35L65, 35Q35, 76N15.
  
\medskip \noindent 
{\bf Keywords:} Advection, Characteristics, Roe matrix, Van Leer method.

\vfill \eject 

%%%%%%%%%%%%%%%%%%%%%%%%%%%%         ajout janvier 2005
\toppagetrue  
\botpagetrue   
%%%%%%%%%%%%%%%%%%%%%%%%%%%%         ajout janvier 2005 
$ \, $ 

\bigskip 
\bigskip 
\bigskip  

\noindent {\smcaps  Contents }

\medskip

\noindent   {\bf 1)  $ \,\, $    Advection equation and  method of
characteristics  }

\noindent $ \qquad 1.1 \quad $  Advection equation $\hfill $ 3

\noindent $ \qquad 1.2 \quad $   Initial-boundary value problems for the advection
equation  $\hfill $ 4

\noindent $ \qquad 1.3 \quad $    Inflow and outflow for the advection equation  $\hfill $ 6

\smallskip  \noindent   {\bf 2)  $ \,\, $    Finite volumes for linear hyperbolic
systems   }

\noindent $ \qquad 2.1 \quad $    Linear advection  $\hfill $ 8

\noindent $ \qquad 2.2 \quad $    Numerical flux boundary conditions  $\hfill $ 13

\noindent $ \qquad 2.3 \quad $    A model system with two equations  $\hfill $ 14

\noindent $ \qquad 2.4 \quad $    Unidimensional linear acoustics  $\hfill $ 17 

\noindent $ \qquad 2.5 \quad $    Characteristic variables  $\hfill $ 21

\noindent $ \qquad 2.6 \quad $    A family of model systems with three equations  $\hfill $ 25

\noindent $ \qquad 2.7 \quad $    First order upwind-centered finite volumes  $\hfill $ 27

\smallskip  \noindent   {\bf 3)  $ \,\, $     Gas dynamics with  the Roe method }

\noindent $ \qquad 3.1 \quad $   Nonlinear acoustics in one space dimension  $\hfill $ 29

\noindent $ \qquad 3.2 \quad $   Linearization of the gas dynamics equations  $\hfill $ 30
 
\noindent $ \qquad 3.3 \quad $   Roe matrix  $\hfill $ 33 

\noindent $ \qquad 3.4 \quad $   Roe flux   $\hfill $ 35

\noindent $ \qquad 3.5 \quad $   Entropy correction   $\hfill $ 38

\noindent $ \qquad 3.6 \quad $   Nonlinear  flux boundary conditions   $\hfill $ 40

\smallskip  \noindent   {\bf 4)  $ \,\, $    Second order and  two space dimensions }

\noindent $ \qquad 4.1 \quad $   Towards second order accuracy   $\hfill $ 42

\noindent $ \qquad 4.2 \quad $  The method of lines  $\hfill $ 43

\noindent $ \qquad 4.3 \quad $  The method of Van Leer  $\hfill $ 45

\noindent $ \qquad 4.4 \quad $  Second order accurate finite volume method for fluid problems
   $\hfill $ 49 

\noindent $ \qquad 4.5 \quad $  Explicit Runge-Kutta integration with respect to time 
$\hfill $ 54

\smallskip  \noindent   {\bf 5)  $ \,\, $    References }  $\hfill $ 55

\bigskip  \bigskip \noindent   $\bullet \qquad \,\,\, $ 
{\it Acknowledgments}. \quad The author thanks Alexandre Gault, listener at the
spring 2000  ``lectures in computational acoustics'' at the  Conservatoire National des
Arts et M\'etiers (Paris, France), for providing his personal  manuscript  notes.

\vfill \eject 

\bigskip 
\bigskip 
\bigskip 

\noindent  {\smcaps 1) $ \quad $ Advection equation and  method of characteristics. }
\smallskip \noindent {\smcaps 1.1 } $ \,  $ { \bf  Advection equation. }

\noindent   $\bullet \qquad \,\,\,\,\, $ 
We consider a given real number $\,\, a > 0 \,\, $ and we wish to solve the so-called
advection equation of unknown function $\,\, u(x,\,t) \,\,$~: 

\noindent  (1.1.1) $\qquad \displaystyle 
{{\partial u}\over{\partial t}} \,\,+\,\, a \,\, {{\partial u }\over{\partial
x}}\,\,=\,\, 0 \,\,, \qquad t \geq 0 \,,\quad x \in \R \,. \,$

\noindent
We  first look to  the homogeneity coherence of the different terms of equation
(1.1.1). On one hand,  the ratio $\,\, {{\partial u}\over{\partial t}} \,\,$ is
homogeneous to the dimension $\, [u] \,$ of function $\,\, u({\scriptstyle \bullet},\,
{\scriptstyle \bullet}) \,\,$ divided by the dimension $\, [t] \,$ of the time and we
have~: $\,\,  {{\partial u}\over{\partial t}} \sim {{[u]}\over{[t]}} . \,\,$ On the
other hand  the expression   $\,\, a \, {{\partial u}\over{\partial x}} \,\,$ is
homogeneous to the  dimension $\,\, [a]\,\,$  of scalar $\, a \,$ multiplied by the
ratio $\,\, {{[u]}\over{[x]}} \,\,$ and we have $\,\, a \,  {{\partial
u}\over{\partial x}} \sim  [a]\, {{[u]}\over{[x]}} \,.\, $ From equation (1.1.1), the
two previous terms  $\,\,  {{\partial u}\over{\partial t}} \,\,$ and $\,\, a \,
{{\partial u }\over{\partial x}}\,\,$ have the same dimension and we deduce from the
previous formulae the equality~: $\,\, {{1}\over{[t]}} \sim  {{[a]}\over{[x]}} \,. \,$
Then we have established that the constant  $\, a \,$ is homogeneous to a {\bf
celerity}~: 

\noindent  (1.1.2) $\qquad \displaystyle 
[a]\,\, \sim \,\,  {{[x]}\over{[t]}} \,. \,$ 

\smallskip \noindent   $\bullet \qquad \,\,\, $ 
The Cauchy problem for the  model equation (1.1.1)   is composed by the equation
(1.1.1) itself and the following initial condition~: 

\noindent  (1.1.3) $\qquad \displaystyle 
u(x, \,0) \,\,= \,\, u_0(x) \,,\qquad x \in \R \,,\,$ 

\noindent
where $\,\, \R \, \ni x \longmapsto u_0(x) \in \R \,\, $ is some given function. 
We observe that the solution of equation (1.1.1) is constant along the characteristic
(straight) lines that satisfy the differential equation 

\noindent  (1.1.4) $\qquad \displaystyle 
{{{\rm d}x}\over{{\rm d}t}} \,\,= \,\, a \,.\, $

\bigskip  \noindent  {\bf Proposition 1.1.    The solution is constant along
the characteristic lines.}

\noindent 
Let $\,\, 0 \leq  \lambda \leq  t  \,\,$ be some given parameter and $\,\,
u({\scriptstyle \bullet},\, {\scriptstyle \bullet}) \,\,$ a  solution of equation
(1.1.1). Then function $\,\,  u({\scriptstyle \bullet},\, {\scriptstyle \bullet})
\,\,$ is constant along the characteristic lines, {\it i.e.} 

\noindent  (1.1.5) $\qquad \displaystyle 
u(x-a \lambda,\, t- \lambda) \,\,= \,\, u(x,\,t) \,\,,\qquad \forall \,\,\, 
x ,\, t ,\, \lambda \,. \,$

\smallskip  \noindent   $\bullet \qquad \,\,\, $ 
The {\bf proof of Proposition 1.1} is obtained as follows.  We consider a fixed point $\,\,
(x,\,t) \,\,$ in space-time $\,\, \R \times [0,\, +\infty[ \,\,$ and the  auxiliary
function $\,\, [0,\, t] \, \ni \lambda \longmapsto v(\lambda) \,=\, u (x-a \lambda,\,
t- \lambda) \,. \,$ We have, due to the usual chain rule for derivation of  operators~: 

\noindent  $ \displaystyle 
{{{\rm d}v}\over{{\rm d}\lambda}} \,\,=\,\, \Bigl[ (-a) \, {{\partial u}\over{\partial
x}} \,-\,  {{\partial u}\over{\partial t}} \, \Bigr] \, (x-a \lambda,\, t- \lambda)
\,\,= \,\,0 \,\,$ \qquad  if function $\,\, u({\scriptstyle \bullet},\, {\scriptstyle
\bullet}) \,\,$ is solution of the advection equation (1.1.1). Then $\,\, v(\lambda)
\,\,$ does not depend on variable $\,\, \lambda \,\,$ and we have in particular $\,\,
v(\lambda) \,=\, v(0) \,,\,$ which exactly expresses the relation (1.1.5). We have in
particular for $\,\, \lambda = t\,$~:  $\,\,u(x,\,t) \,=\, u(x-at,\,0) \,=\, u_0(x-at)
\,\,$ as illustrated on Figure 1.1. $ \hfill \square \kern0.1mm $

\smallskip \centerline {  \epsfysize=5,0cm  \epsfbox  {fig.1.1.epsf} } \smallskip  

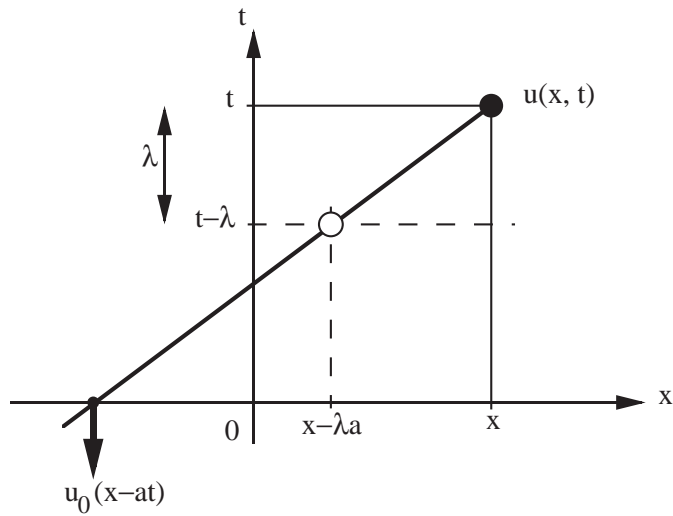
\centerline {\rm  {\bf Figure 1.1.}	\quad The solution  $\,\, u(x, t)\,\,$  of the
advection equation}

\centerline {\rm   is constant along the characteristic lines. }
\bigskip 

\smallskip \centerline {  \epsfysize=4,5cm  \epsfbox  {fig.1.2.epsf} } \smallskip  

\centerline {\rm  {\bf Figure 1.2.}	\quad Initial-boundary value problem for the
advection equation. }
\bigskip

\bigskip \noindent {\smcaps 1.2} $ \,\,\,\,  $ { \bf  Initial-boundary value
problems for the advection equation. }

\noindent   $\bullet \qquad \,\,\,\,\, $ 
The second step is concerned by  the so-called initial-boundary value problem
considered for $\, x>0 \,$ and $t>0 \,$  with some given  initial condition $\,\, 
u_0(x) \,\,$  for $\, t=0 \,$ and a boundary condition $\,\, v_0(t) \,\,$  for $\,
x=0 \,$~: 
  
\noindent  (1.2.1) $\qquad \displaystyle 
{{\partial u}\over{\partial t}} \,\,+\,\, a \,\, {{\partial u }\over{\partial
x}}\,\,=\,\, 0 \,\,, \qquad t > 0 \,,\quad x > 0\,,  \qquad \,$ (equation) 
  
\noindent  (1.2.2) $\qquad \displaystyle 
u(x, \,0) \,\,= \,\, u_0(x) \,,\qquad \quad x > 0 \,,\,$ \qquad  \qquad  (initial
condition) 

\noindent  (1.2.3) $\qquad \displaystyle 
u(0, \,t) \,\,\, = \,\, v_0(t) \,,\qquad \quad \,\, t > 0 \,,\,$ \qquad  \qquad 
(boundary condition).

\smallskip 

\noindent  {\bf Proposition 1.2. $\quad$  Advection in the quadrant 
$\,\, x > 0 \,$  and $\,\, t > 0 \,.$ }

\noindent 
We suppose that $\,\, a > 0 \,.$ Then  the solution of the advection equation (1.2.1)
with the initial condition (1.2.2) and the boundary condition (1.2.3) is given by the
relations 

\noindent  (1.2.4) $\qquad \displaystyle 
u(x, \,t) \,\,\, = \,\, u_0(x - at) \,\,, \qquad \, x-at \,> \, 0 \,$ 

\noindent  (1.2.5) $\qquad \displaystyle 
u(x, \,t) \,\,\, = \,\, v_0 \Bigl( t \,- {{x}\over{a}}\Bigr) \,\,, \qquad x-at \, < \,
0 \,. \,$ 

\noindent 
The initial condition $\,\, u_0({\scriptstyle \bullet}) \,\,$ is advected towards
space-time point $\, (x,\,t) \,$ when $\, x-at > 0 \,$ and the boundary condition
$\,\, v_0({\scriptstyle \bullet}) \,\,$  is activated for  $\, x-at < 0 \,. \,$

\smallskip  \noindent   $\bullet \qquad \,\,\, $ {\bf Proof of Proposition 1.2.}

\noindent 
In order to solve the problem (1.2.1)-(1.2.3), we use the method of characteristics.
We fix a point $\,\, (x,\,t) \,$ of space-time domain that satisfies $\,\, x>0 ,\,
t>0 \,\,$ and we go upstream in time with the help of the characteristic line that goes
through this point (see Figure 1.2)~: 

\noindent  (1.2.6) $\qquad \displaystyle 
 x(\lambda) \,= \, x - a \lambda \,,\,\qquad  t(\lambda) \,=\, t - \lambda \,.\,$  
  
\noindent   $\bullet \qquad \,\,\, $ 
First case~: $\,\, x-at \,> \,0 \,. \, $  When we take the particular value $\, \lambda
= t \,$ in the previous relation (1.2.6), the particular point $\, y \,=\, x(t) \,=\,
x-at \,\,$ on the axis of abscissa  is strictly positive then the initial condition
$\, u_0(y) \,$ is well defined. The solution $\,\, u({\scriptstyle \bullet} ,\,
{\scriptstyle \bullet}) \,\, $ is constant on the characteristic line (see
Proposition 1.1) that contains this particular point. Then relation (1.2.4) is
established. 

\noindent   $\bullet \qquad \,\,\, $ 
Second case~: $\,\, x-at \,< \,0 \,. \, $  We consider the particular value $\,\,
\lambda = {{x}\over{a}} \,$ inside the expression (1.2.6). Then the corresponding
foot of the characteristic belongs to the time axis~: $\,\, \theta \,=\, t-\lambda =
t- {{x}\over{a}} \,$ and $\,\, \theta > 0 \,\,$ due to the inequalities $\, x < a t
\,$ and $\, a > 0 \,.\,$  The solution is constant along the characteristic line
going through this point and the relation (1.2.5) is established.   $ \hfill \square
\kern0.1mm $
  
\smallskip  \noindent   $\bullet \qquad \,\,\, $ 
In the particular case where datum $\,\, u_0(x) \,\,$ is identically equal to zero,
{\it i.e.} 
  
\noindent  (1.2.7) $\qquad \displaystyle 
u_0(x) \,\,= \,\, 0 \,\,,\qquad \qquad \,\, x > 0 \,,\,$

\noindent
and if the boundary condition $\, \, v_0(t) \,\,$ is sinuso\"{\i}dal for time positive
to fix the ideas, 

\noindent  (1.2.8) $\qquad \displaystyle 
v_0(t) \,\,\,= \,\, {\rm sin}(\omega t)  \,\,,\qquad t > 0 \,,\,$

\noindent
the solution of the advection equation in the domain $\,\, x > 0 \,,\,\, t>0 \,\,$
via the relations (1.2.4) and (1.2.5) can be considered with the two following view
points. 

\smallskip \noindent
{\bf (i) } \qquad We take a snap shot  of the solution $\,\, u({\scriptstyle \bullet} ,\,
{\scriptstyle \bullet}) \,\,$ at a fixed time $\, T>0 .\,$ We consider  the partial
function $\,\, [0,\,+\infty[ \, \, \ni x \longmapsto u(x,\,T) \in \R \,\,$ and taking into
account the relations (1.2.4),  (1.2.5), (1.2.7) and   (1.2.8), we have

\setbox21=\hbox {$\displaystyle  {\rm sin} \Bigl[ \omega \bigl( T-{{x}\over{a}}
\bigr)  \Bigr]  \,\,, \qquad x \,< \, a T  \,  $}
\setbox22=\hbox {$\displaystyle  0 \,\,, \qquad \qquad \qquad  \qquad x \,> \, a T\, 
\,.\,   $}
\setbox30= \vbox {\halign{#&# \cr \box21 \cr \box22    \cr   }}
\setbox31= \hbox{ $\vcenter {\box30} $}
\setbox44=\hbox{\noindent  (1.2.9) $\displaystyle  \qquad   u(x,\,T) \,\,= \,\, 
\left\{ \box31 \right. $}  

\noindent $ \box44 $

\noindent
and this function is illustrated on Figure 1.3. 

\smallskip \noindent
{\bf (ii) } \qquad We fix a particular position  $\,\, X \,\,$ in space and we look, as
time is increasing, to the solution $\,\, u({\scriptstyle \bullet} ,\, {\scriptstyle
\bullet}) \,\,$  at this particular point. We show on Figure 1.4 the function $\,\, 
[0,\,+\infty[ \, \ni t \longmapsto u(X,\,t) \in \R \,\,$ and taking into account the
relations (1.2.4),  (1.2.5), (1.2.7) and   (1.2.8), we have 

\setbox21=\hbox {$  0 \,\,, \qquad \qquad \qquad  \qquad t \,< \,
{{X}\over{a}}\,  \,\,   $}
\setbox22=\hbox {$  {\rm sin} \Bigl[ \omega \bigl( T-{{x}\over{a}}
\bigr)  \Bigr]  \,\,, \qquad t \,> \, {{X}\over{a}}\, \,  \,.\,   $}
\setbox30= \vbox {\halign{#&# \cr \box21 \cr \box22    \cr   }}
\setbox31= \hbox{ $\vcenter {\box30} $}
\setbox44=\hbox{\noindent  (1.2.10) $\displaystyle  \qquad   u(x,\,T) \,\,= \,\, 
\left\{ \box31 \right. $}  

\noindent $ \box44 $
 
% version mac 
% \smallskip \vskip 3.7cm  \smallskip 
%  \qquad  \qquad  \qquad  \qquad     \special{illustration fig.1.3.epsf scaled 500}  
% fin de la version version mac 
% version linux
\smallskip  \centerline {  \epsfysize=4cm  \epsfbox  {fig.1.3.epsf} }  \smallskip 
% fin de la version linux

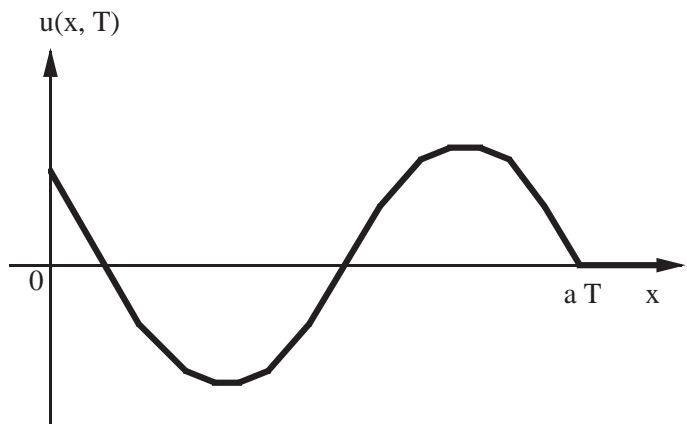
\centerline {\rm  {\bf Figure 1.3.}	\quad Snap shot of the solution of the
advection equation } 

\centerline {\rm at time  $\, t=T .\,$  }
 
\smallskip \centerline {  \epsfysize=4cm  \epsfbox  {fig.1.4.epsf} }  \smallskip
% fin de la version linux

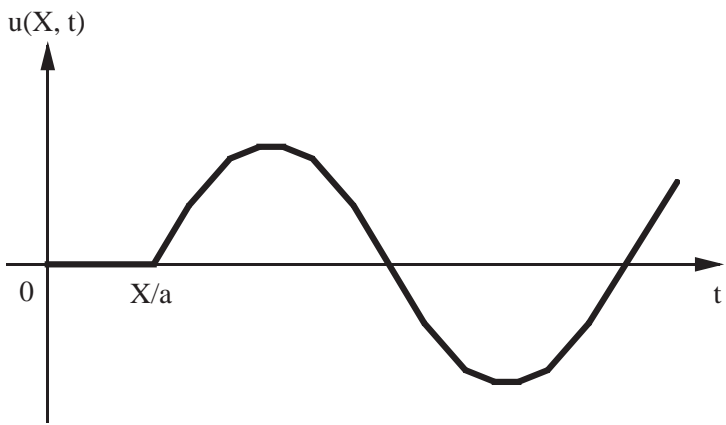
\centerline {\rm  {\bf Figure 1.4.}	\quad  Evolution of the solution at the
particular point  $\, x = X .\,$  }

\bigskip 
\smallskip  \noindent {\smcaps 1.3} $ \,\,\,\,  $ { \bf  Inflow and outflow
for the advection equation. }

\noindent   $\bullet \qquad \,\,\, $ 
We still suppose that celerity $\,a \,$ is positive and we consider the resolution of
the advection (1.2.1) in the space-time domain 

\noindent  (1.3.1) $\qquad \displaystyle 
0 \,< \, x \, < \, L \,, \qquad t \, > \, 0 \,. \,$

\noindent 
The relations (1.2.4) and (1.2.5) can still be applied because the proof of
Proposition 1.2 remains unchanged in this particular case. As a consequence of the
previous property, we remark that {\bf no boundary condition} is necessary at the
particular position $ \,\, x=L \,\,$ for solving the advection problem in the
space-time domain defined in relations (1.3.1). The initial condition (1.2.2) has
simply  to be restricted in domain $\,\, ]0,\, L [ \,\,$~: 

\noindent  (1.3.2) $\qquad \displaystyle 
u(x, \,0) \,\,= \,\, u_0(x) \,,\qquad \quad 0 < x < L  \,,\,$ 

\noindent 
and the boundary condition (1.2.3) at $\,\, x=0 \,\, $ remains unchanged~: 

\noindent  (1.3.3) $\qquad \displaystyle 
u(0, \,t) \,\,\, = \,\, v_0(t) \,,\qquad \quad \,\, t > 0 \,.\,$ 

\smallskip 
% version mac 
% \smallskip \vskip 2.7cm  \smallskip   
% \qquad  \qquad  \qquad  \quad   \special{illustration fig.1.5.epsf scaled 500}  
% fin de la version version mac 
% version linux
\smallskip \centerline {  \epsfysize=3,7cm  \epsfbox  {fig.1.5.epsf} }  \smallskip
% fin de la version linux

\centerline {\rm  {\bf Figure 1.5.}	\quad    Initial-boundary value problem for the
advection equation }

\centerline {\rm  with  $\,\, a > 0 \,\,$  in the domain  $\,\, 0 < x < L \,\,$  and
$\,\, t > 0. \,\,$  }
\bigskip 

% version mac 
% \smallskip \vskip 3.5 cm  \smallskip 
% \qquad  \qquad  \qquad  \quad \special{illustration fig.1.6.epsf scaled 500}  
% fin de la version version mac 
% version linux
\smallskip \centerline {  \epsfysize=4,5cm  \epsfbox  {fig.1.6.epsf} } \smallskip  
% fin de la version linux 

\centerline {\rm  {\bf Figure 1.6.}	\quad    Initial-boundary value problem for the
advection equation }

\centerline {\rm  with  $\,\, a < 0 \,\,$  in the domain  $\,\, 0 < x < L \,\,$  and$\,\,
t > 0. \,\,$  }
\bigskip 

\smallskip \noindent   $\bullet \qquad \,\,\, $ 
The difference between point $\, x=0 \,$ and point $\, x=L \,$ for the resolution of
the advection equation in space-time domain (1.3.1) is due to the fact that we choose
an orientation of the characteristic lines $\,\, x-at \,=\, {\rm constant} \,\,$
associated to an { \bf increase} for the time direction. With this choice of time
direction, the characteristic lines {\bf enter} inside the space-time domain (1.3.1)
at $\,\, x=0 \,\,$ and they go outside at $\,\, x \!=\! L \,. \, $ The boundary
condition (1.3.3) is given at the input of the domain (see Figure 1.5) and at $\,\,x
\!=\! L \,,\,$ there is a free output from space time domain (1.3.1), without necessity
to specify any numerical boundary condition. 

\smallskip \noindent   $\bullet \qquad \,\,\, $ 
If we change the sign of celerity $\, a ,\,$ {\it i.e.} if we suppose now 

\noindent  (1.3.4) $\qquad \displaystyle 
a < 0  \,\,, \,$

\noindent 
the above analysis remains unchanged, but the algebraic relations (1.2.4) and (1.2.5)
have to be modified (see Figure 1.6). We still start from relation (1.1.5) that
expresses that the solution of the advection equation (1.1.1) is constant along the
characteristics lines. The foot of the characteristic line that contains the
particular point $\,\, (x,\,t) \,\,$  in space-time is either the point
$\,\,(y=x-at,\,0) \,\,$ if   $\,\, x-at < L \,,\,$  either the point $\,\, \bigl(L,\,
\theta=t-{{1}\over{a}}(x-L)\bigr) \,\,$ if $\,\, x-at > L \,.\,$ In the first case,
we have $\,\, y>0 \,\,$ and $\,\, \theta <0 \,\,$ then the initial condition (1.3.2)
is  advected inside the domain (1.3.1) and we have~: 

\noindent  (1.3.5) $\qquad \displaystyle 
u(x, \,t) \,\,\, = \,\, u_0(x - at) \,\,, \qquad  \qquad   \,\,\,\,\,  \, x-at \,< \, L
\,. \,$ 

\smallskip \noindent   $\bullet \qquad \,\,\, $ 
On the contrary, if $\,\, x-at > L,\,$ we have $\,\, y>L \,\,$ and $\,\, \theta > 0
\,\,$ then the boundary condition at $\,\,x=L \,\,$ that takes now the expression 

\noindent  (1.3.6) $\qquad \displaystyle 
u(L, \,t) \,\,\, = \,\, w_{L}(t) \,,\qquad \quad \,\, t > 0 \,,\,$ 

\noindent 
is advected inside the domain of study and we have~:  

\noindent  (1.3.7) $\qquad  \displaystyle 
u(x, \,t) \,\,\, = \,\, w_{L} \Bigl( t \,+\, {{L}\over{a}}  \,-\, {{x}\over{a}}\Bigr)
\,\,, \qquad x-at \,>\,L  \,. \,$ 

\noindent
We have established the following 

\smallskip  \noindent 
{\bf Proposition 1.3.  $\quad$  Advection in the domain $\, 0<x<L \,,\,\, a < 0 .\,$ }

 \noindent 
Under the hypothesis (1.3.4), the resolution of the advection equation (1.2.1) in the
space-time domain (1.3.1) conducts to a {\bf well posed} problem when we introduce
the initial condition (1.3.2)  on the interval $\, ]0,\,L [ \,$ and the boundary
condition  (1.3.6) at the {\bf input} region located at $\, x \!=\! L  \,,\, $ without
any boundary condition at the output located at $\, x =0 \,.\,$ The solution of
Problem (1.2.1), (1.3.2) and  (1.3.6) is given by the relations (1.3.5) and (1.3.7).

\bigskip \bigskip 
 
\noindent  {\smcaps 2) $ \,\,\,\,\,\,\, $   Finite volumes for linear hyperbolic
systems.} 
\smallskip \noindent {\smcaps 2.1 } $ \,  $ { \bf  Linear advection. }

\noindent   $\bullet \qquad \,\,\,\,\, $ 
We still study the advection equation parameterized by some celerity $\, a > 0 \,$~: 

\noindent  (2.1.1) $\qquad \displaystyle
{{\partial W}\over{\partial t}} \,\,+\,\, {{\partial }\over{\partial x}} \bigl(
a \, W \bigr)  \,\,=\,\, 0 \,\,,\qquad t>0 \,,\qquad x \in \R \,, \,\,$ 

\noindent
and we search a discrete version of this mathematical model. For doing this, we
introduce a space step $\, \Delta x > 0 \,$ and a space grid composed by 
points $\,\,  x_{j} \,\,$  whose coordinates are multiples of this space step $\,
\Delta x ,\,$ {\it id est }

\noindent  (2.1.2) $\qquad \displaystyle
 x_{j} \,=\, j \, \Delta x \,, \qquad j \in \Z \,.  \,$ 

\noindent
For a finite domain, $\, ]0,\,L[ \,$ to fix the ideas, the above grid is limited to
integer values $\,j \,$ such that

\noindent  (2.1.3) $\qquad \displaystyle
0 \, \leq \, j \, \leq \, J \,=\, {{L}\over{\Delta x}} \,$

\noindent
and the vertices $\,\, (x_j)_{0\leq j \leq J} \,$ are usually used in the context
of the finite difference method. The intervals  $\,\, K_{j+1/2} \,=\,\,
 ]x_j,\,x_{j+1} [\,\,$ between two vertices can be considered as finite elements (or
finite volumes in our study) and they cover the entire domain $\, ]0,\,L[  \,: \,$ 

\noindent  (2.1.4) $\qquad \displaystyle
[ 0,\,L ] \,\,= \,\, \bigcup_{ 0 \leq j \leq J \! - \! 1} [x_{j} \,,\, x_{j+1} ] 
\, \,,\, $

\noindent 
as proposed in the general context of meshes (see {\it e.g.} Ciarlet [Ci78]). We introduce
also a time step $\,\Delta t >0 \,$ and the discrete time values at integer multiples
of the above quantum~: 

\noindent  (2.1.5) $\qquad \displaystyle
t^n \,=\, n \, \Delta t \,, \qquad n \in \N \,.  \,$ 

\noindent
We consider now a space-time volume $\,\,V_{j+1/2}^{n+1/2} \,\,$ obtained by
cartesian  product of the two intervals $\,\, ]x_{j} \,,\, x_{j+1}[ \,\,$ and $\,\,
]t^n ,\,t^{n+1} [\,\,$ (see Figure 2.1)~: 

\noindent  (2.1.6) $\qquad \displaystyle
V_{j+1/2}^{n+1/2} \,\,= \,\,  ]x_{j} \,,\, x_{j+1}[  \, \times \, ]t^n ,\,t^{n+1}
[ \,. \,\,$

\smallskip \noindent   $\bullet \qquad \,\,\,\,\, $ 
The finite volume scheme consists simply in integrating the advection equation
(2.1.1) inside the space-time domain $\,\, V_{j+1/2}^{n+1/2} \,\,$  introduced
previously~: 

\noindent  (2.1.7) $\qquad \displaystyle
\int_{ \displaystyle V_{j+1/2}^{n+1/2}} \, \biggl[ \, 
{{\partial W}\over{\partial t}} \,\,+\,\, {{\partial }\over{\partial x}} \bigl(
a \, W \bigr) \, \biggr] \, {\rm d}x \, {\rm d}t \,\, =\,\, 0 \,\,,\quad 0 \leq j
\leq J \,,\qquad n \geq 0 \,.  \, $

% version mac 
% \smallskip \vskip 2.cm  \smallskip 
% \qquad  \qquad  \qquad  \quad   \special{illustration fig.2.1.epsf scaled 500}  
% fin de la version version mac 

%%%%%%%%%%%%%%%%%%%%%%%%%%%%%%%%%       modif janvier 2011 
%  \titredroite={\pecaps   Finite volumes for linear hyperbolic systems} 

% version linux
\smallskip \centerline {  \epsfysize=3cm  \epsfbox  {fig.2.1.epsf} }   \smallskip
% fin de la version linux 

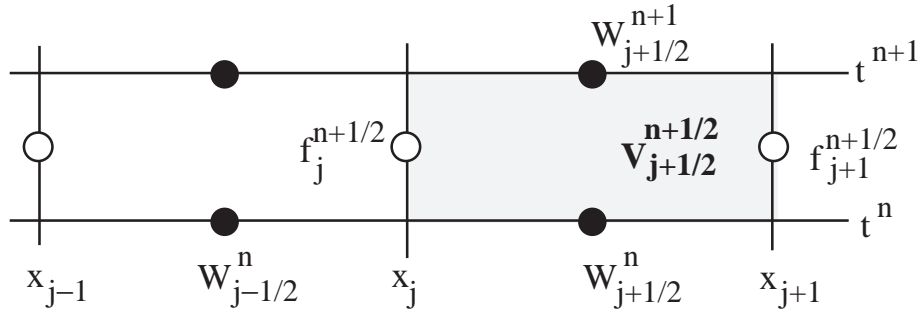
\centerline {\rm  {\bf Figure 2.1.}	\quad   Space-time grid for the finite volume
method.}
\bigskip 

\noindent  {\bf Proposition 2.1.  $\quad$  Finite volume scheme.}

 \noindent 
Let $\,\, \R \times [0,\,+\infty[ \, \ni (x,\,t) \longmapsto W(x,\,t) \in \R \,\,$ be
a solution of the advection equation (2.1.1). We introduce the space mean value
$\,\,W_{j+1/2}^{n} \,\,$   of this solution $\,\, W({\scriptstyle \bullet} ,\,
{\scriptstyle \bullet}) \,\,$ in the cell $\,\, K_{j+1/2} \, : \,$ 

\noindent  (2.1.8) $\qquad \displaystyle
W_{j+1/2}^{n} \,\,= \,\, {{1}\over{\mid  K_{j+1/2} \mid}} \,
\int_{\displaystyle x_{j}}^{\displaystyle x_{j\!+\!1}} \, W(x,\,t^n) \, {\rm d}x \,$ 

\noindent
and the time mean value $\,\, f_{j}^{n+1/2} \,\,$ of the so-called flux $\, a \, W(
{\scriptstyle \bullet} ,\, {\scriptstyle \bullet}) \,\,$   at the space position $\,
x_{j} \,$ and between discrete times $\, t^n \,$ and $\, t^{n+1} \,: \, $ 

\noindent  (2.1.9) $\qquad \displaystyle
f_{j}^{n+1/2} \,\,= \,\, {{1}\over{\Delta t}} \,
\int_ { \displaystyle  t^n }^{ \displaystyle  t^{n\!+\!1}} \, 
(a \,  W) (x_{j},\,t) \, {\rm d}t \,. \,$

\noindent 
Then we have the following constitutive relation of finite volumes schemes~: 

\noindent  (2.1.10) $\qquad \displaystyle
{{1}\over{\Delta t}} \bigl( W_{j+1/2}^{n+1} -  W_{j+1/2}^{n} \bigr) \,\,+\,\, 
{{1}\over{\Delta x}} \bigl( f_{j+1}^{n+1/2} - f_{j}^{n+1/2} \bigr) \,\, = \,\, 0
\,.\,$

\noindent 
This   numerical modelling  characterizes  the so-called finite volume method which
has been developed  thanks to the work of S. Godunov [Go59], Godunov et al [GZIKP79], 
Patankar [Pa80], Harten, Lax and Van Leer [HLV83] or  Faille,  Gallou\"et and  Herbin 
[FGH91] among others.

\smallskip \noindent   $\bullet \qquad \,\,\, $ 
The {\bf proof of Proposition 2.1} consists in a precise  evaluation of the left hand side of
equality (2.1.7). We use Fubini rule for the computation of double integrals and we
begin by integrating in time for the $\,\, {{\partial}\over{\partial t}} \,\,$
term~:  

\noindent  $ \displaystyle
\int_{ \displaystyle V_{j+1/2}^{n+1/2}} \, 
{{\partial W}\over{\partial t}} \, {\rm d}x \, {\rm d}t \,\, =\,\,
\int_{\displaystyle x_{j}}^{\displaystyle x_{j\!+\!1}} \,  \Bigl[ \, 
\int_ {\displaystyle t^n }^{\displaystyle t^{n\!+\!1}} \,
\, {{\partial W}\over{\partial t}}(x,\,t) \,  {\rm d}t  \Bigr] \,{\rm d}x \,\,$

\noindent  $ \displaystyle \qquad  \qquad 
= \,\, \int_{\displaystyle x_{j}}^{\displaystyle x_{j\!+\!1}} \,  
\Bigl[ \, W(x,\,t^{n+1}) -  W(x,\,t^{n}) \,
\Bigr] \, \,{\rm d}x \,\,= \,\, \Delta x   \Bigl[ \,W_{j+1/2}^{n+1} -  W_{j+1/2}^{n}
\, \Bigr]  \,\, $ 

\noindent 
due to the definition (2.1.8).  We proceed in an analogous way with the $\,\,
{{\partial}\over{\partial x}} \,\,$ term and begin now  the Fubini procedure by
integrating in space~; we have 

\noindent  $ \displaystyle
\int_{ \displaystyle V_{j+1/2}^{n+1/2}} \, 
{{\partial }\over{\partial x}} \bigl( a \, W \bigr)  \, {\rm d}x \, {\rm d}t \,\, =\,\,
 \int_ {\displaystyle t^n }^{\displaystyle t^{n\!+\!1}} \, 
\Bigl[ \,  \int_{\displaystyle x_{j}}^{\displaystyle x_{j\!+\!1}} \,  {{\partial
}\over{\partial x}}\,\bigl( a \, W \bigr) (x,\,t) \, {\rm d}x \, \Bigr] \, {\rm d}t
\,\,$

\noindent  $ \displaystyle \qquad \qquad 
= \,\,  \int_ {\displaystyle t^n }^{\displaystyle t^{n\!+\!1}} 
\, \Bigl[ \, ( a \, W \bigr) (x_{j+1},\,t) - ( a
\, W \bigr) (x_{j},\,t) \Bigr] \, {\rm d}t \,\,= \,\, \Delta t \,\, \Bigl[ \,
f_{j+1}^{n+1/2} - f_{j}^{n+1/2} \, \Bigr] \, \,$

\noindent  according to the definition   (2.1.9). We add the two previous results,
use identity (2.1.7) and  divide by $\,\, \Delta t \, \Delta x .\,\,$ We obtain
exactly the relation (2.1.10).  $ \hfill \square \kern0.1mm $

\smallskip \noindent   $\bullet \qquad \,\,\, $ 
The relation (2.1.10) is a very general form for the evolution of the mean values
$\,\, W_{j+1/2} \,\,$ between two time steps. The increment $\,\,  \bigl(
W_{j+1/2}^{n+1} -  W_{j+1/2}^{n} \bigr) \,\, $ is, after correction by
a multiplicative factor,  equilibrated by the flux difference
$\,\, \bigl( f_{j+1}^{n+1/2} - f_{j}^{n+1/2} \bigr) \,.\,$ The idea of a finite
volume scheme is to consider now that the algebraic object $\,\,  W_{j+1/2} \,\,$ is
nomore the mean value of the exact solution but an {\bf approximation} of this mean
value. Then the relation (2.1.10) proposes a numerical scheme for the discrete
evolution of the approximated mean values  $\,\,  W_{j+1/2} \,,\,\, j=0, \cdots,\,
J \!-\! 1 .\,$ Nevertheless, the numerical scheme is {\bf not} entirely defined by the
relation (2.1.10). Starting from mean values at the initial time step, {\it i.e.} 

\noindent  (2.1.11) $\qquad \displaystyle
W^0_{j+1/2} \,\, = \,\, {{1}\over{\Delta x}} \,  \int_{\displaystyle x_{j}}
^{\displaystyle x_{j\!+\!1}}
\,W_0(x) \, {\rm d} x \,\,,\qquad  j=0, \cdots,\, J \!-\! 1 \,, \,$

\noindent 
we are able to increment the time step with relation (2.1.10) only if all the fluxes
$\,\, f_j^{n+1/2} \,,$ $\,\,j=0,\cdots ,\, J \,$ have been {\it  a priori} first
determined as a functional of the previous values. In a very general way, we say that the
finite volume scheme (2.1.10) is an {\bf explicit  scheme} if each   flux  $\, f^{n+1/2}_j
\,\,$ is   a given function $\, \Psi_j \,$ of the mean values $\,\, \bigl( W^n_{k+1/2}
\bigr)_{k=1,\cdots,\, J\!-\!1} \, \,$  at the preceding time step number $\, n
\,$~: 

\noindent  (2.1.12) $\qquad \displaystyle
f^{n+1/2}_j \,\,= \,\,  \Psi_j \bigl(  \, \{ W^n_{k+1/2},\, k=0,\cdots,\, J\!-\!1 \}
\bigr)\,, \qquad j=0,\, \cdots ,\, J  \!-\! 1 \, . \,$

\noindent
The function  $\, \Psi_j \,$ is called the {\bf local numerical flux function} at point
$\, x_{j} \,$  and, joined with the evolution equation (2.1.10),  its choice 
determines   the numerical scheme. 

\smallskip \noindent   $\bullet \qquad \,\,\, $ 
A natural hypothesis claims that we have {\bf translation invariance} for the
evaluation of the flux  if we  move the discrete data in  the same  way~; in
other words, the numerical flux function $\,\, \Psi_j \, \,$ only depends on the $\,p
\,$ first neighbors of the interface $\, x_{j} \,.\,$ Then the explicit numerical
flux  is a given function $\, \Phi \,$ of the $\, p\, $ first  neighbors  and  we have~: 

\noindent  (2.1.13) $\qquad \displaystyle
f^{n+1/2}_j \,\,= \,\,  \Phi \bigl( W^n_{j+1/2-p}  ,\,\cdots ,\, 
W^n_{j-1/2} ,\, W^n_{j+1/2} ,\, \cdots ,\,  W^n_{j+1/2+p-1} \bigr) \,. \,$ 

\noindent
A very important particular case is  one of a two-point scheme for the evaluation of
the numerical flux. We have in this particular case~: 

\noindent  (2.1.14) $\qquad \displaystyle
f^{n+1/2}_j \,\,= \,\,  \Phi \bigl( W^n_{j-1/2} ,\, W^n_{j+1/2}  \bigr) \,. \,$ 

\noindent  
With this particular choice, the numerical scheme for  incrementing  in time of
the mean values takes the form~: 

\setbox21=\hbox {$\displaystyle {{1}\over{\Delta t}} \bigl( W_{j+1/2}^{n+1} - 
W_{j+1/2}^{n} \bigr)  \,\,+ \,\, $}
\setbox22=\hbox {$\displaystyle  \qquad  \,\,+ \,\,
{{1}\over{\Delta x}} \Bigl( \Phi \bigl( W^n_{j+1/2} ,\, W^n_{j+3/2}  \bigr) -  \Phi
\bigl(W^n_{j-1/2} ,\, W^n_{j+1/2} \bigr) \Bigr) \,\,= \,\, 0 \,. \,$}
\setbox30= \vbox {\halign{#&# \cr \box21 \cr \box22    \cr   }}
\setbox31= \hbox{ $\vcenter {\box30} $}
\setbox44=\hbox{\noindent  (2.1.15) $\displaystyle  \qquad \left\{ \box31 \right.$}  

\noindent $ \box44 $

\noindent
It is also a three-point finite difference scheme. The finite volume scheme
(2.1.10) (2.1.13) is said to be {\bf consistent} with the advection equation (2.1.1) 
when the numerical flux function $\, \Phi \,$ satisfies the condition 

\noindent  (2.1.16) $\qquad \displaystyle
\Phi \bigl( W  ,\,\cdots ,\,  W ,\, W ,\, \cdots ,\,  W \bigr) \,\, = \,\, a \, W
\,\,,\qquad \forall \, W \in \R \,.$ 

\smallskip \noindent   $\bullet \qquad \,\,\, $ 
The crucial question is how to choose a numerical finite volume scheme. The simplest
choice consists in a two point explicit scheme such that the finite difference scheme
is identical to the upstream-centered scheme (see {\it e.g.} Richtmyer-Morton [RM67]). It
takes the following expressions~: 

\noindent  (2.1.17) $\qquad \displaystyle
{{1}\over{\Delta t}} \bigl( W_{j+1/2}^{n+1} - 
W_{j+1/2}^{n} \bigr) \,\,+ \,\, a \Bigl( W_{j+1/2}^{n} - W_{j-1/2}^{n} \Bigr) \,\, =
\,\, 0 \,, \qquad a > 0 \,$

\noindent  (2.1.18) $\qquad \displaystyle
{{1}\over{\Delta t}} \bigl( W_{j+1/2}^{n+1} - 
W_{j+1/2}^{n} \bigr) \,\,+ \,\, a \Bigl( W_{j+3/2}^{n} - W_{j+1/2}^{n} \Bigr) \,\, =
\,\, 0 \,, \qquad a < 0\,. \,$

\smallskip \noindent 
The corresponding flux function is called the {\bf first order upstream-centered
flux},  is simply given by the following relations~: 

\setbox21=\hbox {$\displaystyle a \, W_l \,,\qquad a > 0 \, $}
\setbox22=\hbox {$\displaystyle a \, W_r \,,\qquad a < 0 \,.\, $}
\setbox30= \vbox {\halign{#&# \cr \box21 \cr \box22    \cr   }}
\setbox31= \hbox{ $\vcenter {\box30} $}
\setbox44=\hbox{\noindent  (2.1.19) $\displaystyle  \qquad  \Phi(W_l,\, W_r) \,\,=
\,\,  \left\{ \box31 \right.$}  

\noindent $ \box44 $

\noindent
When this flux function acts at a given point $\, x_{j} \,$ of the mesh, we have~: 

\setbox21=\hbox {$\displaystyle a \, W_{j-1/2}^{n} \,,\qquad a > 0 \, $}
\setbox22=\hbox {$\displaystyle a \, W_{j+1/2}^{n} \,,\qquad a < 0 \,.\, $}
\setbox30= \vbox {\halign{#&# \cr \box21 \cr \box22    \cr   }}
\setbox31= \hbox{ $\vcenter {\box30} $}
\setbox44=\hbox{\noindent  (2.1.20) $\displaystyle  \qquad  f_j^{n+1/2}
\,\,= \,\, \Phi \bigl(  W_{j-1/2}^{n} ,\,  W_{j+1/2}^{n} \bigr) \,\,= \,\,   \left\{
\box31 \right.$}  

\noindent $ \box44 $

\noindent
If $\, a>0 , \,$ the exact solution of the advection equation propagates  the
information from the left to the right~; the flux at the interface $\, x_{j} \,$ is
issued from the cell at the left of the interface and this cell at the number $\,
j \!-\! 1/2 \,.$ If $\, a < 0 ,\,$ the propagation of the information with the
advection equation is from right to left~; the interface flux at the abscissa  $\,
x_{j} \,$ is due to the control volume on the right, {\it i.e.} with number $\, j\!+\!1/2
\,$ as depicted on Figure 2.2. 

\bigskip 
\smallskip \centerline {  \epsfysize=3,5cm  \epsfbox  {fig.2.2.epsf} }  \smallskip

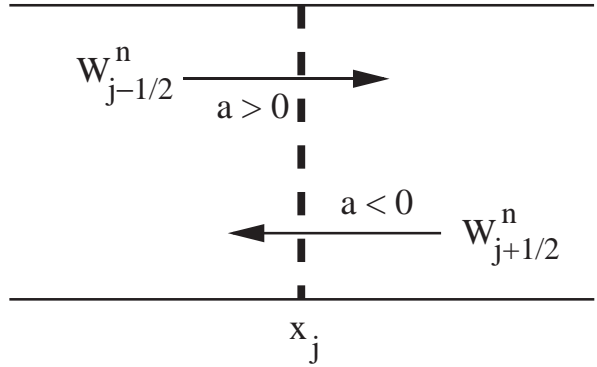
\centerline {\rm  {\bf Figure 2.2.}	\quad   Upwinding of the information for the
advection equation. }
\bigskip 

\smallskip \centerline {  \epsfysize=2,5cm  \epsfbox  {fig.2.3.epsf} }  \smallskip

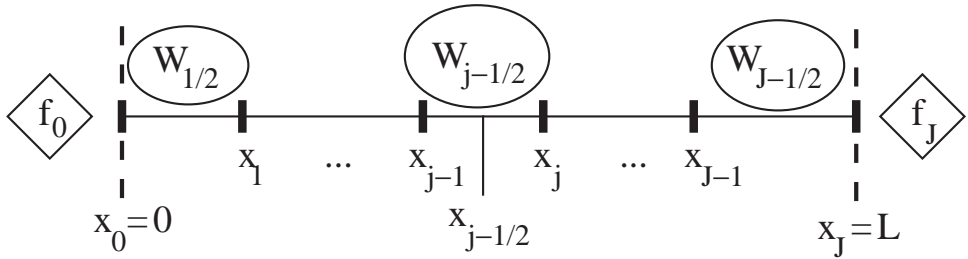
\centerline {\rm  {\bf Figure 2.3.}	\quad   Notations for the one-dimensional finite
volume method. }

\bigskip  
\vfill \eject    %%%%%%%%%%%%%%%%%%%%%%%%%%%%%%  modif janvier 2011 
\noindent   $\bullet \qquad \,\,\, $ 
Recall that practical use of the upwind finite volume scheme like (2.1.17) when $\, a >
0 \,$ or (2.1.18)  if $\, a < 0 \,$  is restricted to the usual
Courant-Friedrichs-Lewy stability condition~: 

\noindent  (2.1.21) $\qquad \displaystyle
a \, {{\Delta t}\over{\Delta x}} \,\, \leq \,\, 1 \,$

\noindent
as developed {\it e.g.} in the book of Richtmyer and Morton [RM67]. 

\bigskip \noindent {\smcaps 2.2 } $ \, $ { \bf Numerical flux boundary conditions }

\noindent   $\bullet \qquad \,\,\,\,\, $ 
In this section, we focus on the problem of the numerical  boundary conditions.
Recall that we study the advection equation in the space domain $\,\, [0,\,L] \,: \,$ 

\noindent  (2.2.1) $\qquad \displaystyle
0 \, \leq \, x \, \leq \, L \, $

\noindent 
and $\,\, J \,= \, {{L}\over{\Delta x}} \, \in \N \,\,$ control cells (or finite
elements) have been used to define a mesh~: 

\noindent  (2.2.2) $\qquad \displaystyle
J \, \Delta x \,\, = \,\, L \,.\,  $

\noindent 
Note that the $\,j^{\rm th} \,$  cell is exactly the interval $\,\, ]x_{j-1} ,\,
x_{j}[ \,\,$ and it is  centered at point $\,\, x_{j-1/2}\,$ as shown on Figure
2.3.

\smallskip \noindent   $\bullet \qquad \,\,\,\,\, $ 
At time step $\, n \, \Delta t \,,\,$ the discrete field  is entirely known and is
composed of all the values $\,\, W^n_{j-1/2} \,$ for $\,\, j\!=\! 1,\, \cdots ,\, J \,.
\, $ With a flux function $\,\, \Phi({\scriptstyle \bullet} ,\, {\scriptstyle \bullet}) 
\,\,$ as proposed at relation (2.1.14), we observe that the two boundary fluxes $\,\,
f^{n+1/2}_0 \,$ and $\, f^{n+1/2}_J \,$  are {\bf not} a priori defined because states
$\, W^n_{-1/2} \,$ or $\, W^n_{J+1/2} \,$ does not exist. The situation is more complex
with numerical fluxes that use four points or more as proposed in  (2.1.13) and   will
not be detailed in this section. Even if the formula giving the numerical flux at the
boundaries has to be specifically studied, the finite volume scheme remains defined by
the relation (2.1.10) and we have for the two cells encountering the boundary~: 

\noindent  (2.2.3) $\qquad \displaystyle
{{1}\over{\Delta t}} \bigl( W_{1/2}^{n+1} -  W_{1/2}^{n} \bigr) \,\,+\,\, 
{{1}\over{\Delta x}} \bigl( f_{1}^{n+1/2} - f_{0}^{n+1/2} \bigr) \,\qquad  = \,\, 0
\,,\,$ 

\noindent  (2.2.4) $\qquad \displaystyle
{{1}\over{\Delta t}} \bigl( W_{J-1/2}^{n+1} -  W_{J-1/2}^{n} \bigr) \,\,+\,\, 
{{1}\over{\Delta x}} \bigl( f_{J}^{n+1/2} - f_{J-1}^{n+1/2} \bigr) \,\, = \,\, 0
\,.\,$ 

\smallskip \noindent   $\bullet \qquad \,\,\,\,\, $ 
The question is now to adapt the relation (1.2.14) in  order to determine the two
{\bf boundary fluxes} $\,\, f_{0}^{n+1/2} \,\,$ at the left of the domain and  $\,\,
f_{J}^{n+1/2} \,\,$ at the right. For the advection equation with celerity $\, a>0
,\,$ we have observed in the first section  that some boundary condition $\, v_0(t)
\,$ has to be assigned 
 at $\, x=0 \,$ and it is not the case for $\, x=L .\,$ It is
therefore natural to take into account this information at the input of the domain
and to set~: 

\noindent  (2.2.5) $\qquad \displaystyle
f_{0}^{n+1/2} \,\, =\,\, {{1}\over{\Delta t}} \, 
\int_{\displaystyle t^n}^{\displaystyle t^{n\!+\!1}}  \!\!\!\! 
\!\! a \, v_0(t) \, {\rm d} t \,$ 

\noindent
or simply 

\noindent  (2.2.6) $\qquad 
f_{0}^{n+1/2} \,\, =\,\,  \, a \, v_0\bigl( (n+{{1}\over{2}})  \Delta t \bigr)
\,,\qquad a > 0 \,, \,$

\noindent 
if function $\, \, t \longmapsto  v_0(t) \,  \,$ has a slow time variation at the scale
defined by the time step. At the output $\, x=L ,\,$ no numerical datum has to be
assigned   to set correctly the continuous mathematical problem. We must maintain this
property if we wish the   numerical method to follow the mathematical
physics  as efficiently    as possible.
 A simple boundary flux is associated with the previous numerical upwind scheme.
For $\, \,x=x_{J}= L \,\,$ and $\, a>0 ,\,$ we observe that the upwind scheme (2.1.20)
is simply written as~: 

\noindent  (2.2.7) $\qquad 
f_{J}^{n+1/2} \,\, =\,\,  \, a \, W^n_{J-1/2} \, \,,\qquad a>0 \,,\,$

\noindent
and this relation (2.2.7) defines a {\bf first order extrapolated} boundary flux.

\smallskip \noindent   $\bullet \qquad \,\,\,\,\, $ 
The roles are reversed when $\, a<0 .\,$ The abscissa $\, x=0 \,$ corresponds to an
output for the advection equation and the right boundary $\, x=L \,$ is an input where
a time field  $\, \, t \longmapsto  w_{L}(t) \,\,$ is given. In the first case, the
upwind scheme (2.1.20) can be applied without modification~: 

\noindent  (2.2.8) $\qquad 
f_{0}^{n+1/2} \,\, =\,\,  \, a \,  W^n_{1/2} \, \,,\qquad a < 0 \,,\,$

\noindent
and it corresponds to a first order extrapolation of the  internal data $\,\, \bigl\{
W^n_{j-1/2} ,\, $ $ \,j = 1,\, \cdots ,\, J \, \bigr\} \,\,$  at the boundary at time step
$\, n \Delta t .\,$ For $\, x=L ,\,$ the boundary flux $\, f_J^{n+1/2} \,$ uses the
given information between the two time steps~: 

\noindent  (2.2.9) $\qquad 
f_{J}^{n+1/2} \,\, =\,\,  \, a \, w_{L}\bigl( (n+{{1}\over{2}})  \Delta t \bigr)
 \, \,,\qquad a < 0 \,. \,$

\smallskip   \noindent  {\bf Proposition 2.2.}

 \noindent $\hfill$  {\bf    Flux boundary conditions for
the advection equation. }
 
\noindent 
When we approach the advection equation (2.1.1) with the finite volume method, the
numerical boundary conditions induces a choice for the two boundary fluxes $\,\,
f_{0}^{n+1/2} \,\,$ and  $\,\, f_{J}^{n+1/2} . \,$  When $\,\,a>0 ,\,$ the boundary
condition $\,\, v_0(t) \,\,$ at the input  can be  introduced  into the boundary
 with  the relation (2.2.6)  and the free output

% 
%  \noindent 

\noindent 
at the right can be treated with an extrapolation of the type  (2.2.7).  When  $\,\,a < 0
,\,$ the free output at the left of the domain can be taken into account with the help of
relation (2.2.8) whereas the input condition $\,\, w_{L}(t) \,\,$ at the right can be
introduced thanks to relation (2.2.9).

\bigskip \noindent {\smcaps 2.3 } $ \, $ { \bf A model system with two equations }

\noindent   $\bullet \qquad \,\,\,\,\, $ 
Let $\,\,  a > 0 \, \,$ and $\,\, b>0 \,\,$ be two positive real number. We study in
this section a model problem that is composed by the juxtaposition of an advection
equation with celerity $\,a \,$ and an advection with celerity $\,-b .\,$ We explicit
the associated algebra~: 

\noindent  (2.3.1) $\qquad \displaystyle
{{\partial u}\over{\partial t}} \,\,+\,\,a \,  {{\partial u}\over{\partial x}}  
\,\,=\,\, 0 \,\,,\qquad t>0 \,,\qquad x \in \R \,, \,\,$ 

\noindent  (2.3.2) $\qquad \displaystyle
{{\partial v}\over{\partial t}} \,\,-\,\,b \,  {{\partial v}\over{\partial x}}  
\,\,=\,\, 0 \,\,,\qquad t>0 \,,\qquad x \in \R \,. \,\,$ 

\noindent
We associate the two equations (2.3.1) and (2.3.2) and consider a unique problem
with a vector field as  unknown. We set~: 

\noindent  (2.3.3) $\qquad \displaystyle
\varphi \,\, = \,\, \pmatrix{u \cr v \cr} \,$

\noindent
and the set of equations  (2.3.1)-(2.3.2) can naturally be written as a system~: 

\noindent  (2.3.4) $\qquad \displaystyle
{{\partial \varphi}\over{\partial t}} \,\,+\,\, \pmatrix{ a & ~0 \cr 0 & -b  \cr}
\,{{\partial \varphi}\over{\partial x}} \,\,=\,\, 0 \,.\,$ 

\noindent
By introducing the {\bf flux function} $\,\,F(\varphi) \,\,$ according to the relation 

\noindent  (2.3.5) $\qquad \displaystyle
F(\varphi) \,\, = \,\, \pmatrix{a \, u \cr -b \, v \cr} \, \,$

\noindent
the system (2.3.4) takes the general conservative form~: 

\noindent  (2.3.6) $\qquad \displaystyle
{{\partial \varphi}\over{\partial t}} \,\,+\,\, {{\partial}\over{\partial x}} \bigl(
F(\varphi) \bigr) \,\,= \,\, 0 \,. \,$ 

\smallskip \noindent  $\bullet \qquad \,\,\,\,\, $ 
The approximation of system (2.3.6) with a grid parameterized by a space step $\,
\Delta x \,$ and a time step $\, \Delta t \,$ is conducted exactly as in the case of
the advection equation. The following property is a straightforward generalization of
Proposition~2.1. We left the proof to the reader.

\bigskip  \noindent  {\bf Proposition 2.3.  $\quad$  Finite volume scheme.}

 \noindent 
Let $\,\, \R \times [0,\,+\infty[ \, \, \ni (x,\,t) \longmapsto \varphi(x,\,t) \in \R
 \times \R \,\,$ be a solution of the linear conservation law   (2.3.6). We define
the   space mean value $\,\,\varphi_{j+1/2}^{n} \,\,$   of this solution $\,\,
\varphi ({\scriptstyle \bullet} ,\, {\scriptstyle \bullet}) \,\,$ in the cell 
$\,\, K_{j+1/2} \, : \,$ 

% 
%  \noindent 

\noindent  (2.3.7) $\qquad \displaystyle
\varphi_{j+1/2}^{n} \,\,= \,\, {{1}\over{\mid  K_{j+1/2} \mid}} \,
\int_{\displaystyle x_{j}}^{\displaystyle x_{j\!+\!1}} \, \varphi(x,\,t^n) \, {\rm d}x \,$ 

% 
%  \noindent 

\noindent
and the time mean value $\,\, f_{j}^{n+1/2} \,\,$ of the  flux function introduced in
(2.3.5)  at the space position $\, x_{j} \,$  between discrete times $\, t^n \,$
and $\, t^{n+1} \,: \, $ 

\noindent  (2.3.8) $\qquad \displaystyle
f_{j}^{n+1/2} \,\,= \,\, {{1}\over{\Delta t}} \, 
\int_ {\displaystyle t^n }^{\displaystyle t^{n\!+\!1}} \,F \bigl(
\varphi (x_{j},\,t ) \bigr)  \, {\rm d}t \,. \,$ 

\noindent 
We have the following   relation that characterizes the  finite volumes
schemes~: 

\noindent  (2.3.9) $\qquad \displaystyle
{{1}\over{\Delta t}} \bigl( \varphi_{j+1/2}^{n+1} -  \varphi_{j+1/2}^{n} \bigr)
\,\,+\,\,  {{1}\over{\Delta x}} \bigl( f_{j+1}^{n+1/2} - f_{j}^{n+1/2} \bigr) \,\, =
\,\, 0  \,.\,$

\smallskip \noindent  $\bullet \qquad \,\,\,\,\, $ 
We have now to propose a precise numerical flux function analogous to the relation
(2.1.12) to transform the conservation property (2.3.9) into a finite volume 
numerical scheme able to propagate the discrete values $\,\, \varphi_{j+1/2}^{n} \,\,$
up to the discrete time $\, t^{n+1} .\,$ For internal interfaces $\,\, x_{j} ,\,
j=1,\cdots \,,J\!-\!1 \,,\,$ it is natural to apply the upwinding scheme (2.1.20)
with a left upwinding for the first equation and  a right upwinding for the equation
(2.3.2).  Figure 2.4 illustrates the associated  algebra~: 

\setbox21=\hbox {$\displaystyle a \, u_{j-1/2}^{n}  $}
\setbox22=\hbox {$\displaystyle -b \, v_{j+1/2}^{n}  $}
\setbox30= \vbox {\halign{#&# \cr \box21 \cr \box22    \cr   }}
\setbox31= \hbox{ $\vcenter {\box30} $}
\setbox44=\hbox{\noindent  (2.3.10) $\displaystyle  \qquad  f_j^{n+1/2}
\, = \, \Phi \bigl(  \varphi_{j-1/2}^{n} ,\,   \varphi_{j+1/2}^{n} \bigr) \, =
\,    \left( \box31 \right) \,,\quad  j=1,\cdots \,,J\!-\!1 \,.\,$  }
\noindent $ \box44 $

\smallskip 
\smallskip \centerline {  \epsfysize=3,5cm  \epsfbox  {fig.2.4.epsf} } \smallskip 

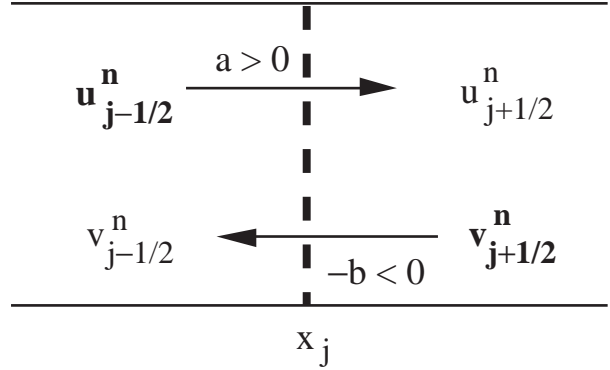
\centerline {\rm  {\bf Figure 2.4.}	\quad   Interface upwind numerical flux  }

\centerline {\rm  for a model problem with two equations. }

\smallskip \centerline {  \epsfysize=4,1cm  \epsfbox  {fig.2.5.epsf} }  \smallskip

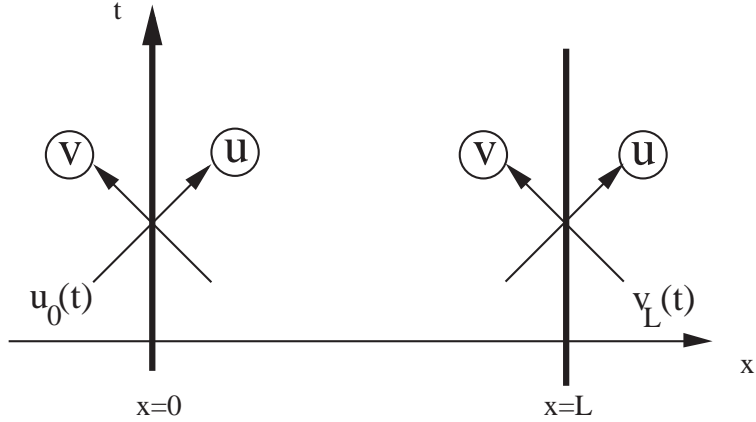
\centerline {\rm  {\bf Figure 2.5.}	\quad   Boundary conditions for a model problem 
with two equations. }

\bigskip 

\noindent  $\bullet \qquad \,\,\,\,\, $ 
At the left boundary $\,\, x=0 \,,\,$ we have an input for the variable $\, u \,$ and
we suppose given the associated datum $\, [0,\, +\infty [ \, \ni t \longmapsto  u_0(t)
\in \R \,$~: 

\noindent  (2.3.11) $\qquad \displaystyle
u(0,\,t) \,\,= \,\, u_0(t) \,\,, \qquad t>0 \,$

\noindent
whereas it is an output for the $\, v \,$ variable. By association of relations
(2.2.6) and (2.2.8), we obtain 

\noindent  (2.3.12) $\qquad \displaystyle
f_0^{n+1/2} \,\,= \,\, \pmatrix {a \, u_0\bigl( (n+{{1}\over{2}}) \Delta t \bigr)
\cr -b \, v_{1/2}^n \cr } \,\,. \,$ 

\noindent 
At the other boundary of the interval$\,\, ]0,\,L[ \,,\,$ we have an output for the
first variable $\, u \,$ and an input for the second one, and an associated boundary
condition  $\,\,  [0,\, +\infty [ \, \ni t \longmapsto   v_{L}(t) \in \R \, \,$ is
supposed to have been given~: 

\noindent  (2.3.13) $\qquad \displaystyle
v(L,\,t) \,\,= \,\, v_{L}(t) \,\,, \qquad t>0 \,$

\noindent
as illustrated on Figure 2.5. The numerical flux at the right is evaluated by
association of the relations (2.2.7) and (2.2.9)~: 

\noindent  (2.3.14) $\qquad \displaystyle
f_L^{n+1/2} \,\,= \,\, \pmatrix {a \, u_{J-1/2}^n \cr -b \, v_{L}\bigl(
(n+{{1}\over{2}}) \Delta t \bigr) \cr } \,. \,$ 

\bigskip \noindent {\smcaps 2.4 } $ \, $ { \bf Unidimensional linear acoustics} 

\noindent   $\bullet \qquad \,\,\,\,\, $ 
We consider a gas in a pipe of uniform section at normal 
 conditions of temperature and
pressure. The reference density  is denoted by $\,\, \rho_0 \,\,$ and the
reference pressure is named $\,\, p_0 .\,$ The {\bf sound celerity} $\, c_0 \,$ of
this gas satisfies the relation 

\noindent  (2.4.1) $\qquad \displaystyle
c_0 \,\,= \,\, \sqrt{ {{\gamma p_0}\over{\rho_0}}} \,\,$

\noindent 
with $\,\, \gamma \!=\! 1.4 \, \,$ as proved {\it e.g.} in the book of Landau and Lifchitz
[LL54]. A sound wave is a small perturbation of this reference state. The
differences of density, pressure and velocity  fields are denoted respectively 
by $\, \rho , \,p \,$ and $\, u .\,$ The hypothesis of a small perturbation implies
that the entropy of the reference state is maintained for all the time evolution and
in consequence, it is easy to establish the following relation between the
perturbations of density and pressure~: 

\noindent  (2.4.2) $\qquad \displaystyle
p \,\,= \,\, c_0^2 \, \rho \,. \,$

\smallskip \noindent   $\bullet \qquad \,\,\,\,\, $ 
The conservation of mass leads  to a first order linear conservation law~: 

\noindent  (2.4.3) $\qquad \displaystyle
{{\partial \rho}\over{\partial t}} \,\, + \,\, \rho_0 \, 
{{\partial u}\over{\partial x}} \,\, = \,\, 0 \,\,$ 

\noindent
and the conservation of   momentum links the time evolution of velocity with the
spatial gradient of pressure~: 

\noindent  (2.4.4) $\qquad \displaystyle
\rho_0 \, {{\partial u}\over{\partial t}} \,\, + \,\, 
{{\partial p}\over{\partial x}} \,\, = \,\, 0 \,. \,$ 

\noindent
We introduce the vector $\,\, W = \pmatrix{p \cr u \,} \,$ of unknowns.
Then the equations (2.4.3) and (2.4.4) can be written as a linear hyperbolic system
of conservation laws~: 

\noindent  (2.4.5) $\qquad \displaystyle
{{\partial W}\over{\partial t}} \,\, + \,\, A \, 
{{\partial W}\over{\partial x}} \,\, = \,\, 0 \, \,$ 

\noindent
with 

\noindent  (2.4.6) $\qquad \displaystyle
A \,\,= \,\, \pmatrix {0 & \rho_0 \, c_0^2 \cr  { {{1}\over{\rho_0}}} & 0 \cr} \,.
\,$ 

\smallskip \noindent   $\bullet \qquad \,\,\,\,\, $ 
When we consider the eigenvalues and eigenvectors of matrix $\, A ,\,$  it is  natural
to introduce the {\bf  characteristic variables}  defined respectively by 

\noindent  (2.4.7) $\qquad \displaystyle
\varphi_+ \,\,= \,\, p \,+\, \rho_0 \, c_0 \, u \,\,$

\noindent  (2.4.8) $\qquad \displaystyle
\varphi_- \,\,= \,\, p \,-\, \rho_0 \, c_0 \, u \,\,$

\noindent 
and the quantity $\,\,  \rho_0 \, c_0 \,\,$ is named the {\bf acoustic impedance}. We
have from the relations (2.4.3) and (2.4.4)~: 

\noindent $ \displaystyle
{{\partial \varphi_+}\over{\partial t}} \, + \, c_0 
{{\partial \varphi_+}\over{\partial x}} \,\,\, = \,\, \, \Bigl( 
{{\partial p}\over{\partial t}}\,+\,   \rho_0 \, c_0 \,{{\partial u}\over{\partial t}}
\Bigr) \,\,+\,\, \Bigl(   c_0 \, {{\partial p}\over{\partial x}} \,+\,  \rho_0 \,
c_0^2 \, \,{{\partial u}\over{\partial x }} \Bigr) \,\,$ 

\noindent $ \displaystyle \qquad \qquad \qquad  \quad \,\,\,
= \,\,  c_0^2 \, \, \Bigl( {{\partial \rho}\over{\partial t}}\,+\, \rho_0 \,
 \,{{\partial u}\over{\partial x }} \Bigr) \,\,+\,\,   c_0 \, \Bigl( \rho_0
\,  {{\partial u}\over{\partial t}} \,+\,  {{\partial p}\over{\partial x}}
\Bigr) \,\,\,\, = \,\, 0 \,, \, $ 

\noindent $ \displaystyle
{{\partial \varphi_-}\over{\partial t}} \, - \, c_0 
{{\partial \varphi_-}\over{\partial x}} \,\,\, = \,\, \, \Bigl( 
{{\partial p}\over{\partial t}}\,-\,   \rho_0 \, c_0 \,{{\partial u}\over{\partial t}}
\Bigr) \,\,-\,\,  c_0 \, \Bigl(  {{\partial p}\over{\partial x}} \,-\,  \rho_0
\, c_0 \, \,{{\partial u}\over{\partial x }} \Bigr) \,\,$ 

\noindent $ \displaystyle \qquad \qquad \qquad  \quad \,\,\,
= \,\,  c_0^2 \, \, \Bigl( {{\partial \rho}\over{\partial t}}\,+\, \rho_0 \,
 \,{{\partial u}\over{\partial x }} \Bigr) \,\,-\,\,   c_0 \, \Bigl( \rho_0
\,  {{\partial u}\over{\partial t}} \,+\,  {{\partial p}\over{\partial x}}
\Bigr) \,\,\,\, = \,\, 0 \,, \, $ 

\noindent
and we recover a system of the type (2.3.4) studied previously~: 

\noindent  (2.4.9) $\qquad \displaystyle
{{\partial}\over{\partial t}}  \pmatrix{\varphi_- \cr \varphi_+ \cr}  \,\, +\,\,
\pmatrix{-c_0 & 0 \cr 0 & c_0} \,  {{\partial}\over{\partial t}}   \pmatrix{\varphi_-
\cr \varphi_+ \cr}  \,\,= \,\, 0 \,. \,$ 

\smallskip \noindent   $\bullet \qquad \,\,\,\,\, $ 
A typically   physical problem is the following~: a given acoustic  pressure wave $\,\,
[0,\,+\infty[ \,\, \ni t \longmapsto \Pi(t) \, > 0 \,\, $ is injected at the left
$\, x=0 \,$ of the pipe and the waves go away freely at the right boundary $\, x
\!=\! L \,.\,$ At $\, t=0 ,\,$ the velocity and pressure of the fluid are given~: 

\noindent  (2.4.10) $\qquad \displaystyle
u(x,\,0) \,\, = \,\, u_0(x) \,\,,\qquad 0 < x < L \,$

\noindent  (2.4.11) $\qquad \displaystyle
p(x,\,0) \,\, = \,\, p_0(x) \,\,,\qquad 0 < x < L  \, . \,$

\noindent
From a mathematical viewpoint,  the boundary conditions have to respect the dynamics of
this system of acoustic equations written in   diagonal form (2.4.9)~: the
variable $\, \varphi_+ \,$ must be given at $\, x \!=\! 0 \,$ and the variable $\, 
\varphi_- \,$ at the abscissa $\, x \!=\! L .\,$ From (2.4.7) and (2.4.8), we
determine the pressure as a function of the two characteristics variables $ \, 
\varphi_+  \,$ and  $ \,  \varphi_-  \, :
\,$ 

\noindent  (2.4.12) $\qquad \displaystyle
p \,\, = \,\, {{1}\over{2}} \, \bigl(   \varphi_+ \,+\,   \varphi_- \bigr) \,$ 

\noindent 
and if the pressure is imposed at $\, x=0 ,\,$ the relation (2.4.12) can be written
under the form~: 

\noindent  (2.4.13) $\qquad \displaystyle
\varphi_+(0,\,t) \,\,= \,\, - \, \varphi_-(0,\,t)  \,+\, 2 \, \Pi(t) \,  \,\,,\qquad
x=0 \,\,, \quad t> 0 \,,\,$ 

\noindent
that makes in evidence a reflection operator~:  the input variable $\,  \varphi_+ \,$
is a given affine function of the output variable $\,  \varphi_-  \, .\,$ At the other
boundary $\, x \!=\! L \,,\,$ the notion of free output expresses that the waves that go
outside of the domain of study have no reflection at the boundary. When $\, x \!=\! L
\,,\,$ the characteristic variable  $\,   \varphi_+ \,$ is going outside and there is
no boundary condition for this variable. We have to express also that this wave has no
influence on   the characteristic $\,  \varphi_-  \,$  that wish to go inside the domain
$\, ]0,\,L[ .\,$  In other terms, the input value  $\,  \varphi_-  \,$ is independent of
the variable  $\,  \varphi_+  \,$ and also of  time. We have in consequence 

\noindent  (2.4.14) $\qquad \displaystyle
{{\partial}\over{\partial t}} \varphi_-(L,\,t) \,\, = \,\, 0 \, . \,$

\smallskip \noindent 
We have established 

% \bigskip  
\vfill  \eject               %%%%%%%%%%%%  modif janvier 2011 
\noindent  {\bf Proposition 2.4.  $\quad$  Boundary conditions for
acoustic problem. }

 \noindent 
The mathematical boundary conditions associated with the datum of a given acoustic 
pressure wave $\,\, [0,\,+\infty[ \,\, \ni t \longmapsto \Pi(t) \, > 0 \,\, $ at the
left of the domain $\, ]0,\, L [ \,$ admits the expression  (2.4.13)  and a condition
of  free output of the waves at the   right boundary $\, x \!=\! L \,$  can be
expressed by the relation (2.4.14).

\smallskip \centerline {  \epsfysize=4,3cm  \epsfbox  {fig.2.6.epsf} }  \smallskip
% fin de la version linux  

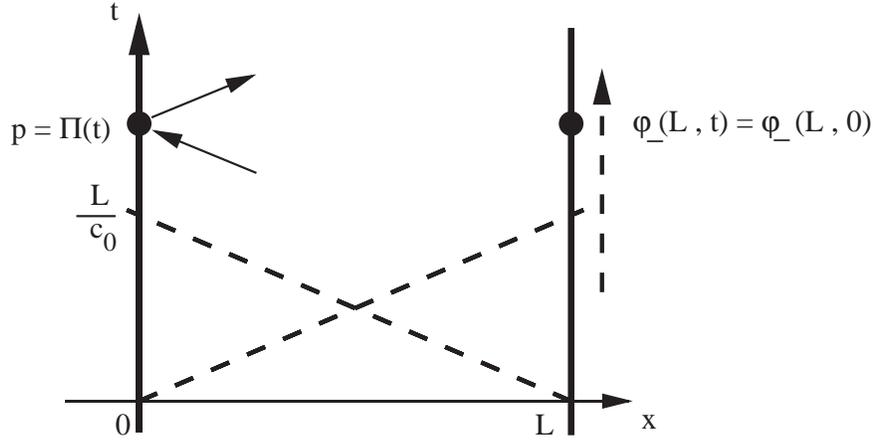
\centerline {\rm  {\bf Figure 2.6.}	\quad   Solution of the acoustic equations 
in one space dimension }

\centerline {\rm  for a model problem with two equations}
\bigskip

\noindent   $\bullet \qquad \,\,\, $ 
The above acoustic problem associated with the first order partial differential
equations (2.4.3) (2.4.4), the initial conditions (2.4.10) (2.4.11) and the boundary
conditions (2.4.13) (2.4.14) is illustrated on Figure 2.6. The initial conditions are
active in the beginning of the evolution in time $\,\, (t \leq {{L}\over{c_0}} )\,\,$
and have a trace for higher times due to the boundary conditon (2.4.13), that gives,
due to (2.4.8) and (2.4.13)~: 

\noindent  (2.4.15) $\quad \displaystyle
\varphi_-(x,\,t) \,\,\equiv \,\, p(x,\,t) \,-\, \rho_0 \, c_0 \, u(x,\,t)   \,\, =
\,\,  p_0(L) \,-\, \rho_0 \, c_0 \, u_0(L)  \, , \quad t \geq{{L}\over{c_0}} \,. \,$

\noindent
On the other hand, the inflow boundary condition (2.4.12) and the second row of
matrix equation (2.4.9) implies~: 

\setbox21=\hbox {$\displaystyle 
\varphi_+(x,\,t) \,\equiv \,p(x,\,t) \,+\, \rho_0 \, c_0 \, u(x,\,t)   \,\, = \,\, $}
\setbox22=\hbox {$\displaystyle  \qquad   \qquad   \qquad  
\,\, = \,\,  2 \, \Pi \bigl( t-{{x}\over{c_0}} \bigr) \,-\, \varphi_- \bigl( 0,\, 
t-{{x}\over{c_0}} \bigr)  \, , \quad t \geq{{L}\over{c_0}} \,. \,\,$}
\setbox30= \vbox {\halign{#&# \cr \box21 \cr \box22    \cr   }}
\setbox31= \hbox{ $\vcenter {\box30} $}
\setbox44=\hbox{\noindent  (2.4.16) $\displaystyle  \qquad \left\{ \box31 \right.$}  
\noindent $ \box44 $

\noindent
We deduce from the relations (2.4.15) (2.4.16) joined with the definitions (2.4.7) and 
(2.4.8)~: 

\noindent  (2.4.17) $\quad 
p(x,\,t) \,\,= \,\,  \Pi \bigl( t-{{x}\over{c_0}} \bigr)  \,, \hfill  \qquad
0 \leq x \leq L \,,\quad  t \geq{{L}\over{c_0}} \, \,$

\noindent  (2.4.18) $\quad  
u(x,\,t) \,\,= \,\, u_0(L)  \,+\, {{1}\over{ \rho_0 \, c_0 }}  \Bigl(  \Pi \bigl(
t-{{x}\over{c_0}} \bigr)    -    p_0(L)  \Bigr) \,, \quad 0 \leq x \leq
L \,,\quad  t \geq{{L}\over{c_0}} \,. \,$

\smallskip \noindent   $\bullet \qquad \,\,\, $ 
We turn now to the numerical finite volume scheme. We have to determine the internal
fluxes $\, f_j^{n+1/2} \,,$ $\,\,j=1,\cdots ,\, J \!-\! 1 \,$ and the boundary fluxes 
$\, f_0^{n+1/2} \,$ and $\, f_J^{n+1/2} .\, $ Recall first that the physical flux
$\,\, F(W) \,\,$  function for the acoustic equation (2.4.5) is equal to 

\noindent  (2.4.19) $\qquad \displaystyle 
F(W) \,\,= \,\, \pmatrix{\rho_0 \, c_0^2 \, u \cr {{1}\over{\rho_0}} \, p} \quad \,$
with $\quad W \,\,= \,\,  \pmatrix{p \cr u} \,. \,$ 

\bigskip  \noindent  {\bf Proposition 2.5.  $\quad$  Upwind scheme for
computational acoustics.  }

\noindent 
The extension of the upwind finite volume scheme (2.3.10), (2.3.12) and (2.3.14) is
determined by the following relations~: 

\noindent  (2.4.20) $\qquad \displaystyle 
f_j^{n+1/2} \,\,= \,\, \pmatrix { {{\rho_0 \, c_0^2}\over{2}}  \, \bigl( u_{j-1/2}^{n}
+  u_{j+1/2}^{n}\bigr)  \,-\, {{c_0}\over{2}} \bigl( p_{j+1/2}^{n} -  p_{j-1/2}^{n}
\bigr) \cr {{1}\over{2 \, \rho_0}} \, \bigl( p_{j-1/2}^{n} +  p_{j+1/2}^{n}\bigr)
\,-\,  {{c_0}\over{2}} \bigl( u_{j+1/2}^{n} -  u_{j-1/2}^{n} \bigr) \cr } \, $

\noindent
for the internal fluxes, {\it i.e.} for indexes $j$ that satisfy $\,\,1 \leq j \leq  J
\!-\! 1 \,. \,$ The two boundary fluxes follow the following relations~: 

\noindent  (2.4.21) $\qquad \displaystyle 
f_0^{n+1/2} \,\,= \,\, \pmatrix { \rho_0 \, c_0^2 \,  u_{1/2}^{n} \,+\, c_0 \, \Bigl(
\Pi \bigl( (n+{1\over2}) \Delta t \bigr) -  p_{1/2}^{n} \Bigr) \cr  {{1}\over{
\rho_0}} \, \Pi \bigl( (n+{1\over2}) \Delta t \bigr) \cr} \,$

\noindent  (2.4.22) $\qquad \displaystyle 
f_J^{n+1/2} \,\,= \,\, \pmatrix { {{\rho_0 \, c_0^2}\over{2}}  \, \bigl( u_{J-1/2}^{n}
+  u_{J-1/2}^{0}\bigr)  \,-\, {{c_0}\over{2}} \bigl( p_{J-1/2}^{0} -  p_{J-1/2}^{n}
\bigr) \cr {{1}\over{2 \, \rho_0}} \, \bigl(  p_{J-1/2}^{n} +  p_{J-1/2}^{0}  \bigr)
\,-\,  {{c_0}\over{2}} \bigl( u_{J-1/2}^{0} -  u_{J-1/2}^{n} \bigr) \cr } \, .\, $

\smallskip \noindent   $\bullet \qquad \,\,\, $ 
The internal fluxes are determined with the scheme (2.3.10)  applied with the diagonal
form of relation (2.4.9). We have 

\noindent  (2.4.23) $\qquad \displaystyle 
\varphi_{+ ,\,j}^{n+1/2} \,\,= \,\, \varphi_{+ ,\,j-1/2}^{n} \,\,\equiv \,\,
p_{j-1/2}^{n} \,+\, \rho_0 \, c_0 \,\,  u_{j-1/2}^{n} \,$ 

\noindent  (2.4.24) $\qquad \displaystyle 
\varphi_{- ,\,j}^{n+1/2} \,\,= \,\, \varphi_{- ,\,j+1/2}^{n} \,\,\equiv \,\,
p_{j+1/2}^{n} \,-\, \rho_0 \, c_0 \,\,  u_{j+1/2}^{n} \,$ 

\noindent 
then the relation (2.4.20) is established. 

\noindent 
The left boundary flux uses the extension of relation (2.3.12). We first determine
the characteristic variables on the left boundary according to relation (2.4.13)

\noindent  (2.4.25) $\qquad \displaystyle 
\varphi_{+ ,\,0}^{n+1/2} \,\,= \,\, 2 \Pi \bigl( (n+{1\over2}) \Delta t \bigr) \,-\,
\varphi_{- ,\,0}^{n+1/2} \,$

\noindent 
and use a first order  extrapolation of the outgoing characteristic variable~: 

\noindent  (2.4.26) $\qquad \displaystyle 
\varphi_{- ,\,0}^{n+1/2} \,\,= \,\, \varphi_{- ,\,1/2}^{n} \,\, \equiv \,\, 
p_{1/2}^{n} \,-\, \rho_0 \, c_0 \, \, u_{1/2}^{n} \,. \,$ 

\noindent 
Then we solve the system (2.4.25) (2.4.26) and find finally the relation (2.4.21).
The process   is analogous for the right boundary. 
The input datum is imposed according to the relation (2.4.14)~: 

\noindent  (2.4.27) $\qquad \displaystyle 
\varphi_{- ,\,J}^{n+1/2} \,\,= \,\, \varphi_{- ,\,J}^{0} \,\, \equiv \,\,
p_0(L) \,-\, \rho_0 \, c_0 \,  u_0(L) \,\, \approx \,\, p_{J-1/2}^0 \,-\, \rho_0 \, 
c_0 \,\, u_{J-1/2}^0  \,$

\noindent 
and the output characteristic variable is extrapolated from the interior of the
domain~: 

\noindent  (2.4.28) $\qquad \displaystyle 
\varphi_{+ ,\,J}^{n+1/2} \,\,= \,\, \varphi_{+ ,\,J-1/2}^{n} \,\, \equiv \,\, 
p_{J-1/2}^{n} \,+\, \rho_0 \, c_0 \,\,  u_{J-1/2}^{n} \,. \,$ 

\noindent  
The relation (2.4.22) follows after two steps  of elementary algebra.  $ \hfill
\square \kern0.1mm $

\smallskip \noindent   $\bullet \qquad \,\,\, $ 
We remark that both relations (2.4.20) and (2.4.22) are identical, except that the
boundary state $\,\, W_0(L) \approx W_{J-1/2}^0 \,\,$  has replaced the right state
$\,\, W_{j+1/2}^n \,.\,$ Moreover the  flux boundary condition (2.4.21) that
involves the pressure is a natural discretization of the exact characteristic
solution (2.4.17) (2.4.18) at $ \, x \!=\! 0 \,.\,$

\bigskip \noindent {\smcaps 2.5 } $ \, $ { \bf  Characteristic variables. }

\noindent   $\bullet \qquad \,\,\,\,\, $ 
We suppose now to fix the ideas that the unknown vector $\,\,  W({\scriptstyle \bullet}
,\, {\scriptstyle \bullet} )\,\,$ 

\noindent  (2.5.1) $\qquad \displaystyle 
[0,\,L] \times [0,\, +\infty[ \,\, \ni (x,\,t) \longmapsto W(x,\,t) \in \R^3 \,\,$ 

\noindent 
has three real components $\, w_1 ,\, w_2 \, $ and $ \, w_3 .\,$ We suppose also that
the function $\,\, W({\scriptstyle \bullet} ,\, {\scriptstyle \bullet} )\,\,$ is
solution of a  conservation law of the type

\noindent  (2.5.2) $\qquad \displaystyle 
{{\partial W}\over{\partial t}} \,+\, {{\partial }\over{\partial x}} F(W)  \,\,= \,\,
0  \,$

\noindent
where the  flux $\, F(W) \,$ is a {\bf linear} function of vector $\, W \,: \,$ 

\noindent  (2.5.3) $\qquad \displaystyle 
F(W) \,\,= \,\, A \,{\scriptstyle \bullet}\, W \,$

\noindent 
and $\, A \,$ is a 3 by 3  diagonalizable  real  matrix. 

\smallskip \noindent   $\bullet \qquad \,\,\,\,\, $ 
We first detail the fact that matrix $\, A \,$ is a diagonalizable matrix. There
exists three non null real vectors $\, r_{1} \,,\, $  $\, r_{2}
\,,\, $  $\, r_{3} \, $ and three real scalars $\, \lambda_{1} \,,\, $ 
$\,\lambda_{2} \,,\, $  $\,\lambda_{3} \,\,  $ in such a way that 

\noindent  (2.5.4) $\qquad \displaystyle 
A \,{\scriptstyle \bullet}\,  r_{j}  \,\, = \,\, \lambda_{j} \,  r_{j}
\,\,,\qquad j=1,\, 2 ,\, 3 .\, $

\noindent 
>From a matricial viewpoint, we denote by $\,\, R_{k\,j} \,$ the $\, k^0 \,$
component of the eigenvector $\,  r_{j}  \,$, {\it i.e.} 

\noindent  (2.5.5) $\qquad \displaystyle 
 r_{j}  \,\, = \,\, \pmatrix {  R_{1\,j} \cr  R_{2\,j} \cr  R_{3\,j} \cr}
\,\, \equiv \,\,  \pmatrix { \bigl( r_{j} \bigr)_{1} \cr \bigl( r_{j}
\bigr)_{2} \cr \bigl( r_{j} \bigr)_{3} \cr }  $ 

\noindent
and we introduce the 3 by 3 matrix $\, R \,$ composed by the scalars   $\,\,
R_{k\,j} \,. \,$ The vector $\,  r_{j} \,$ is the $ \, k^0 \,$ column of matrix
$\, R .\,$  The relation (2.5.4) can also be written as 

\noindent  (2.5.6) $\qquad \displaystyle 
A \,{\scriptstyle \bullet}\,  R  \,\, = \,\, R \,{\scriptstyle \bullet}\, \Lambda
\,,\,$ 

\noindent
and $\, \Lambda \,$ is the diagonal matrix whose diagonal terms are equal to the
eigenvalues $\,  \lambda_{j} \,$~: 

\noindent  (2.5.7) $\qquad \displaystyle 
\Lambda \,\,= \,\, \pmatrix { \lambda_{1} & 0 & 0 \cr 0 &  \lambda_{2} & 0 \cr 0
& 0 & \lambda_{3} } \,. \,$

\smallskip \noindent   $\bullet \qquad \,\,\,\,\, $ 
We consider now {\bf two} distinct bases for linear space $\, \R^3 \,:$ on one hand
the canonical basis $\,\, \bigl( e_{j} \bigr)_{j=1,\,2,\,3} \,\,$ defined by 

\noindent  (2.5.8) $\qquad \displaystyle 
e_{1} \,=\, \pmatrix { 1 \cr 0 \cr 0} \,\,,\qquad  e_{2} \,=\, \pmatrix { 0 \cr 1
\cr 0} \,\,,\qquad  e_{3} \,=\, \pmatrix { 0 \cr 0 \cr 1} \,$

\noindent
where the vector $\, W \,$ admits the natural decomposition introduced above~: 

\noindent  (2.5.9) $\qquad  
W \,\,= \,\, \sum_{k=1}^{k=3} \, w_{k} \, e_{k} \,\,, \,$ 

\noindent
and on the other hand the basis of $\, \R^3 \,$ composed by the eigenvectors $\, 
(r_{j})_{j=1,\,2,\,3} .\,$ In the latter, the vector $\, W \,$ can be decomposed
with a formula of the type 

\noindent  (2.5.10) $\qquad  
W \,\,= \,\, \sum_{j=1}^{j=3} \, \varphi_{j} \, r_{j} \, \,$ 

\noindent
and the scalar $\,\, \varphi_{j} \,\,$ define the {\bf characteristic variables}
associated with the system (2.5.2) (2.5.3).  The link between the relations (2.5.9) and
(2.5.10) is classical~: we consider the components $\,\, R_{k\,j} \, \,$ of vector
$\,\,  r_{j} \, \,$ inside the canonical basis and we get from  the relation
(2.5.5)~: 

\noindent  (2.5.11) $\qquad  
w_{k} \,\,=\,\,  \sum_{j=1}^{j=3} \, \varphi_{j} \, R_{k\,j} \, .\,$ 

\noindent
Then the relation (2.5.11) can be re-written under a matricial form~: 

\noindent  (2.5.12) $\qquad \displaystyle 
W \,\,= \,\, R  \,{\scriptstyle \bullet}\, \varphi \,. \,$ 

\smallskip \noindent   $\bullet \qquad \,\,\,\,\, $ 
The relation (2.5.12) proposes to change the unknown function, {\it i.e.} to replace the
research of $\,\, W(x,\,t) \in \R^3 \, \,$ by the  equivalent research of  the
characteristic vector  $\,\, \varphi(x,\,t) \in \R^3 \, \,$ and defined by~:

\noindent  (2.5.13) $\qquad \displaystyle 
\varphi \,\,=\,\, R^{-1}  \,{\scriptstyle \bullet}\, W \,. \,$

\smallskip   \noindent  {\bf Proposition 2.6.} ~ 
% $\hfill $  
{\bf  Characteristic variables satisfy  advection equations.}

\noindent 
The vector $\, \, [0,\,L] \times [0,\, +\infty[ \,\, \ni (x,\,t) \longmapsto \varphi
(x,\,t) \in \R^3 \,\,$ of characteristic variables satisfy the matrix equation 

\noindent  (2.5.14) $\qquad \displaystyle 
{{\partial \varphi}\over{\partial t}} \,+\, \Lambda  \,{\scriptstyle \bullet}\, 
{{\partial \varphi}\over{\partial x}} \,\, = \,\, 0 \,$ 

\noindent
that takes also the equivalent scalar form~: 

\noindent  (2.5.15) $\qquad \displaystyle 
{{\partial \varphi_{j} }\over{\partial t}} \,+\, \lambda_{j} \, 
{{\partial \varphi_{j}}\over{\partial x}} \,\, = \,\, 0 \,\,,\qquad \quad j=1,\,2,\,3
\,.\,$

\smallskip \noindent   $\bullet \qquad \,\,\, $ 
We have from (2.5.2), (2.5.3), (2.5.6) and (2.5.12)~: 

\noindent $ \displaystyle
{{\partial W}\over{\partial t}} \,\, + \,\, A \,  {{\partial W}\over{\partial x}} \,\,
= \,\, R \, {\scriptstyle \bullet} \,  {{\partial \varphi}\over{\partial t}} \,\,+\,\,
A \, {\scriptstyle \bullet} \,  R \, {\scriptstyle \bullet} \, {{\partial
\varphi}\over{\partial x}} \,\, 
=\,\,  R \, {\scriptstyle \bullet} \, \biggl( \,  {{\partial \varphi}\over{\partial t}}
\,\,+\,\, R^{-1} \, {\scriptstyle \bullet} \, A \, {\scriptstyle \bullet} \, R \,
{\scriptstyle \bullet} \,  {{\partial \varphi}\over{\partial x}} \, \biggr) \,$

\noindent $ \displaystyle \qquad \qquad \qquad \quad \,\,
= \,\,  R \, {\scriptstyle \bullet} \, \biggl( \,  {{\partial \varphi}\over{\partial t}}
\,\,+\,\, \Lambda \, {\scriptstyle \bullet} \,  {{\partial \varphi}\over{\partial x}}
\, \biggr) \,\,= \,\, 0 \,, \,$

\noindent
and since the matrix $\, R \,$ is invertible, we deduce from the previous calculus
the relation (2.5.14). The relation (2.5.15) is an immediate consequence of (2.5.14)
and (2.5.7).  $ \hfill \square \kern0.1mm $
 
% version mac 
% \smallskip \vskip 2.7cm  \smallskip
% \qquad  \qquad   \quad   \special{illustration fig.2.7.epsf scaled 500}  
% fin de la version version mac 
% version linux
\smallskip \centerline {  \epsfysize=3,7cm  \epsfbox  {fig.2.7.epsf} }  \smallskip
% fin de la version linux  

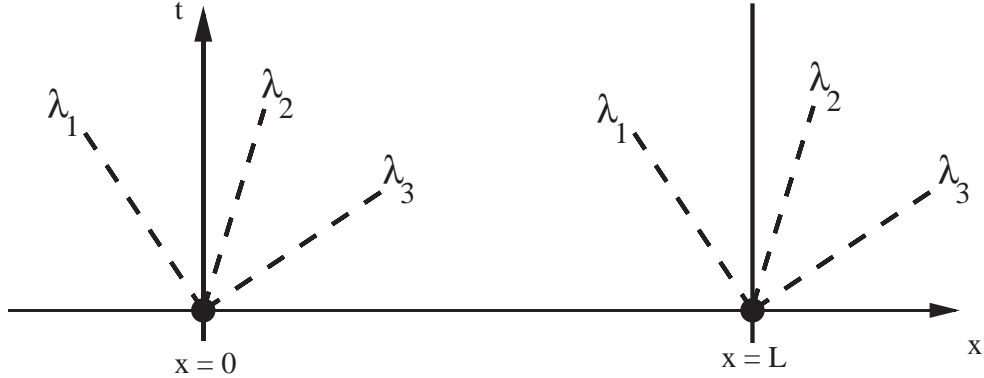
\centerline {\rm  {\bf Figure 2.7.}	\quad    Linear hyperbolic system with three
equations }

\centerline {\rm  and eigenvalues satisfying  $\,\,  \lambda_{1} \,\, < \,\, 0 \,\, <
\,\, \lambda_{2} \,\, < \,\, \lambda_{3} .\,$ }

% \bigskip 
\smallskip \noindent   $\bullet \qquad \,\,\,\,\, $ 
To fix the ideas, we suppose that the eigenvalues $\,\, \lambda_j \,\,$ of matrix $\,
A \,$ are distinct, enumerated with an increasing order and with distinct signs as
illustrated on Figure 2.7~:

\noindent  (2.5.16) $\qquad \displaystyle 
\lambda_{1} \,\, < \,\, 0 \,\, < \,\, \lambda_{2} \,\, < \,\, \lambda_{3}
\,.\,$ 

\noindent 
The propagation of the first variable $\, \varphi_{1} \,$ goes from right to left
(because $\, \lambda_{1} <  0 \,$) with celerity  $\, \abs{\lambda_{1}} , \,$ 
 the second characteristic variable $\, \varphi_{2} \,$ from left to right with
celerity  $\, \lambda_{2} \,$ and the same property holds for variable $\, 
\varphi_{3} \,$ with eigenvalue $\, \lambda_{3} .\,$ 

\smallskip \noindent   $\bullet \qquad \,\,\,\,\, $ 
A set of well posed boundary conditions is a consequence of the diagonal form (2.5.15)
of the equations and of the particular choice (2.5.16) for the signs. The directions
associated with eigenvalues  $\, \lambda_{2} \,$ and  $\, \lambda_{3} \,$ are
ingoing at $\, x \!=\! 0 \,$ and we have to give some boundary condition for $\,
\varphi_{2} \,$ and $\, \varphi_{3} \,$ at this point~: 

\noindent  (2.5.17) $\qquad \displaystyle 
\varphi_{2} (x\!=\!0 ,\, t) \,\,= \,\, \beta_0(t) \,$ 

\noindent  (2.5.18) $\qquad \displaystyle 
\varphi_{3} (x\!=\!0 ,\, t) \,\,= \,\, \gamma_0(t) \,.  \,$ 

\noindent 
The direction associated with the eigenvalue $\, \lambda_{1} \,$ is ingoing at the
abscisssa  $\, x \!=\! L ,\,$ and this condition imposes to have some datum
concerning $\, \varphi_{1} \,$ at this particular point~: 

\noindent  (2.5.19) $\qquad \displaystyle 
\varphi_{1} (x\!=\!L ,\, t) \,\,= \,\, \alpha_L(t) \,. \,$ 

\noindent
The previous boundary conditions (2.5.17) to (2.5.19) define a well posed problem.
Nevertheless, the introduction of physically relevant boundary conditions (as a
pressure condition as seen in the previous section) requires a more general
formulation of the boundary condition. In the linear case, the stability study
developed by Kreiss [Kr70] shows that the ingoing characteristic can be an affine
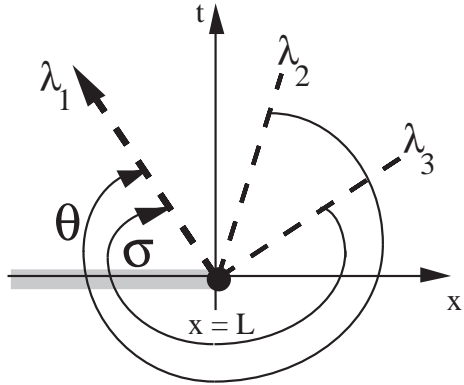
function of the outgoing characteristic through a {\bf reflection operator} at the
boundary. We can explicit the former with the above example. 

\smallskip \centerline {  \epsfysize=3,7cm  \epsfbox  {fig.2.8.epsf} }  \smallskip

\centerline {\rm  {\bf Figure 2.8.}	\quad    Reflection operator at  $\,\, x = 0 .\,$ }

\smallskip \centerline {  \epsfysize=3,7cm  \epsfbox  {fig.2.9.epsf} }  \smallskip

\centerline {\rm  {\bf Figure 2.9.}	\quad    Reflection operator at  $\,\, x = L .\,$ }

\bigskip 
\noindent   $\bullet \qquad \,\,\,\,\, $ 
At $\, x \!=\! 0 \,,\,$ the first characteristic is outgoing and the two last ones
are going inside the domain of study. Then we can replace the conditions (2.5.17)
and (2.5.18) by the following ones~: 

\noindent  (2.5.20) $\qquad \displaystyle 
\varphi_{2} (x\!=\!0 ,\, t) \,\,= \,\, \beta_0(t) \,\,+\,\, p(t) \,\varphi_{1}
(x\!=\!0 ,\, t) \,  $ 

\noindent  (2.5.21) $\qquad \displaystyle 
\varphi_{3} (x\!=\!0 ,\, t) \,\,= \,\, \gamma_0(t) \,\,+\,\, q(t) \,\varphi_{1}
(x\!=\!0 ,\, t) \,\,, \, $  

\noindent 
where $\,\, t \longmapsto p(t)  \, \,$ and  $\,\, t \longmapsto q(t)  \, \,$  are
given fixed real functions of time. The conditions (2.5.20) and (2.5.21) are
illustrated on Figure 2.8. We can also write them 

\noindent  (2.5.22) $\qquad \displaystyle 
\varphi^{in}   (x ,\, t) \,\,= \,\, g(t) \,\,+\,\, S(t) \, {\scriptstyle \bullet}
\, \varphi^{out}   (x ,\, t) \,, \, \quad  x \, $ point on the boundary, 

\noindent  with  $  \,\,\,\,  \displaystyle 
\varphi^{in} =  \pmatrix { \varphi_{2} \cr  \varphi_{3}\cr } \,,\,\, g(t)=
\pmatrix {  \beta_0(t) \cr   \gamma_0(t) } \,,\,\,  S(t) = \pmatrix {  p(t) \cr 
q(t) \cr }   \,,\,\,   \varphi^{out} = \varphi_{1}  \,  $

\noindent 
 when $\, x = 0 \,.\,$  

\smallskip \noindent   $\bullet \qquad \,\,\,\,\, $ 
When $\, x \!=\!L ,\,$ the relation   (2.5.19) is replaced by a more general one 

\noindent  (2.5.23) $\qquad \displaystyle 
\varphi_{1} (x\!=\!L ,\, t) \,\,= \,\, \alpha_L(t) \,\,+\,\,  \theta(t) \,
\,  \varphi_{2} (x\!=\!L ,\, t) \,\, + \,\, \sigma(t) \,\, \varphi_{3} (x\!=\!L ,\,
t) \,\,$ 

\noindent 
illustrated on Figure 2.9 and including  an affine component of the outgoing
characteristic variables. The boundary condition (2.5.23) takes again a form of the
type  (2.5.22) with this time  the  following relations~: 
 $ \displaystyle
\varphi^{in}  \,= \, \varphi_{1} \,\,, \quad g(t) \,=\, \alpha_L(t)  \,\,, \quad $ 
$ \displaystyle  S(t) \,=\, \bigl( \theta(t)  \,\, \sigma(t)  \bigr)  \,\,, \quad
\varphi^{out}  \,=\, \pmatrix { \varphi_{2} \cr  \varphi_{3}\cr }  \,  $  when
$\,\, x\!=\! L \,.\,$  

\bigskip \noindent {\smcaps 2.6 } $ \, $ { \bf   A family of model systems with three
equations }

\noindent   $\bullet \qquad \,\,\,\,\, $ 
We still study a 3 by 3 linear hyperbolic system of the type (2.5.2) (2.5.3) with the
condition (2.5.16) to fix a particular  example.  We suggest in this section to
explicit a way for evaluation of the numerical flux $\,\, f_j^{n+1/2} \, \,$ that is
the key point for the discrete evolution in time of the mean values $\,\, W_{j+1/2}
\,: \,$ 

\noindent  (2.6.1) $\qquad \displaystyle 
{{1}\over{\Delta t}} \Bigl( W_{j+1/2}^{n+1} \,-\,  W_{j+1/2}^{n} \Bigr) \,\,+\,\,
{{1}\over{\Delta x}} \Bigl( f_{j+1}^{n+1/2} \,-\,  f_{j}^{n+1/2} \Bigr) \,\,= \,\, 0
\,.\, $

\noindent
The internal fluxes $\,\, \bigl( f_j^{n+1/2} \bigr)_{j=1,\cdots,\, J \!-\! 1} \,\,$
are evaluated with the help of a two-point numerical flux function $\,\,
\Phi({\scriptstyle \bullet} ,\, {\scriptstyle \bullet}) \,: \,$ 

\noindent  (2.6.2) $\qquad \displaystyle 
 f_{j}^{n+1/2} \,\, = \,\, \Phi ( \,  W_{j-1/2}^{n} ,\,  W_{j+1/2}^{n}) \,\,$ 

\noindent
and the boundary fluxes $\,\,  f_{0}^{n+1/2} \,\, $ and $\,\,  f_{J}^{n+1/2} \,\, $
are detailed in a forthcoming sub-section. 

% version mac 
% \smallskip \vskip 1.3cm  \smallskip 
% \qquad  \qquad  \qquad \qquad  \qquad  \qquad  \special{illustration fig.2.10.epsf scaled 500}  
% fin de la version version mac 
% version linux
\smallskip \centerline {  \epsfysize=2,3cm  \epsfbox  {fig.2.10.epsf} }  \smallskip
% fin de la version linux     

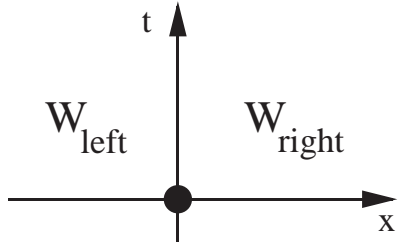
\centerline {\rm  {\bf Figure 2.10.}	\quad   Discontinuity at the interface between two
cells. }

\bigskip 
\noindent   $\bullet \qquad \,\,\,\,\, $ 
We change the notations and wish to determine the numerical flux $\,\,
\Phi(W_l,$ $\,W_r) \,\,$ for $\, W_l = W_{\rm left} \,$ and $\,  W_r = W_{\rm right}
\,$ given respectively at the left and at the right of the interface (see Figure 2.10).
When we consider the advection equation (and in that case the variables  $\, W_l \,$
and  $\, W_r \,$ are real numbers) the relation (2.1.19) gives the result~: $\,
\Phi(W_l,\,W_r) \,= \,a \, W_l \,$ when $ \,   a > 0 \, $ and  $\,
\Phi(W_l,\,W_r) \,= \,a \, W_r \,$ when $ \,   a < 0 \,.\,  $ We have to generalize
this study when the field $\, W({\scriptstyle \bullet},\, {\scriptstyle \bullet}) \,$
is three-dimensional. We first decompose the vector $\,\, \Phi (W_l,$ $\,W_r) \,\,$
with the basis $\, r_{j} \,$ of eigenvectors and introduce its (scalar)
components $\,\, \psi_{j}(W_l,\,W_r) \,: \,$ 

\noindent  (2.6.3) $\qquad  
\Phi(W_l,\,W_r) \,\,= \,\, \sum_{j=1}^{j=3} \,\psi_{j} (W_l,\,W_r) \,\,  r_{j}
\,$ 

\noindent
{\it i.e. }

\noindent  (2.6.4) $\qquad  
\Phi_{k}(W_l,\,W_r) \,\,= \,\, \sum_{j=1}^{j=3} \,  R_{k\,j} \, \psi_{j} (W_l,
\,W_r) \,. \,$

\noindent 
For $\, j \!=\! 1 ,\,$ we have $\, \lambda_{1} < 0 \,$ then the numerical scheme has
to be upwinded in the right direction~: 

\noindent  (2.6.5) $\qquad \displaystyle 
\psi_{1} (W_l, \,W_r) \,\,= \,\, \lambda_{1} \,  \varphi_{1,\,r} \,$

\noindent 
whereas for $\, j\!=\! 2 \,$ or $\, j\!=\! 3 ,\,$ we have  $\, \lambda_{2} > 0 \,$  and
$\, \lambda_{3} > 0 \,$ and the scheme must be upwinded to the left. It comes 

\noindent  (2.6.6) $\qquad \displaystyle 
\psi_{2} (W_l, \,W_r) \,\,= \,\, \lambda_{2} \, \varphi_{2,\,l} \,, \qquad 
\psi_{3} (W_l, \,W_r) \,\,= \,\, \lambda_{3} \, \varphi_{3,\,l} \,.\, $

\noindent
In consequence of the relations (2.6.3) to (2.6.6), the numerical flux function $\,\,
\Phi({\scriptstyle \bullet},\,{\scriptstyle \bullet}) \,\,$ can be written globally~: 

\noindent  (2.6.7) $\qquad \displaystyle 
\Phi(W_l,\,W_r) \,\,= \,\,\lambda_{1} \, \varphi_{1,\, r} \, r_{1} \,\,+\,\, 
\lambda_{2} \,  \varphi_{2,\,l} \, r_{2} \,\,+\,\, \lambda_{3} \,
 \varphi_{3,\,l} \, r_{3} \, \,,\,$

\noindent
or in an equivalent way with introducing the Cartesian components~: 

\noindent  (2.6.8) $\quad \displaystyle 
\Phi_{k} (W_l,\,W_r) \,\,= \,\,\lambda_{1} \,  \varphi_{1,\,r} \,  R_{k\,1} \, 
\,\,+\,\,  \lambda_{2} \, \varphi_{2,\,l} \,  R_{k\,2} \,\,+\,\, \lambda_{3} \,
\varphi_{3,\,l} \,  R_{k\,3} \,, \,\,\, k\!=\! 1,\, 2 ,\, 3 \,. \,$

\smallskip \noindent   $\bullet \qquad \,\,\,\,\, $ 
We can also re-write the relation (2.6.8) for the particular interface  $\,\, x_{j}
\,:\,$ 

\noindent  (2.6.9) $\qquad \displaystyle 
W_l \,\, = \,\, W_{\rm left} \,\,= \,\, W_{j-1/2}^n \,\,, \qquad  W_r \,\, = \,\,
W_{\rm right} \,\,= \,\, W_{j+1/2}^n \,.\, $

\noindent 
We first  decompose the vector $\, W \,$ on the eigenvectors of matrix $\, A \,$ as
in (2.5.11)~: 

\noindent  (2.6.10) $\quad  
\bigl( W_{j+1/2}^n \bigr)_k \,\, = \,\, \sum_{i=1}^{i=3} \, \varphi_{i,\, j+1/2}^n \, \,
 R_{k\,i} \,, \quad  k = 1,\, 2 ,\, 3 \,, \quad j = 1 ,\, \cdots ,\, J \! - \! 1 \,, \,$
 
\noindent
then we introduce the component  number $\, k \,$ of the flux $\,\, f_j^{n+1/2}  \,,\,$
{\it i.e.} $\,\,\,  (  f_j^{n+1/2} )_k \, $ 
$ =\, \Phi_{k} (W_{j-1/2}^n ,\,W_{j+1/2}^n) \, $
at the interface  $\,\, x_{j} \,:\,$  
 
\noindent  (2.6.11) $\qquad \displaystyle 
\bigl(  f_j^{n+1/2} \bigr)_k \,=\,\lambda_{1} \,  \varphi_{1,\, j+1/2}^n \, 
R_{k\,1} \,+\,  \lambda_{2} \,   \varphi_{2,\, j-1/2}^n \,  R_{k\,2} \,+\,
\lambda_{3} \, \varphi_{3,\, j-1/2}^n  \,  R_{k\,3} \,. \,$

\smallskip \noindent   $\bullet \qquad \,\,\,\,\, $ 
We detail in this sub-section  the determination of the numerical flux $\,\, 
f_0^{n+1/2} \,\,$  at the boundary $\, x \!=\! 0 .\, $ We first recall that the
continuous boundary conditions at this point take the form given in (2.5.20) (2.5.21).
The idea is  to try to apply the upwind scheme (2.6.11) at the particular vertex $\,
j\!=\! 0 \,: \, \,  $   $  f_0^{n+1/2}  \,\,=\,\,\lambda_{1} \,  \varphi_{1,\,
1/2}^n \, r_{1}  \,+\,$ $+\,   \lambda_{2} \,   \varphi_{2,\, -1/2}^n \,r_{2} 
\,+\, \lambda_{3} \, \varphi_{3,\, -1/2}^n   \,r_{3}  \,\,\,\, $  and then to
replace the characteristic values $\,\,  \varphi_{2,\, -1/2}^n \,\, $ and $\,\,
\varphi_{3,\, -1/2}^n \,\,$ (that are not defined on the mesh) by their values
evaluated after a rough discretization of relations (2.5.20) and (2.5.21)~: 
 $  \,\,\,  \varphi_{2,\, -1/2}^n \,=\,  \beta_0^{n+1/2} \,+\,$ $+\,  p^{n+1/2} \, 
\varphi_{1,\, 1/2}^n \,, \, $ $\quad  \varphi_{3,\, -1/2}^n \,=\, 
\gamma_0^{n+1/2} \,+\, q^{n+1/2} \,  \varphi_{1,\, 1/2}^n \,.\,\,\,   $ We obtain in
consequence the following expression for the {\bf boundary flux} at $\, x \!=\! 0 \,:
\, $

\setbox21=\hbox {$\displaystyle \lambda_{1} \,  \varphi_{1,\, 1/2}^n \, r_{1}  \,+\, 
\lambda_{2} \,  \bigl(\, \beta_0^{n+1/2} \,\,+\,\, p^{n+1/2} \, 
\varphi_{1,\, 1/2}^n \, \bigr)  \,r_{2}  \,+\,  $}
\setbox22=\hbox {$\displaystyle    \qquad \,+\,
\lambda_{3} \,  \bigl(\, \gamma_0^{n+1/2} \,\,+\,\, q^{n+1/2} \,  \varphi_{1,\,
1/2}^n \, \bigr) \,r_{3}  \, \,$}
\setbox30= \vbox {\halign{#&# \cr \box21 \cr \box22    \cr   }}
\setbox31= \hbox{ $\vcenter {\box30} $}
\setbox44=\hbox{\noindent  (2.6.12) $\displaystyle  \qquad  
 f_0^{n+1/2}  \,\,=\,\,  \left\{ \box31 \right.$}  
\noindent $ \box44 $

\noindent
or in an equivalent way~: 

\setbox21=\hbox {$\displaystyle  \varphi_{1,\, 1/2}^n \, \bigl( \, 
\lambda_{1} \, r_{1}  \,+\, \lambda_{2} \, p^{n+1/2}  \, r_{2} \,+\, 
\lambda_{3} \,  q^{n+1/2} \,r_{3}  \, \bigr) \,\,+ \, $} 
\setbox22=\hbox {$\displaystyle   \qquad \,+\,\,  \lambda_{2} \,\beta_0^{n+1/2}
\,r_{2} \,\,+\,\,\lambda_{3} \, \gamma_0^{n+1/2} \,r_{3} \,. \,$} 
\setbox30= \vbox {\halign{#&# \cr \box21 \cr \box22    \cr   }}
\setbox31= \hbox{ $\vcenter {\box30} $}
\setbox44=\hbox{\noindent  (2.6.13) $\displaystyle  \qquad  
 f_0^{n+1/2}  \,\,=\,\,  \left\{ \box31 \right.$}  
\noindent $ \box44 $

\smallskip \noindent   $\bullet \qquad \,\,\,\,\, $ 
The determination of the {\bf  boundary flux} $\,\,  f_J^{n+1/2} \,\,$ can be conducted
in the same way. Starting from the  expression of the upwind scheme (2.6.11) when
$\, j \!=\! J \,,\,$ {\it i.e.} formally  $\,\,\,\,\,   f_J^{n+1/2}  \,=\,\lambda_{1} \, 
\varphi_{1,\, J+1/2}^n \, r_{1} \,+\,  \lambda_{2} \,   \varphi_{2,\, J-1/2}^n
\,  r_{2}  \,+\, \lambda_{3} \, \varphi_{3,\, J-1/2}^n  \,  r_{3} \,\,,\,\, $
we replace the first characteristic variable that appears external of the domain by
its value given by the boundary condition (2.5.23)~: $\,\, \varphi_{1,\, J+1/2}^n \,=
\,\alpha_L^{n+1/2} \,+\, \theta^{n+1/2} \, \varphi_{2,\, J-1/2}^n \,+\, $ $ +\, 
\sigma^{n+1/2} \, \varphi_{3,\, J-1/2}^n \,.\,$ We deduce~: 

\setbox21=\hbox {$\displaystyle  \lambda_{1} \, \bigl( \, \alpha_L^{n+1/2}
\,+\,\theta^{n+1/2} \, \varphi_{2,\, J-1/2}^n \,+\, 
\sigma^{n+1/2} \, \varphi_{3,\, J-1/2}^n \, \bigr) \,  r_{1} \,   $}
\setbox22=\hbox {$\displaystyle   \qquad  \qquad  \qquad 
\,+\,\,  \lambda_{2} \,    \varphi_{2,\, J-1/2}^n \,  r_{2}
 \,+\, \lambda_{3} \, \varphi_{3,\, J-1/2}^n  \,  r_{3} \,  $}
\setbox30= \vbox {\halign{#&# \cr \box21 \cr \box22    \cr   }}
\setbox31= \hbox{ $\vcenter {\box30} $}
\setbox44=\hbox{\noindent  (2.6.14) $\displaystyle  \quad   f_J^{n+1/2} \,\,= \,\,
 \left\{ \box31 \right.$}  
\noindent $ \box44 $

\noindent
or in an equivalent manner~: 

\setbox21=\hbox {$\displaystyle  \lambda_{1} \, \alpha_L^{n+1/2} \,  r_{1}  \, 
+\,  \varphi_{2,\, J-1/2}^n \, \bigl( \, \lambda_{1} \, \theta^{n+1/2} \,  r_{1}
\,+\,   \lambda_{2} \,   r_{2} \, \bigr) \,+\,    $}
\setbox22=\hbox {$\displaystyle  \qquad  \qquad  \qquad 
\,+\,\,  \varphi_{3,\, J-1/2}^n \, \bigl( \,
\lambda_{1} \,\sigma^{n+1/2} \,  r_{1}  \,  +\,    \lambda_{3} \,   r_{3} \, 
\bigr) \,\,. \,  $}
\setbox30= \vbox {\halign{#&# \cr \box21 \cr \box22    \cr   }}
\setbox31= \hbox{ $\vcenter {\box30} $}
\setbox44=\hbox{\noindent  (2.6.15) $\displaystyle  \quad   f_J^{n+1/2} \,\,= \,\,
 \left\{ \box31 \right.$}  
\noindent $ \box44 $

\bigskip \noindent {\smcaps 2.7 } $ \, $ { \bf  First order upwind-centered finite
volumes }

\noindent   $\bullet \qquad \,\,\,\,\, $ 
We consider now a general system of conservation laws 

\noindent  (2.7.1) $\qquad \displaystyle 
{{\partial W}\over{\partial t}} \,+\, {{\partial }\over{\partial x}} F(W)  \,\,= \,\,
0  \,$

\noindent
with an  unknown vector $\,\,  W({\scriptstyle \bullet} ,\, {\scriptstyle \bullet}
)\,\,$ that belongs to linear space $\,\, \R^m \,: \,$ 

\noindent  (2.7.2) $\qquad \displaystyle 
[0,\,L] \times [0,\, +\infty[ \,\, \ni (x,\,t) \longmapsto W(x,\,t) \in \R^m \,\,$ 

\noindent 
and a {\bf linear} flux function $\,\, F( {\scriptstyle \bullet} ) \,\,$ 

\noindent  (2.7.3) $\qquad \displaystyle 
F(W) \,\,= \,\, A \,{\scriptstyle \bullet}\, W \,$

\noindent 
associated with  a diagonalizable matrix $\, A \,$ with eigenvalues $\,\, \lambda_{j}
\,  \,$ and eigenvectors $\,\,   r_{j} \,\,$ 

\noindent  (2.7.4) $\qquad \displaystyle 
A \,{\scriptstyle \bullet}\,  r_{j}  \,\, = \,\, \lambda_{j} \,  r_{j}
\,\,,\qquad j=1,\, 2 ,\, \cdots,\, m \,. \, $

\noindent 
Introducing the $\, m \times m \,$ matrix $\, R \,$ as in relation (2.5.5) and the
diagonal matrix $\, \Lambda \,$ of eigenvalues as in relation (2.5.7), we have~: 

\noindent  (2.7.5) $\qquad \displaystyle 
A \,{\scriptstyle \bullet}\,  R  \,\, = \,\, R \,{\scriptstyle \bullet}\, \Lambda
\,.\,$ 

\smallskip \noindent   $\bullet \qquad \,\,\,\,\, $ 
We propose here to determine a first order upwind flux $\,\, \Phi(W_l ,\, W_r ) \,\,$
between the two states $\, W_{\rm left} = W_l \,  $ and $ \, W_{\rm right} = W_r  \,$
that generalizes the relation (2.6.7) when we have not done any  hypothesis of the
type (2.5.16) concerning the sign of the eigenvalues $\, \lambda_{j} .\,$ We
decompose any state $\, W \,$ on the basis of space $\, \R^m \,$ characterized  by the
eigenvectors~$\, r_{j} \,:\,$ 

\noindent  (2.7.6) $\quad  
W \,\,= \,\, \sum_{j=1}^{j=m} \, \varphi_{j} \, r_{j} \,\,,\quad
W_l\,\,= \,\, \sum_{j=1}^{j=m} \, \varphi_{j,\,l} \, r_{j} \,\,,\quad
W_r\,\,= \,\, \sum_{j=1}^{j=m} \, \varphi_{j,\,r} \, r_{j} \,,\,$

\noindent
and due to  the structure introduced at Proposition 2.6, we obtain an advection
equation for the $\,j^o \,$ characteristic variable $\,\,  \varphi_{j} \,:\,$ 

\noindent  (2.7.7) $\qquad \displaystyle 
{{\partial \varphi_{j} }\over{\partial t}} \,+\, \lambda_{j} \, 
{{\partial \varphi_{j}}\over{\partial x}} \,\, = \,\, 0 \,\,,\qquad \quad j=1,\,2
,\, \cdots,\, m \,. \, $

\noindent
Therefore it is natural to introduce the components $\,\, \psi_{j} (W_l,\,W_r)
\,\,$ of the numerical flux on the basis of the eigenvectors~: 

\noindent  (2.7.8) $\qquad  
\Phi(W_l,\,W_r) \,\,= \,\, \sum_{j=1}^{j=3} \, \psi_{j} (W_l,\,W_r) \,  r_{j}
\,$ 

\noindent
and the first order upwind finite volume scheme is defined by the way we evaluate the
coefficient $\,\,  \psi_{j} (W_l,\,W_r) \,\, $ with the upwind scheme associated
with the advection equation (2.7.7)~: 

\setbox21=\hbox {$\displaystyle  \lambda_{j} \,\,  \varphi_{j,\,l} \,\qquad {\rm
if} \,\,  \lambda_{j} > 0 \,$}
\setbox22=\hbox {$\displaystyle    \lambda_{j} \,\,  \varphi_{j,\,r} \qquad {\rm
if} \,\,  \lambda_{j} < 0 \,.\, $}
\setbox30= \vbox {\halign{#&# \cr \box21 \cr \box22    \cr   }}
\setbox31= \hbox{ $\vcenter {\box30} $}
\setbox44=\hbox{\noindent  (2.7.9) $\displaystyle  \qquad   \psi_{j} (W_l,\,W_r)
\,\,=\,\,   \left\{ \box31 \right.$}  
\noindent $ \box44 $

\smallskip \noindent   $\bullet \qquad \,\,\,\,\, $ 
For any real number $\,\, \mu \,,\,$  we introduce the positive part $\, \, \mu^+ \,\,$
and the negative part $\,\, \mu^- \,\,$ by the relations

\setbox21=\hbox {$\displaystyle  \mu \qquad {\rm if } \,\, \mu \geq 0 \,$}
\setbox22=\hbox {$\displaystyle  0 \, \qquad {\rm if } \,\, \mu \leq 0 \,$}
\setbox30= \vbox {\halign{#&# \cr \box21 \cr \box22    \cr   }}
\setbox31= \hbox{ $\vcenter {\box30} $}
\setbox51=\hbox {$\displaystyle  0 \, \qquad  {\rm if } \,\, \mu \geq 0 \,$}
\setbox52=\hbox {$\displaystyle  \mu \qquad {\rm if } \,\, \mu \leq 0 \,. \,$}
\setbox53= \vbox {\halign{#&# \cr \box51 \cr \box52    \cr   }}
\setbox54= \hbox{ $\vcenter {\box53} $}
\setbox44=\hbox{\noindent  (2.7.10) $\displaystyle  \qquad    \mu^+ \,\,=\,\,  
\left\{ \box31 \right. \,, \qquad   \mu^- \,\,=\,\, \left\{ \box54 \right. $}  
\noindent $ \box44 $

\noindent
We remark that we have 

\noindent  (2.7.11) $\qquad  \displaystyle 
\,\,  \mu \,  \,\, \equiv \, \, \,  \mu^+ \,+ \, \mu^-  \,\,,\qquad \forall \, \mu
\in \R \,$ 
 
\noindent  (2.7.12) $\qquad  \displaystyle 
\abs{\mu} \,\, \equiv \, \,  \mu^+ \,- \, \mu^-  \,\,,\qquad \forall \, \mu \in \R
\,. \,$ 
 
\noindent
We introduce also the absolute value $\,\, \abs{ \Lambda } \,\, $  of the diagonal
matrix $\,\, \Lambda \,\,$ by the condition~: 

\noindent  (2.7.13) $\qquad  \displaystyle 
\abs{ \Lambda } \,\, \equiv \,\,  \abs{ {\rm diag} \bigl( \, \lambda_{1} ,\cdots,\,
\lambda_{m} \,\bigr) } \,\,= \,\, {\rm diag} \bigl( \, \abs{\lambda_{1}} ,\cdots,\,
\abs{\lambda_{m}} \,\bigr) \,$
 
\noindent
and due to the relation (2.7.5), the  absolute value $\,\, \abs{ A } \,\, $  of the
matrix $\, A \,$ is defined by~: 

\noindent  (2.7.14) $\qquad  \displaystyle  
\abs{ A } \,\,= \,\, R \, {\scriptstyle \bullet} \, \abs{ \Lambda } \, {\scriptstyle
\bullet} \, R^{-1} \,. \,$

\bigskip  \noindent  {\bf Proposition 2.7.  $\quad$  Three expressions of the
upwind first order scheme.}

\noindent 
Let $\,\, \Phi(W_l,\,W_r) \,\,$ the upwind flux defined by the relations (2.7.8) and
(2.7.9). Then we have the three equivalent expressions~: 

\noindent  (2.7.15) $\qquad  
\Phi(W_l,\,W_r) \,\,=\,\, F(W_l) \,+\, \sum_{j=1}^{j=m} \, \lambda_{j}^-\,\,  \bigl(
 \varphi_{j,\,r} \,- \,  \varphi_{j,\,l}   \bigr) \,\,   r_{j} \,$
 
\noindent  (2.7.16) $\qquad  
\Phi(W_l,\,W_r) \,\,=\,\, F(W_r) \,-\, \sum_{j=1}^{j=m} \, \lambda_{j}^+\,\,  \bigl(
 \varphi_{j,\,r} \,- \,  \varphi_{j,\,l}   \bigr) \,\,   r_{j} \,$
 
\noindent  (2.7.17) $\qquad  
\Phi(W_l,\,W_r) \,\,=\,\, {1\over2} \, \bigl(  F(W_l) \,+\,F(W_r) \bigr) \,-\,
{1\over2} \,\abs{A}\, {\scriptstyle \bullet} \, (W_r \,-\, W_l \bigr) \,. \,$

\smallskip \noindent   $\bullet \qquad \,\,\, $ 
We write the relation (2.7.9) under the form~: 

\noindent  (2.7.18) $\qquad  \displaystyle
\psi_{j} (W_l,\,W_r) \,\,= \,\,  \lambda_{j}^+ \,\, \varphi_{j,\,l} \,\,+\,\, 
\lambda_{j}^- \,\, \varphi_{j,\,r} \,$

\noindent 
and we have~: 

\noindent $ 
\Phi(W_l,\,W_r) \,\,= \,\, \sum_{j=1}^{j=m} \,\, \bigl( \, \lambda_{j}^+ \,\,
\varphi_{j,\,l} \,\,+\,\,  \lambda_{j}^- \,\, \varphi_{j,\,r} \, \bigr) \, 
r_{j} \,$

\noindent $  \qquad \qquad \quad
\,\,= \,\, \sum_{j=1}^{j=m} \,\, \bigl( \, (\lambda_{j} \,- \, \lambda_{j}^-) \,
\, \varphi_{j,\,l} \,\,+\,\, 
\lambda_{j}^- \,\, \varphi_{j,\,r} \, \bigr) \,  r_{j} \,$ \quad due to (2.7.11)

\noindent $  \qquad \qquad \quad
\,\,= \,\, \sum_{j=1}^{j=m} \,\, \lambda_{j} \,\,  \varphi_{j,\,l} \,\,r_{j}
\,\,+ \,\,  \sum_{j=1}^{j=m} \,\, \lambda_{j}^- \,\,  \bigl( \varphi_{j,\,r} \,- \,
\varphi_{j,\,l} \bigr) \, r_{j} \,$ 

\noindent
and the relation (2.7.15) is established. In an analogous way, we have~: 

\noindent $ 
\Phi(W_l,\,W_r) \,\,= \,\, \sum_{j=1}^{j=m} \,\, \bigl( \, \lambda_{j}^+ \,\,
\varphi_{j,\,l} \,\,+\,\,  \lambda_{j}^- \,\, \varphi_{j,\,r} \, \bigr) \, 
r_{j} \,$

\noindent $  \qquad \qquad \quad
\,\,= \,\, \sum_{j=1}^{j=m} \,\, \bigl( \, \lambda_{j}^+ \,\,
\varphi_{j,\,l} \,\,+\,\, ( \lambda_{j} \,-\,  \lambda_{j}^+ ) \, \,
\varphi_{j,\,r} \, \bigr) \,  r_{j} \,$ \quad due to (2.7.11)

\noindent $  
\Phi(W_l,\,W_r) \,\,=  \,\, \sum_{j=1}^{j=m} \,\, \lambda_{j} \,\,  \varphi_{j,\,r}
\,\,r_{j} \,\,- \,\,  \sum_{j=1}^{j=m} \,\, \lambda_{j}^+ \,\, 
\bigl( \varphi_{j,\,r} \,- \, \varphi_{j,\,l} \bigr) \, r_{j} \,$ 

\noindent
and the relation (2.7.16) holds. We remark that

\vfill \eject                %%%%%%%%%%%%%  15 janvier 2011 
\noindent $ 
\abs{A} \, {\scriptstyle \bullet} \, (W_r - W_l) \,\,= \,\, R  \, {\scriptstyle
\bullet} \,\abs{\Lambda} \,  {\scriptstyle \bullet} \, R^{-1} \,  {\scriptstyle
\bullet} \, R \,  {\scriptstyle \bullet} \, (\varphi_r - \varphi_l) \quad $ due to
(2.7.14) and (2.5.12) 

\noindent $  \qquad \qquad \qquad \quad \,\,\,\,\,
= \,\,  R  \, {\scriptstyle \bullet} \,\abs{\Lambda} \,  {\scriptstyle \bullet} \,
 (\varphi_r - \varphi_l) \quad $

\noindent $  \qquad \qquad \qquad \quad \,\,\,\,\,
= \,\,  \sum_{k=1}^{k=m} \,  \sum_{j=1}^{j=m} \,  R_{k \,j} \, \abs{\lambda_j} \,
(\varphi_{j ,\, r} - \varphi_{j,\,l}) \,e_{k} \,\quad  $ then 

\noindent  (2.7.19) $\qquad  
\abs{A} \, {\scriptstyle \bullet} \, (W_r - W_l) \,\,= \,\,  \sum_{j=1}^{j=m} \, 
\abs{\lambda_j} \, (\varphi_{j ,\, r} - \varphi_{j,\,l}) \,r_{j} \,.\,  $

\noindent
We add the  previous results (2.7.15) with (2.5.16), and we divide by two. 
We obtain~: 

\noindent $ 
\Phi(W_l,\,W_r) \,\,= \,\,  {1\over2} \, \bigl(  F(W_l) \,+\,F(W_r)\bigr) \,\,-\,\,
 {1\over2} \, \sum_{j=1}^{j=m} \,\, \bigl( \lambda_{j}^+ \,-\, \lambda_{j}^- \bigr)
\, \bigl( \varphi_{j,\,r} \,- \,  \varphi_{j,\,l}   \bigr) \,\,   r_{j} \,$

\noindent $  
\,\,= \,\,   {1\over2} \, \bigl(   F(W_l) \,+\,F(W_r) \bigr) \,\,-\,\,  {1\over2} \, 
\sum_{j=1}^{j=m} \,\, \abs{ \lambda_{j}} \, \bigl( \varphi_{j,\,r} \,- \, 
\varphi_{j,\,l}   \bigr) \,\,   r_{j}   \hfill$ \quad due to (2.7.12)

\noindent $  
\,\,= \,\,   {1\over2} \, \bigl(  F(W_l) \,+\,F(W_r) \bigr) \,\,-\,\, {1\over2} \,
\abs{ A} \, {\scriptstyle \bullet} \, ( W_r \,-\, W_l ) \,$ \quad 

\noindent
due to the relation (2.7.19). Then the relation (2.7.17) is established
and the proposition 2.7 is proven.  $ \hfill \square \kern0.1mm $

\bigskip   \bigskip   \bigskip 
\noindent  {\smcaps 3) $ \,\,\,\,\,\,\, $   Gas dynamics with  the Roe method.} 
\smallskip \noindent {\smcaps 3.1 } $ \,  $ { \bf Nonlinear acoustics in one space
dimension.   }

\noindent   $\bullet \qquad \,\,\,\,\, $ 
We propose here to describe quickly  a physical problem that comes from the
theoretical modelling of trombone, detailed for instance 
 in the work of  Hirsch\-berg {\it et al} 
[HGMW96]   or in our study [MD99] with R. Msallam. In a first approximation, the duct
of a  trombone is a long cylinder with a constant section  and the acoustic waves 
propagate  only 
in the  longitudinal direction. We can use a one-dimensional description
of the geometry (see Figure 3.1)  and in what follows, the trombone is modelled by a
real space variable $\, x \,$ that ranges from 
 $\, x \!=\! 0 \,$ at the input to  $\, x \!=\! L \,$ at the output. 

\smallskip \noindent   $\bullet \qquad \,\,\,\,\, $ 
At the input $\, x \!=\! 0 , \,$ a given non-stationary   pressure wave $\,\, t
\longmapsto \Pi(t) \,$ is emitted~; this wave is a perturbation of the ambiant
pressure $\,\, p_0 \,\,$ of the air~: 

\noindent  (3.1.1) $\qquad \displaystyle
\abs { \Pi(t) - p_0} \,\,\, < < \, \,\, p_0 \,\,, \qquad t > 0 \,. \,$

\noindent
At the output $\, x \!=\! L ,\,$ the waves go outside without any reflection due to
the presence of a pavilion and the boundary condition is a  ``free output'' and a
{\bf nonreflecting boundary condition} has to be used. At the initial time $\, t
\!=\! 0 ,\,$ we can consider that the air satisfies the usual conditions of
pressure $\, p(x,\,0) \equiv p_0 \, , \,$  temperature $\, T(x,\,0) \equiv T_0  \,$ and
density $\, \rho(x,\,0) \equiv \rho_0  . \, \,$ We study in this section  a finite
volume  method able to treat nonlinearities in the acoustic modelling and based on
the characteristic  decompositions developed in the previous section.

% \bigskip  \smallskip 
\vfill \eject                %%%%%%%%%%%%%  15 janvier 2011 
\centerline {  \epsfysize=1,6cm  \epsfbox  {fig.3.1.epsf} }  \smallskip

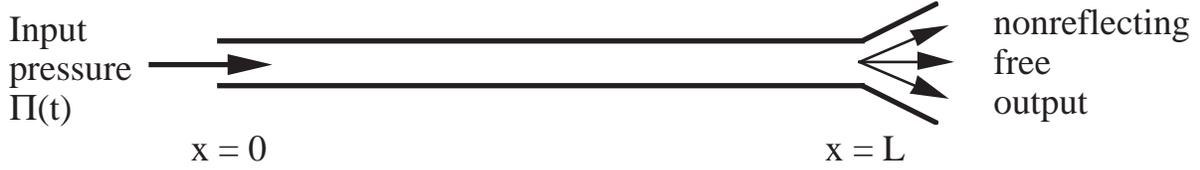
\centerline {\rm  {\bf Figure 3.1.}	\quad    Long unidimensional pipe for the modelling
of a trombone. }

\bigskip 
\noindent {\smcaps 3.2 } $ \,  $ { \bf  Linearization of the gas dynamics
equations.  }

\noindent   $\bullet \qquad \,\,\,\,\, $ 
We study a perfect gas subjected  to a motion with variable velocity in space and
time. We have noticed   that  the primitive  unknowns of this problem are  the
scalar fields that  characterize the thermodynamics  of the gas, {\it i.e.}  density 
$\, \rho ,\, $ internal energy~$\, e ,\,$ temperature $\, T ,\,$  and
pressure~$\,p.\,$ In what follows, we suppose that the gas is a polytropic perfect
gas~; it has constant specific heats at constant volume $\, C_v \,$ and at constant
pressure  $\, C_p .\,$ These two quantities do not depend on any thermodynamic
variable like temperature or pressure~; we denote by  $\, \gamma \,$ their ratio~: 

\noindent  (3.2.1) $\qquad   \displaystyle 
\gamma \,\,= \,\, {{C_p}\over{C_v}} \,\, (= \,\, {\rm constant}) \,. \,$ 

\noindent
We suppose that  the gas satisfies the law of perfect gas that can be written
with the following form~: 

\noindent  (3.2.2) $\qquad \displaystyle 
p \,\,= \,\, (\gamma-1) \, \rho \, e \,.\, $

\noindent
As usual, internal energy and temperature are linked together  by the Joule-Thomson
relation~: 

\noindent  (3.2.3) $\qquad \displaystyle 
e \,\,= \,\, C_v \, T \,. \,$

%%%%%%%%%%%%%%%%%%%%%%%%%%%%%%%%%       modif janvier 2011 
%%  \titredroite={\pecaps    Gas dynamics with  the Roe method } 

\smallskip \noindent   $\bullet \qquad \,\,\, $ 
In the formalism proposed by Euler during the 18$^{\rm th}$ century, the motion is
described with the help of an unknown vector field  $\, u \,$ which is a function of
space $ \,x \,$ and time $ \,t \,$~: 

\noindent  (3.2.4) $\qquad \displaystyle 
u \,\,= \,\, u(x,\,t) \,.\,$

\noindent
In the following, we will suppose that space $ \,x \,$ has only one dimemsion $\, (x
\in \R).$ We have four unknown functions (density, velocity, pressure and internal
energy) linked together by the state law (3.2.2). In consequence, we need three
complementary equations in order to define a unique solution of the problem. The
general laws of Physics assume that mass, momentum and total energy are conserved
quantities, at least in the context of classical physics associated to the paradigm
of  invariance for the  Galileo group of space-time transformations (see {\it e.g.} Landau
and Lifchitz [LL54]). When we write the conservation of mass, momentum and energy
inside an infinitesimal volume  $\, {\rm d}x \,$ advected with celerity   $\, u(x,\,t)
,\,$ which is exactly the mean velocity of particules that compose the gas, it is
classical [LL54] to write the fundamental  conservation laws of Physics with the help
of divergence operators~: 

\noindent  (3.2.5) $\qquad \displaystyle 
{{\partial \rho}\over{\partial t}} \,\,+\,\, {{\partial}\over{\partial x}}\,
\bigl( \rho \, u \bigr) \,\,= \,\, 0  \,$

\noindent  (3.2.6) $\qquad \displaystyle 
{{\partial  }\over{\partial t}} \bigl( \rho \, u \bigr)  \,\,+\,\,
{{\partial}\over{\partial x}}\, \bigl( \rho \, u^2 \,+\, p \bigr) \,\,= \,\, 0  \,$

\noindent  (3.2.7) $\qquad \displaystyle 
{{\partial}\over{\partial t}} \Bigl( {{1}\over{2}} \rho \, u^2 \,+\, \rho \,e \, 
\Bigr)  \,\,+\,\, {{\partial}\over{\partial x}}\,\Bigl( \, \bigl( {{1}\over{2}} \rho
\, u^2 \,+\, {{p}\over{\gamma - 1 }} \bigr) \, u \,
+\, p \, u \,  \Bigr) \,\,= \,\, 0  \,. \,$

\smallskip 
%%% \vfill \eject %%%%%%%%%  janvier 2011
\noindent   $\bullet \qquad \,\,\, $ 
We introduce the specific total energy  $\, E \,$ by unity of volume

\noindent  (3.2.8) $\qquad  
E \,\,= \,\, {{1}\over{2}}  u^2 \,+\, e \,\,, \,$ 

\noindent
the sound celerity $\,\, c \,\, $ following the classical expression~: 

\noindent  (3.2.9) $\qquad  \displaystyle 
c \,\,= \,\, \sqrt{{\gamma \, p} \over{\rho}} \, ,\,  $ 

\noindent
and total enthalpy $\, H \, $ defined according to 

\noindent  (3.2.10) $\qquad  
H \,\,\equiv \,\, E \,+\, {{p}\over{\rho}} \,\,= \,\, 
{{1}\over{2}}  u^2 \,+\, {{1}\over{\gamma - 1}} \, c^2  \,. \,$ 

\noindent
The vector  $\, W \,$ is therefore composed by the ``conservative variables'' or more
precisely by the ``conserved variables''~: 

\noindent  (3.2.11) $\qquad  \displaystyle
W \,\, = \,\, \bigl( \, \rho \,,\, \rho \, u \,,\, \rho \, E \, \bigr) ^{\displaystyle
\rm t} \,\, \equiv \,\, \bigl( \, \rho \,,\, q \,,\, \epsilon \, \bigr) ^{\displaystyle
\rm t} \,. \,$ 

\noindent 
The conservation laws (3.2.5)-(3.2.7) take the following general  form of a  system of
conservation laws~: 

\noindent  (3.2.12) $\qquad \displaystyle 
{{\partial W}\over{\partial t}} \,\,+\,\, {{\partial}\over{\partial x}}\, 
F(W)  \,\,= \,\, 0 \,$

\noindent
where the flux vector $\, W \longmapsto F(W) \,$ satisfies the following algebraic
expression~: 

\noindent  (3.2.13) $\qquad \displaystyle 
F(W)  \,\, = \,\, \bigl( \, \rho \,u \,,\, \rho \, u^2 \,+\,p  \,,\, \rho \, u \, H
\, \bigr) ^{\displaystyle \rm t} \, \,$

\noindent
that can be explicited as a {\it true} function of state vector $\,  W ,\,$ on one
hand  with the pressure law $\,\, P(W) \,\,$ computed with (3.2.2), (3.2.8)  and
(3.2.11)~: 

\noindent  (3.2.14) $\qquad \displaystyle 
P(W) \,\, = \,\, (\gamma \!-\!1 ) \, \biggl( \epsilon - {{q^2}\over{2 \, \rho}}
\biggr)  \,$

\noindent
and on the other hand with an explicit use of the conserved variables $\, \rho ,\,$
$\, q \,$ and $\, \epsilon .\,$ We obtain~: 

\noindent  (3.2.15) $\qquad \displaystyle 
F(W)  \,\, = \,\, \Bigl( \, q \,,\, {{q^2}\over{\rho}} + P(W) \,,\, 
{{q\, \epsilon }\over{\rho}} + P(W)  \, {{q}\over{\rho}} \, \Bigr) \,. \,$ 

\bigskip  \noindent  {\bf Proposition 3.1.  $\quad$  Jacobian matrix of gas
dynamics. }

\noindent 
 $\bullet \quad $
The Jacobian matrix $\,\, {\rm d}F(W) \,$ of the flux function $\,\, W \longmapsto
F(W) \,\,$ for the Euler equations of the gas dynamics admits the following
expression~: 

\noindent  (3.2.16) $\qquad \displaystyle 
{\rm d}F(W) \,\, = \,\, \pmatrix { 0 & 1 & 0 \cr (\gamma \!-\! 1)\,  H - u^2 - c^2 &
(3-\gamma) \, u & \gamma - 1  \cr (\gamma \!-\! 2) \, u \, H - u \,c^2 & H - (\gamma
\!-\! 1) \, u^2 & \gamma u \cr  } \,. \,$ 

\noindent   $\bullet \quad $
The matrix $\,\, {\rm d}F(W) \,$ is diagonalizable~; the eigenvalues $\,\,
\lambda_{j}(W) \,$ satisfy the relations 

\noindent  (3.2.17) $\qquad \displaystyle 
\lambda_{1}(W) \equiv u - c \quad < \quad  \lambda_{2}(W) \equiv u  \quad < \quad 
\lambda_{3}(W) \equiv u + c \,. \,$ 

\noindent 
and the associated eigenvectors $\,\, r_{j}(W) \,\,$  are proportional to the
following ones~: 

\noindent  (3.2.18) $\quad \displaystyle 
r_{1}(W)=  \pmatrix{ 1 \cr u - c \cr H - u \, c \cr } \,, \,\,\,\,  
r_{2}(W)=  \pmatrix{ 1 \cr u  \cr{1\over2}u^2 } \,, \,\,\,\,  
r_{3}(W)=  \pmatrix{ 1 \cr u + c \cr H + u \, c \cr } \,. \,$

\smallskip \noindent   $\bullet \qquad \,\,\, $ 
We first differentiate the pressure law $\,\, W \longmapsto P(W) \,\,$ given in
(3.2.14)~:

\noindent  (3.2.19) $\quad \displaystyle 
{{\partial P}\over{\partial \rho}} \,=\, {{\gamma \!-\!1}\over{2}} \, u^2 \,=\,
(\gamma \!-\!1)\, H - c^2 \,,\quad  {{\partial P}\over{\partial q}} \,=\, -(\gamma
\!-\!1)\, u  \,,\quad  {{\partial P}\over{\partial \epsilon}} \,=\,(\gamma \!-\!1) \, $

\noindent
and the second row of the matrix (3.2.16) is a direct consequence of the  relations
$  \displaystyle \,\,\,  {{\partial }\over{\partial \rho}} 
\bigl( {{q^2}\over{\rho}} \bigr) = -u^2 \,\,$ and 
$\,\,  \displaystyle  {{\partial }\over{\partial q}} 
\Bigl( {{q^2}\over{\rho}} \Bigr) = 2 \, u \,.\, $ 

\smallskip \noindent   $\bullet \qquad \,\,\, $ 
The calculus of the third row of matrix in (3.2.16) demands first  evaluation  of the
gradient of $\,\, \rho \, u \, E \,=\, u \, \epsilon  \,\,$ relatively to the state $\,
W .\,$ We get 

\noindent  (3.2.20) $\qquad \displaystyle 
{{\partial}\over{\partial \rho}} \Bigl( {{q\, \epsilon }\over{\rho}} \Bigr) \,= \, -u\,
E  \,,\quad  {{\partial}\over{\partial q}} \Bigl( {{q\, \epsilon }\over{\rho}} \Bigr)
\,= \, E \,,\quad {{\partial }\over{\partial \epsilon}}  \Bigl( {{q\, \epsilon
}\over{\rho}} \Bigr) \,=\,u \,. \,$

\noindent 
We have also $\quad  {{\partial }\over{\partial W}} ( P \,u ) = 
{{\partial P}\over{\partial W}}\,u \,+\, p \, {{\partial}\over{\partial W}}
\bigl({{q}\over{\rho}} \bigr) \quad $ then we deduce from (3.2.19) and the following 
expressions for the gradient of velocity  $ \,\,\,  {{\partial }\over{\partial \rho}}
\bigl( {{q}\over{\rho}} \bigr) = -{{u}\over{\rho}}   \,\,$ and $\,\,   {{\partial
}\over{\partial q}} \Bigl( {{q}\over{\rho}} \Bigr) = {{1}\over{\rho}}\,:\, $ 

\setbox21=\hbox {$\displaystyle 
{{\partial}\over{\partial \rho}} \Bigl( {{P \, q}\over{\rho}} \Bigr) \,= \,{{\gamma
\!-\!1}\over{2}}\,u^3 \, - \, {{u\,p}\over{\rho}}  \,,\,\,   $}
\setbox22=\hbox {$\displaystyle  
{{\partial}\over{\partial q}} \Bigl( {{P \, q}\over{\rho}} \Bigr) \,= \,-(\gamma
\!-\!1)\,u^2 \,+\, {{p}\over{\rho}}  \,\,,\quad {{\partial }\over{\partial \epsilon}}
\Bigl( {{P \, q}\over{\rho}} \Bigr) \,= \,(\gamma \!-\!1)\,u \,.\,$  }
\setbox30= \vbox {\halign{#&# \cr \box21 \cr \box22    \cr   }}
\setbox31= \hbox{ $\vcenter {\box30} $}
\setbox44=\hbox{\noindent   (3.2.21)  $\displaystyle  \qquad \left\{ \box31
\right.$}  
\noindent $ \box44 $

\noindent
We add the relations (3.2.20) and (3.2.21)~; then the third row of matrix (3.2.16)
admits the following expression~:   $  \,\,\, 
\bigl( \, {{\gamma \!-\!1}\over{2}}\,u^3 \,- \, u \, H \,\,,\,\,\,  H \,-\,
(\gamma\!-\!1)\,u^2 \,\,,\,\,\, \gamma \,u \, \bigr) \,\,\, $ and this result is
exactly the third row of the right hand side of (3.2.16) when we take into account
the relation (3.2.10) between $ \,H ,\,$ $\, u^2 \,$ and $\, c^2 .\,$  The relations
(3.2.17) and (3.2.18) are elementary to satisfy~; they express simply the three
relations~: 

\noindent  (3.2.22) $\qquad \displaystyle 
{\rm d}F(W) \, {\scriptstyle \bullet} \, r_{j}(W) \,\,= \,\, \lambda_{j}(W) \, 
r_{j}(W) \,\,,\quad j=1,\,2,\,3 \,$ 
 
\noindent 
and Proposition 3.1 is established.  $ \hfill \square \kern0.1mm $

\smallskip \noindent   $\bullet \qquad \,\,\, $ 
We keep into memory the following expression of the Jacobian matrix $\,\, {\rm d}F(W)
\,: \, $

\noindent  (3.2.23) $\qquad \displaystyle 
{\rm d}F(W) \,\, = \,\, \pmatrix { 0 & 1 & 0 \cr {{\gamma \!-\! 3}\over{2}} \,  u^2  &
(3-\gamma) \, u & \gamma - 1  \cr{{\gamma \!-\! 1}\over{2}}  \, u^3 \, - u \, H & H -
(\gamma \!-\! 1) \, u^2 & \gamma u \cr  } \, \,$ 

\noindent 
that needs only the datum of velocity $\, u \,$ and total enthalpy $\, H \,$ of the
state $\, W .\,$

 \bigskip  \noindent {\smcaps 3.3 } $ \,  $ { \bf   Roe matrix. }

\noindent   $\bullet \qquad \,\,\,\,\, $ 
We consider two states $\,W_{\rm left} \equiv W_l \,$ and $\,W_{\rm right} \equiv W_r
\,$ relatively to the gas dynamics, {\it i.e.} they both belong to space $\, \R^3 \,$ and
have an expression of the form (3.2.11).  By {\bf definition}, a Roe matrix $\,\, A(
W_l ,\,  W_r) \,\,$ between these two states is a 3 by 3 matrix that satisfy the
three  following properties~: 

\noindent  (3.3.1) $\qquad \displaystyle 
 A( W_l ,\,  W_r) \,\,$ is a diagonalizable matrix on the field $\, \R \,$ of real
numbers

\noindent  (3.3.2) $\qquad \displaystyle 
 A( W,\,W) \,\,= \,\, {\rm d}F(W) \,$ 

\noindent  (3.3.3) $\qquad \displaystyle 
F(W_r) - F(W_l) \,\,= \,\, A( W_l ,\,  W_r) \, {\scriptstyle \bullet} \, (W_r -  W_l)
\, . \,$

\noindent 
In his original article, P. Roe   [Roe81] has proposed a very simple algebraic way to
construct a Roe matrix for the dynamics of polytropic gas. We propose it in the
following Proposition. 

\smallskip   \noindent  {\bf Proposition 3.2.  $\quad$  Algebraic construction of a
Roe matrix [Roe81]. }

\noindent 
Let $\,W_l \,$ and $\, W_r \,$ be two states for gas dynamics, defined by their
densities $\, \rho_l \,$ and $\, \rho_r ,\,$ their velocities $\, u_l \,$ and $\, u_r
\,$ and their total enthalpies $\, H_l \,$ and $\, H_r .\,$ We introduce an {\bf
intermediate state} $\,\, W^*(W_l ,\,  W_r) \,\,$ by its density $\, \rho^* ,\,$ its
velocity $\, u^* \,$ and its total enthalpy $\, H^* \,$  according to  the following
relations~: 

\noindent  (3.3.4) $\qquad \displaystyle 
\rho^* \,\,= \,\, \sqrt{\rho_l \, \rho_r} \,$ 

\noindent  (3.3.5) $\qquad \displaystyle
u^* \,\,= \,\, {{\sqrt{\rho_l} \, u_l \,+\, \sqrt{\rho_r} \, u_r}\over{ \sqrt{\rho_l}
 \,+\, \sqrt{\rho_r} }} \,$ 
 
\noindent  (3.3.6) $\qquad \displaystyle
H^* \,\,= \,\, {{\sqrt{\rho_l} \, H_l \,+\, \sqrt{\rho_r} \, H_r}\over{ \sqrt{\rho_l}
 \,+\, \sqrt{\rho_r} }} \,. \,$ 

\noindent 
Then the matrix $\,\, A( W_l ,\,  W_r) \, \,$ defined as the Jacobian matrix of the
flux for the intermediate state $\, \,  W^*(W_l ,\,  W_r) ,\,\,$ {\it i.e.}

\noindent  (3.3.7) $\qquad \displaystyle
 A( W_l ,\,  W_r) \,\,= \,\, {\rm d}F \bigl(  W^*(W_l ,\,  W_r)  \bigr) \,$ 

\noindent 
is a Roe matrix.

\smallskip \noindent   $\bullet \qquad \,\,\, $ 
Due to the expression (3.2.23) of the Jacobian matrix of gas dynamics, we remark that
the formula (3.3.4) giving the density $\, \rho^* \,$ is not necessary for the
determination of the matrix $ \,\, {\rm d}F (  W^*(W_l ,\,  W_r)  ) \,\,$ and an
entire family of states $\,  W^*(W_l ,\,  W_r) \,$ define a Roe matrix according to
the  relations (3.3.5), (3.3.6) and  (3.3.7).  Nevertheless, we keep this definition of
density  $\, \rho^* \,$  by convenience and  simplicity for future algebraic
expressions. The proof of Proposition 3.2 needs some algebraic developments. We begin
by the following technical lemma. 

\smallskip   
% \vfill \eject 

\noindent  {\bf Proposition 3.3.  }

\noindent 
Under the hypotheses of Proposition 3.2, we have the following relations~: 

\noindent  (3.3.8) $\qquad \displaystyle
(u^*)^2 \,(\rho_r - \rho_l) \,- \, 2 \, u^* \, (\rho_r \, u_r - \rho_l \,u_l) \, +\,
(\rho_r \, u_r^2 - \rho_l \,u_l^2) \,\,=\,\, 0 \,$

\setbox21=\hbox {$\displaystyle 
-u^* \, H^* \, (\rho_r - \rho_l) \,+\, H^* \,(\rho_r \, u_r - \rho_l \,u_l) \,+\,u^*
\, (\rho_r \, H_r - \rho_l \,H_l) \,\,=\, \, $}
\setbox22=\hbox {$\displaystyle  \qquad  \qquad 
\,\, =\,\,   \rho_r \, u_r \, H_r - \rho_l \,u_l \,H_l \,. \,$  }
\setbox30= \vbox {\halign{#&# \cr \box21 \cr \box22    \cr   }}
\setbox31= \hbox{ $\vcenter {\box30} $}
\setbox44=\hbox{\noindent   (3.3.9)  $\displaystyle  \qquad \left\{ \box31
\right.$}  
\noindent $ \box44 $

\smallskip 
\vfill \eject    %%%%%  janvier 2011  
\noindent   $\bullet \qquad \,\,\, $ 
We first evaluate the left hand side of relation (3.3.8)~: 

\noindent $ \displaystyle
(u^*)^2 \,(\rho_r - \rho_l) \,- \, 2 \, u^* \, (\rho_r \, u_r - \rho_l \,u_l) \, +\,
(\rho_r \, u_r^2 - \rho_l \,u_l^2) \,\,=\,\, $ 

\noindent $ \displaystyle \qquad \qquad =\,\, 
u^* (\sqrt{\rho_r} - \sqrt{\rho_l}) \, (\sqrt{\rho_r} \, u_r + \sqrt{\rho_l} \, u_l )
\,-\, 2 \, u^*  (\rho_r \, u_r - \rho_l \, u_l ) \,+\, (\rho_r \, u_r^2 - \rho_l
\,u_l^2) \,$

\noindent $ \displaystyle \qquad \qquad =\,\, 
u^* \, \bigl(  \sqrt{\rho_l} \, ( \sqrt{\rho_l} +  \sqrt{\rho_r}) \, u_l \,-\, 
\sqrt{\rho_r} \, ( \sqrt{\rho_l} +  \sqrt{\rho_r}) \, u_r \bigr) \,+\, (\rho_r \, 
u_r^2 - \rho_l \,u_l^2) \,$

\noindent $ \displaystyle \qquad \qquad =\,\, 
 ( \sqrt{\rho_l} \, u_l +  \sqrt{\rho_r} \, u_r )  \,  ( \sqrt{\rho_l} \, u_l - 
\sqrt{\rho_r} \, u_r ) \,+\, (\rho_r \,  u_r^2 - \rho_l \,u_l^2) \,$

\noindent $ \displaystyle \qquad \qquad =\,\,  0 \qquad \qquad  $ and the relation
(3.3.8) is established. 

\smallskip \noindent   $\bullet \qquad \,\,\, $ 
We work on the left hand side of (3.3.9)  as follows~: 

\noindent $ \displaystyle
-u^* \, H^* \, (\rho_r - \rho_l) \,+\, H^* \,(\rho_r \, u_r - \rho_l \,u_l) \,+\,u^*
\, (\rho_r \, H_r - \rho_l \,H_l) \,\,=\, \, $ 

\noindent $ \displaystyle  =\,\, 
-u^* ( \sqrt{\rho_r} -  \sqrt{\rho_l}) \,(  \sqrt{\rho_l}\,H_l +  \sqrt{\rho_r} \,
H_r) \,+\, H^* \, (\rho_r \, u_r - \rho_l \,u_l) \,+\, u^* \, (\rho_r \, H_r - \rho_l
\,H_l) \,$ 

\noindent $ \displaystyle  =\,\, 
\sqrt{\rho_l \, \rho_r} \, u^* \, (H_r - H_l) \,+\, H^* \, (\rho_r \, u_r - \rho_l
\,u_l) \,$

\noindent $ \displaystyle  =\,\, 
{{ \sqrt{\rho_l}\, \sqrt{\rho_r} \, (\sqrt{\rho_l}\,u_l + \sqrt{\rho_r}\,u_r) \,( -H_l
+ H_r)  + (-\rho_l \, u_l + \rho_r \, u_r) \,(\sqrt{\rho_l}\,H_l
+ \sqrt{\rho_r}\,H_r)  } \over{ \sqrt{\rho_l}+\sqrt{\rho_r}}} \,$

\noindent $ \displaystyle  =\,\, 
{{1}\over{\sqrt{\rho_l}+\sqrt{\rho_r}}} \, \bigl[ \, -\rho_l \,(
\sqrt{\rho_l}+\sqrt{\rho_r})\, u_l \, H_l \,+\, \rho_r \,(
\sqrt{\rho_l}+\sqrt{\rho_r})\, u_r \, H_r \, \bigr] \,$ 

\noindent $ \displaystyle  =\,\, 
\rho_r \, u_r \, H_r - \rho_l \,u_l \,H_l \,$ 

\noindent
and the proposition 3.3 is established.  $ \hfill \square \kern0.1mm $

\smallskip \smallskip \noindent   $\bullet \qquad \,\,\, $ 
The {\bf proof of Proposition 3.2} consists in satisfying the three hypotheses that define
a Roe matrix. First, due to the fact that the relation  (3.3.7) defines the matrix
$\,\, A(W_l,\,W_r) \,\,$ as a Jacobian of some state, this matrix is diagonalizable
with real elements due to the result of Proposition 3.1 and the first property
(3.3.1) is satisfied.  The second property (3.3.2) is a simple consequence of the fact
that if $\,\, W_l \,=\, W_r \,=\, W ,\, \,$ then we have  from the relations (3.3.4)
to (3.3.6)~: $\,\, W^*(W_l ,\,W_r) \, = \,  W \, \,$ and the property results
from (3.3.7). 

\smallskip \noindent   $\bullet \qquad \,\,\, $ 
The third property (3.3.3) needs more work. We  remark that the first row of
this matricial relation is clear. For the second row, we have~: 

\noindent 
Second row of matrix $\,\, A(W_l ,\,W_r) \, {\scriptstyle \bullet} \, (W_r - W_l)
\,\,= \,$

\noindent  $ \displaystyle  =\,\, 
{{\gamma \!-\!3}\over{2}} \, (u^*)^2 \,(\rho_r - \rho_l) \,+\, (3 \!-\! \gamma) \,
u^* \,(\rho_r \, u_r - \rho_l \, u_l ) \,+\, (\gamma \!-\!1) \,(\rho_r \, E_r - \rho_l
\, E_l ) \,$

\noindent  $ \displaystyle  =\,\, 
{{\gamma \!-\!3}\over{2}} \, \bigl[ (u^*)^2 \,(\rho_r - \rho_l) \,- 2 \, u^*
\, (\rho_r \, u_r - \rho_l \, u_l )  \bigr] \,+\,  {{\gamma \!-\!1}\over{2}} \, 
(\rho_r \, u_r^2 - \rho_l \, u_l^2 )  \,+\, (p_r - p_l) \,$

\noindent  $ \displaystyle  =\,\, 
(\rho_r \, u_r^2 - \rho_l \, u_l^2 )  \,+\, (p_r - p_l) \qquad \qquad \qquad \qquad 
\qquad \qquad  \,   \, $ due to (3.3.8) 

\noindent  $ \displaystyle  =\,\, $
second row of the flux difference $\,\, F(W_r) - F(W_l) \,. $ 

\smallskip \noindent   $\bullet \qquad \,\,\, $ 
We have also, in consequence of  (3.2.23), 

\noindent 
Third  row of matrix $\,\, A(W_l ,\,W_r) \, {\scriptstyle \bullet} \, (W_r - W_l)
\,\,= \,$

\noindent  $ \displaystyle  =\,\, 
u^* \, \bigl( {{\gamma \!-\!1}\over{2}} \, (u^*)^2 - H^* \bigr) \, (\rho_r - \rho_l)
\,+\, (H^* - (\gamma \!-\!1) \, (u^*)^2 ) \, (\rho_r \, u_r  - \rho_l \, u_l ) \,+\,$

\noindent  $ \displaystyle  \qquad \qquad \qquad \,+\,
\gamma \, u^* \, (\rho_r \, E_r  - \rho_l \, E_l ) \,$

\noindent  $ \displaystyle  =\,\, 
{{\gamma \!-\!1}\over{2}} \,u^* \, \bigl[ (u^*)^2 \,  (\rho_r - \rho_l) \, - \, 2 \,
u^* \,  (\rho_r \, u_r  - \rho_l \, u_l ) \bigr] \,+ \,$

\noindent  $ \displaystyle  \qquad \qquad \qquad \,+\,
 \bigl[ -H^* \, u^* \, (\rho_r - \rho_l) \,+\, H^* \,  (\rho_r \, u_r  - \rho_l \, 
 u_l ) \bigr] \,+\, \gamma \, u^* \, (\rho_r \, E_r  - \rho_l \, E_l ) \,$ 

\noindent  $ \displaystyle  =\,\, 
{{\gamma \!-\!1}\over{2}} \,u^* \, \bigl[ (u^*)^2 \,  (\rho_r - \rho_l) \, - \, 2 \,
u^* \,  (\rho_r \, u_r  - \rho_l \, u_l ) \bigr]  \,-\, u^* \,  (\rho_r \, H_r  -
\rho_l \, H_l ) \,+\, $

\noindent  $ \displaystyle  \qquad \qquad \qquad   \,+\, 
(\rho_r \, u_r  \, H_r  - \rho_l \,u_l \,   H_l ) \,+\, \gamma \, u^* \, (\rho_r \,
E_r  - \rho_l \, E_l ) \qquad \qquad $ due to (3.3.9) 

\noindent  $ \displaystyle  =\,\, 
{{\gamma \!-\!1}\over{2}} \,u^* \, \bigl[ \, (u^*)^2 \,  (\rho_r - \rho_l) \, - \, 2 \,
u^* \,  (\rho_r \, u_r  - \rho_l \, u_l ) \,+\, \rho_r \, u_r^2  - \rho_l \, u_l^2 \,
\bigr]  \,+\, $

\noindent  $ \displaystyle  \qquad \qquad \qquad   \,+\, 
 u^* \,  (- \gamma \, \rho_r \, e_r  + \gamma \, \rho_l \, e_l \, + \gamma \, \rho_r
\, e_r - \gamma \, \rho_l \, e_l ) \,+\, \rho_r \, u_r \, H_r - \rho_l \, u_l \, H_l
\,$

\noindent  $ \displaystyle  =\,\, 
\rho_r \, u_r \, H_r - \rho_l \, u_l \, H_l  \qquad \qquad \qquad \qquad \qquad \qquad  
 \qquad \qquad \qquad \qquad \,  $  due to (3.3.8)

\noindent  $ \displaystyle  =\,\, $
third row of the flux difference $\,\, F(W_r) - F(W_l) \,$  

\noindent
in the view of relation  (3.2.13). The proposition 3.2 is established.  $ \hfill
\square \kern0.1mm $

\bigskip \noindent {\smcaps 3.4 } $ \,  $ { \bf   Roe flux. }

\noindent   $\bullet \qquad \,\,\,\,\, $ 
The principal interest of the Roe matrix is to be able to use all what has been 
developed for  {\bf linear} hyperbolic systems in Section 2. In particular, the
following  linear hyperbolic system defined with a given Roe matrix $\,\,
A(W_l,\,W_r) \,\,$ 

\noindent  (3.4.1) $\qquad \displaystyle
{{\partial W}\over{\partial t}} \,\, + \,\,  A(W_l,\,W_r) \, {\scriptstyle \bullet} \,
{{\partial W}\over{\partial x}} \,\,= \,\, 0 \,$ 

\noindent 
can be treated with the upwind scheme defined at  proposition 2.7. We obtain by doing
this the following

\smallskip   \noindent  {\bf Proposition 3.4.  $\quad$  Three formulae for a flux. }

\noindent 
 $\bullet \quad $
Let $\,\, W_l \,\,$ and $\,\, W_r \,\,$ be two fluid states and $\,\, W^* \,\, $ the
intermediate state defined by the relations  (3.3.4) to (3.3.6). The sound celerity
$\,\, c^* \,\,$ of state $\, W^* \,$ is defined with the help of relation (3.2.10),
{\it i.e.}  

\noindent  (3.4.2) $\qquad \displaystyle 
c^* \,\, = \,\, \sqrt{(\gamma \!-\!1) \, \Bigl(\, H^* - {{(u^*)^2}\over{2}} \,\Bigr) }
\,,\,$ 

\noindent 
and the eigenvalues $\,\, \lambda_{j}^* \,\, $ of the Roe matrix $\, \,
A(W_l,\,W_r) \,\equiv\,  {\rm d}F (  W^*(W_l ,\,  W_r)  ) \, \,$ are given by a
relation analogous to (3.2.17). 

\noindent  (3.4.3) $\qquad \displaystyle 
\lambda_{1}^* \equiv u^* - c^* \quad < \quad  \lambda_{2}^* \equiv u^* 
\quad < \quad  \lambda_{3}^* \equiv u ^* + c^* \,. \,$ 

\noindent 
The associated eigenvectors $\,\, r_{j}^* \,\equiv\,  r_{j}(W^*) \,\,$  are
proportional to the following ones~: 

\noindent  (3.4.4) $\quad \displaystyle 
r_{1}^* =  \pmatrix{ 1 \cr u^* - c^* \cr H^* - u^* \, c^* \cr } \,, \,\,\,\,  
r_{2}^* =  \pmatrix{ 1 \cr u^*  \cr{1\over2}(u^*)^2 } \,, \,\,\,\,  
r_{3}^* =  \pmatrix{ 1 \cr u^* + c^* \cr H^* + u^* \, c^* \cr } \,. \,$

\noindent  $\bullet \quad $
We introduce the decomposition of vector $\,\, W_r - W_l \,\, $ in the basis $\,\,
r_{j}^* \, : \,$

\noindent  (3.4.5) $\quad  
 W_r - W_l \,\, = \,\, \sum_{j=1}^{j=3} \, \alpha_{j} \,  r_{j}^* \,. \,$

\noindent 
The three following relations define a unique numerical flux $\,\, \Phi(W_l,\,W_r)
\,\,$ named the {\bf Roe flux}  between the two states $\, W_l \,$ and $ \, W_r \,:
\,$ 

\noindent  (3.4.6) $\qquad  
\Phi(W_l,\,W_r) \,\,=\,\, F(W_l) \,+\, \sum_{j=1}^{j=3} \, (\lambda_{j}^*)^-\,\, 
\alpha_{j}\,\,   r_{j}^* \,$
 
\noindent  (3.4.7) $\qquad  
\Phi(W_l,\,W_r) \,\,=\,\, F(W_r) \,-\, \sum_{j=1}^{j=3} \, (\lambda_{j}^*)^+\,\,  
\alpha_{j}\,\,   r_{j}^* \,$
 
\noindent   (3.4.8) $\qquad  
\Phi(W_l,\,W_r) \,\,=\,\, {1\over2} \, \bigl(  F(W_l) \,+\,F(W_r) \bigr) \,-\,
{1\over2} \,\abs{A(W_l,\,W_r)}  \, {\scriptstyle \bullet} \, (W_r \,-\, W_l \bigr) \,.
\,$

\smallskip \noindent   $\bullet \qquad \,\,\, $ 
The first non-obvious point is to verify that the relation (3.4.2) defines a real
number $\, c^* .\,$ We have 

\noindent  $ \displaystyle 
H^* \,- \, {{(u^*)^2}\over{2}} \,\,=\,\, 
  {{\sqrt{\rho_l} \, H_l \,+\, \sqrt{\rho_r} \, H_r}\over{ \sqrt{\rho_l}
 \,+\, \sqrt{\rho_r} }}\,-\, {1\over2} \, \Bigl( \, {{\sqrt{\rho_l} \, u_l \,+\, \sqrt{\rho_r} \, u_r}\over{ \sqrt{\rho_l}
 \,+\, \sqrt{\rho_r} }}  \, \Bigr)^2 \,\,= \, \,$

\noindent  $ \displaystyle  = \,\,
{{1}\over{(\sqrt{\rho_l} + \sqrt{\rho_r})^2}} \,\Bigl[ \,(\sqrt{\rho_l} +
\sqrt{\rho_r}) \, \Bigl( \sqrt{\rho_l}\,\bigl( {1\over2}\,
u_l^2 + {{1}\over{\gamma \!-\!1}}\, c_l^2\bigr) \,+\, \sqrt{\rho_r}\,\bigl( {1\over2}\,
u_r^2 + {{1}\over{\gamma \!-\!1}}\, c_r^2\bigr) \, \Bigr) \,\, \,$

\noindent  $ \displaystyle  \qquad \qquad \qquad \qquad \qquad   \,-\, 
{1\over2} \, \Bigl( \rho_l \, u_l^2 \,+\, 2 \, \rho^* \, u_l \, u_r \,+\, \rho_r \,
u_r^2  \Bigr) \, \Bigr] \,$

\noindent  $ \displaystyle  = \,\,
{{1}\over{(\sqrt{\rho_l} + \sqrt{\rho_r})^2}} \,\Bigl[ \,{1\over2} \, \bigl( 
\rho^* \, u_l^2  - 2 \, \rho^* \, u_l \, u_r + \rho^* \, u_r^2 \bigr) \,\,+\,\,
{{\rho_l + \rho^*}\over{\gamma \!-\!1}} \, c_l^2 +  {{\rho^* + \rho_r}\over{\gamma
\!-\!1}} \, c_r^2 \, \Bigr] \,$ 

\noindent  $ \displaystyle  = \,\,
{{1}\over{(\sqrt{\rho_l} + \sqrt{\rho_r})^2}} \,\Bigl[ \,{1\over2} \,\rho^* \, (u_r
- u_l)^2 \,\,+\,\, {{\rho_l + \rho^*}\over{\gamma \!-\!1}} \, c_l^2 +  {{\rho^* +
\rho_r}\over{\gamma \!-\!1}} \, c_r^2 \, \Bigr] \qquad > \,\,  0 \,$ 

\noindent
and 

\noindent  (3.4.9) $\qquad \displaystyle 
c^* \,\, = \,\, {\sqrt{ \,{{\gamma \!-\!1} \over{2}} \,\rho^* \, (u_r
- u_l)^2 \,\,+\,\, (\rho_l + \rho^*) \, c_l^2 +  (\rho^* + \rho_r) \, c_r^2
\,}\over{\sqrt{\rho_l} + \sqrt{\rho_r}}} \,. \,$

\smallskip \noindent   $\bullet \qquad \,\,\, $ 
We make the difference between the right hand sides of (3.4.6) and (3.4.7). We get~:

\noindent  $   
F(W_r) - F(W_l) \,\,-\,\, \sum_{j=1}^{j=3} \, \bigl(  (\lambda_{j}^*)^+ + 
(\lambda_{j}^*)^- \bigr) \,\, \alpha_{j} \,\, r_{j}^* \,\,= \,$ 

\noindent  $    = \,\,
A(W_l,\,W_r) \, { \scriptstyle \bullet} \,  (W_r - W_l) \,\,- \,\, \sum_{j=1}^{j=3}
\,  \lambda_{j}^*\,\, \alpha_{j} \,\, r_{j}^* \quad \,\, $  \qquad due to
(3.3.3) and (2.7.11) 

\noindent  $    = \,\,
A(W_l,\,W_r) \, { \scriptstyle \bullet} \, \bigl( 
\sum_{j=1}^{j=3} \, \alpha_{j} \,  r_{j}^*  \bigr) \,\,- \,\, \sum_{j=1}^{j=3}
\,  \lambda_{j}^*\,\, \alpha_{j} \,\, r_{j}^* \, \hfill $  due to (3.4.5) 

\noindent  $    = \,\, 0 \,  \hfill $  
because $ \,\,\, A(W_l,\,W_r) \, { \scriptstyle \bullet} \,r_{j}^* 
\,=\, \lambda_{j}^*\,r_{j}^* \,\,$ for each integer $\,j .\, $  

 \noindent
The proof of relation (3.4.8) is obtained by taking
the half sum of  (3.4.6) and (3.4.7). It is  analogous to the one done for
Proposition 2.7. The proof of Proposition 3.4 is completed. 
~  $ \hfill \square \kern0.1mm $

\smallskip \noindent   $\bullet \qquad \,\,\,\,\, $ 
We make  explicit the parameters $\,\, \alpha_{j} \,\,$ introduced in relation
(3.4.5) in order to be complete for the implementation of the above formulae  on a
computer. 

\smallskip   \noindent  {\bf Proposition 3.5.  $\quad$  New acoustic impedance.}

\noindent 
With the notations introduced at Proposition 3.4, and denoting by $\, p_l \,$ and $\,
p_r \,$ the respective pressures of states $\, W_l \,$ and $ \, W_r ,\,$ we have the
following relations for the scalar components $\,\, \alpha_{j} \,\,$ of the state
difference $\,\, W_r - W_l \,\,$ in relation (3.4.5)~:

\noindent   (3.4.10) $\qquad  \displaystyle 
\alpha_{1} \,\,= \,\, {{1}\over{2 \, (c^*)^2}} \, \bigl[ \, (p_r - \rho^* \, c^* \,
u_r) \,- (p_l - \rho^* \, c^* \, u_l) \, \bigr] \,$ 

\noindent   (3.4.11) $\qquad  \displaystyle 
\alpha_{2} \,\,= \,\, -{{1}\over{ (c^*)^2}} \, \bigl[ \, (p_r - (c^*)^2 \,
\rho_r) \,- (p_l - (c^*)^2 \, \rho_l) \, \bigr] \,$ 

\noindent   (3.4.12) $\qquad  \displaystyle 
\alpha_{3} \,\,= \,\, {{1}\over{2 \, (c^*)^2}} \, \bigl[ \, (p_r + \rho^* \, c^* \,
u_r) \,- (p_l + \rho^* \, c^* \, u_l) \, \bigr] \, \,$ 

\noindent 
with  {\bf acoustic impedance} $\,\, \rho^* \, c^* \,\,$  that is nomore the one
$\,\, \rho_0 \, c_0 \,\,$  of a reference state as in traditional  acoustics  but an
impedance  associated with the Roe  intermediate state $\,\, W^*(W_l,\,W_r) \,\,$ of
relations (3.3.4) to (3.3.6).

\smallskip \noindent   $\bullet \qquad \,\,\,\,\, $ 
We have just to explicit the three components of the relation (3.4.5). It comes~: 

\noindent  $  \displaystyle \qquad \quad 
\pmatrix { \rho_r - \rho_l \cr \rho_r \, u_r - \rho_l \, u_l \cr  \rho_r \, E_r -
\rho_l \, E_l \cr } \,=\, \pmatrix { \alpha_{1} +  \alpha_{2} +  \alpha_{3} \cr 
 \alpha_{1} \, (u^* - c^*) \,+\,  \alpha_{2} \, u^* \,+\, \alpha_{3} \, (u^* +
c^*) \cr \alpha_{1} \, (H^* - u^* \, c^*)   \,+\,  \alpha_{2} {{(u^*)^2}\over{2}} 
\,+\,  \alpha_{3} \, (H^* + u^* \, c^*)   } \,$

\noindent   (3.4.13) $\qquad  \displaystyle 
\alpha_{1} +  \alpha_{2} +  \alpha_{3} \,\, = \,\,  \rho_r - \rho_l \,, \,$

\noindent
and we deduce after multiplying the  equation (3.4.13)  by $\,\, -u^* \,\,$ and adding
to the second equation of the above matrix equality~: 

\noindent  $  \displaystyle 
c^* \, (\alpha_{3} - \alpha_{1}) \,\,= \,\, \rho_r \, u_r \,-\,  \rho_l \, u_l 
\,-\, u^* \, (\rho_r - \rho_l) \,$ 

\noindent  $  \displaystyle \qquad \qquad \quad \,\,\, \, = \,\,
\rho_r \, u_r \,-\,  \rho_l \, u_l  \,-\, (\sqrt{\rho_r} - \sqrt{\rho_l}) \,
(\sqrt{\rho_l} \, u_l + \sqrt{\rho_r} \, u_r) \,\,,\, $ 

\noindent
then 

\noindent   (3.4.14) $\qquad  \displaystyle 
c^* \, (- \alpha_{1} + \alpha_{3} ) \,\,= \,\, \rho^* \, (u_r - u_l) \,.\,$ 

\smallskip \noindent   $\bullet \qquad \,\,\,\,\, $ 
We deduce from the third equation of relation (3.4.5)~: 

\noindent  $  \displaystyle 
{{(c^*)^2}\over{\gamma \!-\!1}} \, ( \alpha_{1} +  \alpha_{3} ) \,\,= \,\, 
\rho_r \, E_r \,-\, \rho_l \, E_l \,-\, {1\over2} \, (u^*)^2 \,  (\rho_r - \rho_l)
\,-\, u^* \, c^* \, (\alpha_{3} - \alpha_{1}) \,$

\noindent  $  \displaystyle \qquad  = \,\,
{{1}\over{\gamma \!-\!1}} \,(p_r - p_l ) \,+\,  {1\over2} \, (\rho_r \, u_r^2 -
\rho_l \, u_l^2 ) \,+\, {1\over2} \, (u^*)^2 \,  (\rho_r - \rho_l) \, - \, u^* \, 
(\rho_r \, u_r -  \rho_l \, u_l ) \,$

\noindent  $  \displaystyle \qquad  = \,\,
{{1}\over{\gamma \!-\!1}} \,(p_r - p_l ) \qquad  \qquad  \qquad $ due to (3.3.8).
Then we have~:  

\noindent   (3.4.15) $\qquad  \displaystyle 
(c^*)^2 \, (\alpha_{1} + \alpha_{3}) \,\,= \,\, p_r - p_l  \,.\,$

\smallskip \noindent   $\bullet \qquad \,\,\,\,\, $ 
The solution of the 2 by 2 linear system with unknowns $\, \alpha_{1} \,$ and $\,
\alpha_{3} \,$ defined by the relations (3.4.14) and (3.4.15) directly gives the
relations (3.4.10) and (3.4.12). The expression (3.4.11) of variable $\,\,
\alpha_{2} \,\,$ is a direct consequence of the relations   (3.4.10), (3.4.12) and
(3.4.13) and  Proposition 3.5 is proven.  $ \hfill \square \kern0.1mm $

\smallskip
\smallskip   \noindent  {\bf Proposition 3.6.  $\quad$  An algorithm for the Roe
flux.}

\noindent 
Let $\,\, W_l\,\,$ and $\,\, W_r \,\,$ be two compressible fluid states. The
computation of the Roe flux $\,\, \Phi(W_l,\,W_r) \,\,$ of relations  (3.4.6)-(3.4.8) 
between these two states is summarized by the following points~: 

% 
%  \noindent 

\noindent   $\bullet  \quad $
Evaluation of density $\,\, \rho^* ,\,\,$ velocity $\,\, u^* \,\,$ and total enthalpy
$\,\, H^* \,\,$ of the intermediate state $\, W^* \,$ with the relations
(3.3.4) to (3.3.6), 

\noindent   $\bullet  \quad $
Determination of the sound celerity $\,\, c^* \,\,$ of the intermediate state $\, W^*
\,$ from the previous data with the relation (3.4.2),  

\noindent   $\bullet  \quad $
Eigenvectors $\, r_{j}^* \,$ of the Roe matrix with the relations  (3.4.4), 

\noindent   $\bullet  \quad $
Computation of the characteristic variables $\,\, \alpha_{j} \,\,$ in (3.4.5) 
for the difference $\, W_r - W_l \,$ with the relations (3.4.10) to (3.4.12),  

\noindent   $\bullet  \quad $
Final computation of the Roe flux $\,\, \Phi(W_l,\,W_r) \,\,$ with the minimum of
work~: 

\noindent   (3.4.16) $\qquad  \displaystyle 
\Phi(W_l,\,W_r) \,\,= \,\, F(W_l)  \qquad \qquad \qquad \qquad \qquad  \quad \,
{\rm if} \,\, u^* - c^* \geq 0 \,$ 

\noindent   (3.4.17) $\qquad  \displaystyle 
\Phi(W_l,\,W_r) \,\,= \,\, F(W_l) \,+\, (u^* - c^*) \,\, \alpha_{1} \,\, r_{1}^* 
\qquad \quad \, {\rm if} \,\, u^* - c^* \leq 0 \,< u^* \,$

\noindent   (3.4.18) $\qquad  \displaystyle 
\Phi(W_l,\,W_r) \,\,= \,\, F(W_r) \,-\, (u^* + c^*) \,\, \alpha_{3} \,\, r_{3}^* 
\qquad  \quad {\rm if} \,\, u^* \leq 0 \, < \, u^* +  c^* \,$

\noindent   (3.4.19) $\qquad  \displaystyle 
\Phi(W_l,\,W_r) \,\,= \,\, F(W_r) \qquad \qquad \qquad \qquad \qquad  \quad 
{\rm if} \,\, u^* + c^* \leq 0 \,.\,$

\smallskip \noindent   $\bullet \qquad \,\,\,\,\, $ 
The proof of the relations (3.4.16)-(3.4.19) is obtained by starting from the
expression of the Roe flux given in (3.4.6). We know that $\,\,
\lambda_{1} = u^* - c^* ,\,\,$  $\,\, \lambda_{2} = u^* ,\,\,$ $\,\, \lambda_{3}
= u^* + c^* .\,\,$ If $\,\, u^* - c^* \geq 0 ,\,\,$ then $ \, u^* \geq 0 \,$ and $\,
u^* + c^* \geq 0 ,\,\,$ so the relation (3.4.6) reduces to (3.4.16) because, due to
(2.7.10),  $\, \mu^- = 0 \,$ if $\, \mu \,$ is a positive real number. If $\,\,   u^* -
c^* \leq 0 \,< u^* \,< u^* + c^* ,  \,\,$ the term containing $\, (\lambda_{1}^*)^-
\,$ is the only one that contributes in relation (3.4.6) and the relation (3.4.17) is a direct
consequence of this remark. If $\,\,   u^* - c^* \, < \, u^* \leq 0 \, < u^* + c^*
,\,$ the term that contains $\, (\lambda_{3}^*)^+ \,$ is the only nonzero element 
among the three  inside the relation (3.4.7) and we deduce the relation (3.4.18)
from this property. When $\,\,   u^* - c^* \, < \, u^* \, < \,  u^* + c^* \leq 0 ,\,\,$
the term $\,\,  F(W_r) \,\,$ is the only  to subsist inside the relation (3.4.7)
and the relation (3.4.19) is established. We remark also that the algebraic
expression (3.4.11) for $\,\, \alpha_{2} \,\,$ is not necessary for the implementation
of the algorithm.   $ \hfill \square \kern0.1mm $

\bigskip \noindent {\smcaps 3.5 } $ \, $ { \bf   Entropy correction. }

\noindent   $\bullet \qquad \,\,\,\,\, $ 
The Roe flux replaces the nonlinear waves of the gas dynamics, {\it i.e.} the rarefactions
and the shock waves by {\bf  linear}  waves that are the {\bf  contact
discontinuities}. If  sufficiently weak shock waves occur for a given discontinuity
between two states $\, W_{\rm left} \,$ and $\, W_{\rm right} ,\,$ the Roe flux
presented above is a good approximation, but if a rarefaction  containing a {\bf
sonic} point is present among the nonlinear waves that solves the discontinuity
problem between $\, W_{\rm left} \,$ and $\, W_{\rm right} ,\,$  it has been early
remarked that for this very  particular situation, the Roe flux does not satisfy the
entropy condition (see {\it e.g.} Godlewski and Raviart [GR96]). 

\smallskip \noindent   $\bullet \qquad \,\,\,\,\, $ 
A popular response has been proposed by Harten [Ha83] with a tuning parameter that
plays in fact the role of an artificial viscosity and P.~Roe himself [Roe85] has
proposed a nonparameterized entropy correction for his flux.  With G. Mehlman,
we have treated the same subject by the introduction of hyperbolic nonlinear models
with nonconvex flux functions and have proved a discrete entropy inequality if
sufficiently weak nonlinear waves are present in the problem [DM96]. We detail here
the modification of the algorithm that we have proposed and tested numerically for
various gas dynamics problems.

\smallskip \noindent   $\bullet \qquad \,\,\,\,\, $ 
We introduce as above two states $\,W_l \equiv W^0  \,$ and $\, W_r  \equiv W^3  \,$
and the Roe matrix $\, A(W_l,\,W_r)\,$ described  in the preceding sub-sections. We
have in particular the relation 

\noindent  (3.5.1) $\quad  
 W_r - W_l \,\,   \equiv \,\, W^3 - W^0 \,\, = \,\,   \sum_{j=1}^{j=3} \,
\alpha_{j} \,  r_{j}^* \, \,$

\noindent 
and we do not make the confusion between the eigenvalues $\, \lambda_{j}(W^0) \,$
of the left state,  $\, \lambda_{j}(W^3) \,$ of the right state, and   $\,
\lambda_{j}^*  \,$ of the Roe matrix. We introduce the following two intermediate
states $\,W^1 \,$ and $\, W^2 \,$ according to 

\noindent  (3.5.2) $\qquad  
W^1 \,\,= \,\, W^0 + \alpha_{1}  \,  r_{1}^* \,,\qquad W^2 \,\,= \,\,  W^0  +
\alpha_{1}  \,  r_{1}^* \,+\, \alpha_{2}  \,  r_{2}^* \,\,= \,\, W^3 -
\alpha_{3}  \,  r_{3}^* \, \, $ 

\noindent 
and illustrated on Figure 3.2. Note that $\,\, \lambda_{j}(W^k) \,\,$ is well
defined for $ \, j = 1 ,\, 2 ,\, 3 \,$ and $\,\, k = 0,\, 1,\, 2 ,\, 3 \, : \,$ it is
the $\, j^0 \,$ eigenvalue of the $\, k^0 \,$ intermediate state $\, W^k .\,$  We
define now the set $\, S \,$ of {\bf sonic indices} by the condition that the sign of
the $\, j^0 \,$ eigenvalue is increasing from negative to positive values accross some
$j$-wave~:

\noindent  (3.5.3) $\qquad  
S \,\, =\,\, \bigl\{ \, j \in \{1,\,2,,3 \} ,\quad  \lambda_{j}(W^{j-1}) \,<\, 0
\,<\, \lambda_{j}(W^{j}) \, \bigr\} \,. \,$

\noindent
The modification of the Roe flux is active only for the sonic indices and we introduce
a polynomial $\, p_{j} \,$ {\bf of degree 3} by the classical Hermite
interpolation conditions 

\setbox21=\hbox {$\displaystyle 
p_{j}(0) \,=\, 0 \,,\qquad   \qquad   \,\,
p'_{j}(0) \,=\, \lambda_{j}(W^{j-1}) \,,\,\,   $}
\setbox22=\hbox {$\displaystyle  
p_{j}(\alpha_{j}) \,=\, \lambda_{j}^*  \, \alpha_{j}  \,,\qquad   
p'_{j}(\alpha_{j}) \,=\, \lambda_{j}(W^{j}) \,,\,\,\, $  }
\setbox30= \vbox {\halign{#&# \cr \box21 \cr \box22    \cr   }}
\setbox31= \hbox{ $\vcenter {\box30} $}
\setbox44=\hbox{\noindent   (3.5.4)  $\displaystyle  \qquad \left\{ \box31
\right. \,\qquad   j \in S \,$ }  
\noindent $ \box44 $

\noindent
that defines explicitely the polynomial $\,\,  p_{j}({\scriptstyle \bullet})  \,\,$ by
the algebraic relation

\setbox21=\hbox {$\displaystyle 
p_{j}(\xi) \,\,= \,\, {{\lambda_{j}(W^{j}) + \lambda_{j}(W^{j-1}) - 2\,
\lambda_{j}^*}\over{(\alpha_{j})^2 }} \,\, \xi^3 \quad + \,   $}
\setbox22=\hbox {$\displaystyle  \qquad  \qquad + \quad 
 {{3 \, \lambda_{j}^* - 2 \,  \lambda_{j}(W^{j-1}) - 
\lambda_{j}(W^{j})}\over{\alpha_{j}}} \,\, \xi^2 \,\, + \,\, \lambda_{j}(W^{j-1})
\,\, \xi \,. \,$ }
\setbox30= \vbox {\halign{#&# \cr \box21 \cr \box22    \cr   }}
\setbox31= \hbox{ $\vcenter {\box30} $}
\setbox44=\hbox{\noindent   (3.5.5)  $\displaystyle  \qquad \left\{ \box31
\right. \,$ }  
\noindent $ \box44 $
  
% version mac 
% \smallskip \vskip 3.cm  \smallskip 
% \qquad  \qquad   \qquad  \qquad  \qquad  \special{illustration fig.3.2.epsf scaled 500}  
% fin de la version version mac 
% version linux
\smallskip \centerline {  \epsfysize=4,0cm  \epsfbox  {fig.3.2.epsf} }  \smallskip
% fin de la version linux     

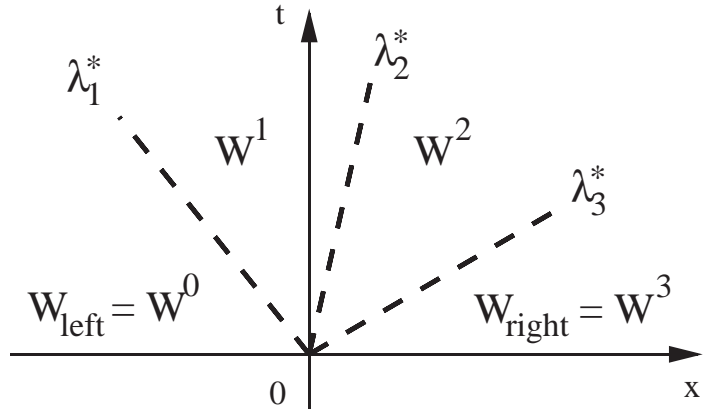
\centerline {\rm  {\bf Figure 3.2.}	\quad   Intermediate states for the entropy
 correction  }

\centerline {\rm  of the Roe upwind scheme. }

\bigskip 
\noindent   $\bullet \qquad \,\,\,\,\, $ 
With the hypothesis that $\, j \in S ,\,$ it is not difficult  to see [DM96]  that the
polynomial $\, p_{j}({\scriptstyle \bullet}) \,$ has a unique minimum  inside  the
interval $\,\, (0,\, \alpha_{j}) .\,$ The argument $\,\, \xi_{\,j}^* \,$ 
of this point of minimum is given according to~:

\setbox21=\hbox {$\displaystyle  
\bigl( \, 3 \, \lambda_{j}^* - 2\, \lambda_{j}(W^{j-1}) - \lambda_{j}(W^{j}) 
\, \bigr) \,\,+  \,  $}
\setbox22=\hbox {$\displaystyle  \, 
  \,+ \,  \sqrt { \, \bigl( \,  3 \, \lambda_{j}^* - \lambda_{j}(W^{j})
- \lambda_{j}(W^{j-1}) \,  \bigr)^2 - \lambda_{j}(W^{j-1}) \, \lambda_{j}(W^{j}) 
\,} \,  $}
\setbox25= \vbox {\halign{#&#  \cr \box21 \cr \box22 \cr   }}
\setbox26= \hbox{ $\vcenter {\box25} $}
\setbox27=\hbox{\noindent   $\displaystyle \, \left( \box26 \right)  \,$ }
\noindent  (3.5.6) $\,  \displaystyle  
\xi_{\,j}^* \,\,=\,\, {{-\lambda_{j}(W^{j-1}) \, \, \alpha_{j}}\over{\box27}}
\,. \,$

\noindent
Since $\,\, p_{j}(\xi_{\,j}^*) \,\,$ is the unique minimum of the polynomial $\,\,
p_{j}({\scriptstyle \bullet}) \,\,$ on the interval $\, (0,\, \alpha_{j}) ,\,$ we
have $\,\,\,  {{ p_{j}(\xi_{\,j}^*) }\over{\alpha_{j}}} \leq 0 \,\,$ and $\,\, 
{{ p_{j}(\xi_{\,j}^*) }\over{\alpha_{j}}} \leq  \lambda_{j}^*  \,.\,$ Then the
modified flux $\,\, \Phi^{\rm modif}(W_l, \,$ $ \,W_r) \,\,$ is defined from the Roe
flux $\,\,  \Phi(W_l,\,W_r) \,\,$ by the relation 

\noindent  (3.5.7) $\quad  \displaystyle
\Phi^{\rm modif}(W_l,\,W_r) =  \Phi(W_l,\,W_r) \,+\, \sum_{j \in S} \,
{\rm max} \biggl( \, {{ p_{j} \bigl(  \xi_{\,j}^* \bigr)   }\over{\alpha_{j}}} ,\,
{{ p_{j} \bigl(  \xi_{\,j}^* \bigr)  }\over{\alpha_{j}}} -  \lambda_{j}^* \,
\biggr) \, \alpha_{j} \,\, r_{j}^* \,$

\noindent
that makes the added numerical viscosity explicit.

\bigskip  \noindent {\smcaps 3.6 } $ \,  $ { \bf Nonlinear flux boundary conditions. }

\noindent   $\bullet \qquad \,\,\,\,\, $ 
At the two extremities $\, x \!=\! 0 \,$ and  $\, x \!=\! L  \,$ of the pipe, 
we have to express on one hand  the datum of a given nonstationary pressure $\,\, \Pi(t)
\,\,$ at $\, x \!=\! 0 \,$ and on the other hand a free output of the waves at  $\, x
\!=\! L . \,$

\smallskip \noindent   $\bullet \qquad \,\,\,\,\, $ 
For the numerical boundary condition for pressure, we follow a general approach
founded on the so-called  {\bf  partial Riemann problem} [Du01] that generalizes to 
nonlinear hyperbolic systems  the reflection operator of relation (2.5.22). For a given
discrete time $\,\, t^n = n \, \Delta t ,\,\,$ and a given state $\,\, W_{1/2}^n \equiv
W_r  \,\,$ in the first cell of the unidimensional mesh, we construct a boundary state
$\,\, W_{0}^n \equiv W_l \,\,$ that satisfies the boundary constraint 

\noindent  (3.6.1) $\qquad  \displaystyle
p(W_0^n) \,\, = \,\, \Pi^{n+1/2} \,,\quad n \geq 0 \,,\,$ 

\noindent 
and moreover, we impose that the state $\,\, W_{1/2}^n \,\,$ present in the first cell 
is issued from the boundary state $\,\, W_{0}^n \,\,$ with an ingoing  3-wave, {\it i.e.}
we impose the relation 

\noindent  (3.6.2) $\qquad  \displaystyle
 W_{1/2}^n -  W_{0}^n \,\, = \,\, \alpha_{3} \,\, r_{3}^* \,\,$ 

\noindent 
as illustrated on Figure 3.3. This problem has a unique solution, as claims the 

\smallskip   \noindent  {\bf Proposition 3.7.  ~   Pressure flux boundary
condition with Roe matrix.  }

 \noindent 
We consider a left  boundary condition associated with a pressure $\, \Pi \,$ and a
right datum defined by  state $\, W_r .\, $ Then there exists a {\bf unique}
left state $\, W_l \,$  that satisfies the pressure condition 

\noindent  (3.6.3) $\qquad  \displaystyle
p(W_l) \,\, = \,\, \Pi \,$ 

\noindent
and such that when we construct the Roe intermediate state $\, W^* \,$ according to
the relations (3.3.4) to (3.3.6), the difference $\,\, W_r - W_l \,\,$ has only one
component over the third eigenvector of the Roe matrix $\,\, {\rm d}F(W^*) ,\,\,$ {\it i.e.}
the relation (3.4.5) can been written under the form 

\noindent  (3.6.4) $\qquad  \displaystyle
W_r - W_l  \,\, = \,\, \alpha_{3} \,\, r_{3}^* \,. \,$ 

\noindent 
The density $\,\, \rho_l \,\, $ and the velocity  $\,\, u_l \,\, $ of the left state
$\, W_l \,$ are given according to the relations 

\noindent  (3.6.5) $\qquad  \displaystyle
\rho_l \,\, = \,\, {{ (\gamma \!+\!1) \, \Pi \,+ \,  (\gamma \!-\!1) \, p_r} \over{
(\gamma \!-\!1) \,\Pi \,+\, (\gamma \!+\!1) \, p_r}} \, \rho_r \,$ 

\noindent  (3.6.6) $\qquad  \displaystyle
u_l \,\,= \,\, u_r \,+\,  (\Pi - p_r ) \,\,  \sqrt{ {{2}\over{ \rho_r \, \bigl( (\gamma
\!-\!1) \, \Pi \,+ \,  (\gamma \!+\!1) \, p_r \bigr) }}} \,. \,$

% version mac 
% \smallskip \vskip 2.7cm  \smallskip 
% \qquad  \qquad   \qquad  \qquad  \special{illustration fig.3.3.epsf scaled 500}  
% fin de la version version mac 
% version linux
\smallskip \centerline {  \epsfysize=3,7cm  \epsfbox  {fig.3.3.epsf} }  \smallskip
% fin de la version linux     

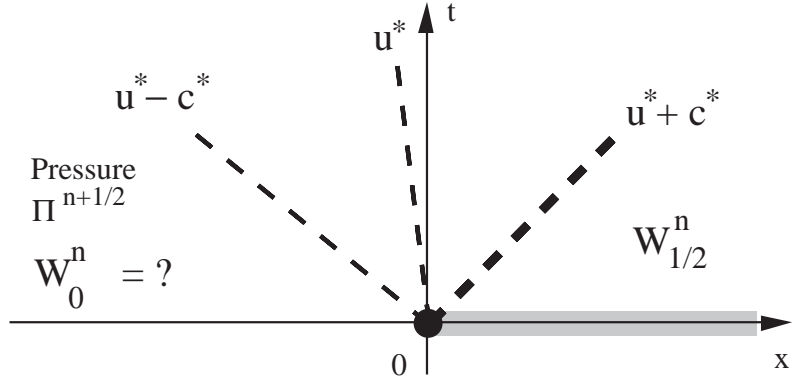
\centerline {\rm  {\bf Figure 3.3.}	\quad   Nonlinear boundary condition for 
given pressure at  $x=0$  }

\centerline {\rm  with  the Roe upwind scheme. }

\bigskip 
\noindent   $\bullet \qquad \,\,\,\,\, $ 
The {\bf proof of Proposition 3.7} is a consequence precisely of the preceding
subsections  about the Roe flux. We first remark that the relations (3.4.5) and
(3.6.4) are absolutly {\bf identical}. Then we deduce that necessarily $\,\,
\alpha_{1} =
\alpha_{2} = 0 \,\,$ and according to the relations (3.4.11) and (3.4.12), we get 

\noindent  (3.6.7) $\qquad  \displaystyle
\Pi - \rho^* \, c^* \, u_l \,\,= \,\, p_r -  \rho^* \, c^* \, u_r \,\,$  

\noindent  (3.6.8) $\qquad  \displaystyle
\Pi - (c^*)^2 \, \rho_l \,\, = \,\, p_r -  (c^*)^2 \, \rho_r  \, . \,$ 

\noindent 
We deduce simply $\quad \displaystyle u_r - u_l \,=\, {{p_r - \Pi}\over{\rho^* \,
c^*}} \,=\, {{c^* \,(\rho_r - \rho_l)}\over{\rho^*}} \quad $ and due to the  relation
(3.4.9), we get~: 

\noindent  $  \displaystyle
(c^*)^2 =	{{1}\over{ (\sqrt{\rho_r} + \sqrt{\rho_l})^2}} \, \Bigl[ \,
{{\gamma \!-\!1}\over{2}} \,  \rho^* \, \Bigl(  {{c^* \,(\rho_r
- \rho_l)}\over{\rho^*}} \, \Bigr)^2 \,+ \, \gamma \, \Bigl( \, {{\rho_l +
\rho^*}\over{\rho_l }} \, \Pi  +  {{\rho^* + \rho_r}\over{\rho_r }} \, p_r \,
\Bigr) \, \Bigr] \,.  $

\noindent
Then after multiplication by $\, (\rho_r - \rho_l) ,\,$ we obtain with the help of
(3.6.8)~:  

\noindent   $  \displaystyle  0 \,\,=  \,\,
(p_r - \Pi) \,-\, {{\gamma \!-\!1}\over{2}} \, (p_r - \Pi) \, {{(\sqrt{\rho_r} -
\sqrt{\rho_l})^2 }\over{\rho^*}}   \, $ 

\noindent   $  \displaystyle 
\,-\,  \gamma \,  {{\sqrt{\rho_r} -
\sqrt{\rho_l}}\over{\sqrt{\rho_r} + \sqrt{\rho_l}}} \,  \Bigl( \, {{\rho_l +
\rho^*}\over{\rho_l }} \, \Pi  \,\,\, + \,\, 
{{\rho^* + \rho_r}\over{\rho_r }} \, p_r \, \Bigr) \qquad = \qquad
\gamma \, (p_r - \Pi) \, $ 

\noindent   $  \displaystyle 
\,-\,  {{\gamma \!-\!1}\over{2}} \, (p_r - \Pi) \, \biggl( \, 
{\sqrt{\rho_r}\over{\sqrt{\rho_l}}} + {\sqrt{\rho_l}\over{\sqrt{\rho_r}}} \, \biggr) \,
\,-\,  \gamma \,\biggl( \, {\sqrt{\rho_r}\over{\sqrt{\rho_l}}} - 1\,  \biggr) \, \Pi
\,-\,  \gamma \,\biggl( \, 1 - {\sqrt{\rho_l}\over{\sqrt{\rho_r}}} \,  \biggr) \, p_r
\,. \,$

\noindent 
We multiply the previous equality by $\,\, \sqrt{ {{\rho_l}\over{\rho_r}} } \,\,$ and
we get~: 

\noindent   $  \displaystyle 
-{{\gamma \!-\!1}\over{2}} \, (p_r - \Pi) \, (\, 1 +
{{\rho_l}\over{\rho_r}}  \, ) \,-\,  \gamma \,  \Pi \,+\,   
\gamma \, {{\rho_l}\over{\rho_r}} \, p_r \,\,= \,\, 0 \,, \,\, $ 

\noindent   {\it id est}   $  \displaystyle \quad 
\Bigl[ - {{\gamma \!-\!1}\over{2}} \, (p_r - \Pi) \, + \, \gamma \, p_r  \Bigr] \,
{{\rho_l}\over{\rho_r}} \,\,= \,\, \gamma \, \Pi \,+\, {{\gamma \!-\!1}\over{2}} \,
(p_r - \Pi)  \,$

\noindent 
and the relation (3.6.5) is established. 

\smallskip \noindent   $\bullet \qquad \,\,\,\,\, $ 
We deduce from the previous relation~:

\noindent  $  \displaystyle 
\rho^* \,\,= \,\, \sqrt{ \, {{  (\gamma \!+\!1) \, \Pi \,+ \,  (\gamma \!-\!1) \,
p_r} \over{ (\gamma \!-\!1) \,\Pi \,+\, (\gamma \!+\!1) \, p_r}} \,} \, \, \rho_r
\,,\,$   $  \displaystyle \quad \rho_r - \rho_l \,\,= \,\, {{2 \, (p_r - \Pi) \, 
}\over{  (\gamma \!-\!1) \,\Pi \,+\, (\gamma \!+\!1) \, p_r }} \, \rho_r\,$

\noindent 
and due to the relation (3.6.8)~: \qquad 

\noindent  $  \displaystyle 
(c^*)^2 \,\,= \,\, {{  (\gamma \!-\!1) \,\Pi \,+\, (\gamma \!+\!1) \, p_r
}\over{2 \, \rho_r }} \,,\,$   $  \displaystyle \quad \rho^* \, c^* \,\,= \,\, \sqrt{\,
{{ (\gamma \!-\!1) \,\Pi \,+\, (\gamma \!+\!1) \, p_r}\over{2}} \,} \, \, \sqrt{\rho_r}
\,$

\noindent
and the relation (3.6.6) is an easy consequence of the last equality joined with
(3.6.7). The proposition 3.7 is established.   $ \hfill \square \kern0.1mm $

\smallskip 
\smallskip  
\smallskip
\noindent   $\bullet \qquad \,\,\,\,\, $ 
The determination of a  {\bf nonlinear nonreflecting} boundary condition at $\, x \!=\!
L \,$ is still an open mathematical problem. We recommand for deriving a flux boundary
condition for such a situation  to impose that {\bf no wave} are present at
 the interaction for the last interface  $\, j \!=\! J . \,$ We just write  

\noindent  (3.6.9) $\qquad  \displaystyle
f_J^{n+1/2} \,\,= \,\, F(W_{J-1/2}^n) \,\,,\qquad n \geq 0 \, \,$ 

\noindent 
which is equivalent of introducing a right boundary state $\,\, W_{J}^n \,\,$ 
according  to the simple relation $\,\,  W_{J}^n = W_{J-1/2}^n \,\,$ and then making
these two states interacting with the Roe flux~: $\,\,  f_J^{n+1/2} \,=\,
\Phi(W_{J-1/2}^n ,\,  W_{J}^n ) .\,$ This last definition is equivalent to the one
proposed in (3.6.9) due to the property  (3.3.2)  of the Roe matrix. 

\bigskip \bigskip 
\noindent  {\smcaps 4) $ \,\,\,\,\,\,\, $  Second order and  two space dimensions. } 
\smallskip \noindent {\smcaps 4.1 } $ \,  $ { \bf Towards second order accuracy. }

\noindent   $\bullet \qquad \,\,\, $ 
The finite volume method described in the previous sections is a natural method for
the discretization of systems of $\, m \,$ conservation laws. It conducts to an
explicit scheme in time~: the evaluation of the field $\,\, W^{n+1} \,\, $ at time
step $\,\, (n \!+\! 1 ) \, \Delta t \,\,$   needs  only  the knowledge of the field
$\,\, W^n_{j\!+\!1 \! / \! 2} \,\, $ for $\,\,  j = 0, \cdots, J \!-\! 1 \,\ $ at the
preceding time step $\,\, n \, \Delta t .\,\,$ This evaluation needs  a certain number
of auxiliary computations without the resolution of any linear system involving  the
new field. The method is parameterized by the choice of a numerical flux and a great
flexibility can be adopted at this level. We have proposed two fluxes for nonlinear
problems related to nonlinear acoustics and gas dynamics,  the Roe flux $\,\,
\Phi({\scriptstyle \bullet},\, {\scriptstyle \bullet})\,\,$ of relations (3.4.6)-(3.4.8) 
that conduct to a discrete scheme according to the relation 

\noindent  (4.1.1) $\qquad \displaystyle 
f_j^{n+1/2} \,\,= \,\, \Phi( W_{j-1/2}^n ,\,  W_{j+1/2}^n ) \, , \,$ 

\noindent
and the modified Roe flux $\,\, \Phi^{\rm modif}({\scriptstyle \bullet},\,
{\scriptstyle \bullet})\,\,$ of relation (3.6.7) that enforces the entropy condition. 
This explicit version of the finite volume method is submitted to  a {\bf stability}
condition that can be written as  a first approximation for linear cases as~: 

\noindent  (4.1.2) $\qquad \displaystyle 
c_0 \, {{\Delta t}\over{\Delta x}} \,\, \leq \,\, 1 \,. \,$

%%%%%%%%%%%%%%%%%%%%%%%%%%%%%%%%%       modif janvier 2011 
%%%  \titredroite={\pecaps      Second order and  two space dimensions } 

\smallskip \noindent   $\bullet \qquad \,\,\, $ 
Nevertheless,  the above finite volume method  is only {\bf first order accurate}. If
we insert an exact solution $\,\, W(x,\,t) \,\,$ of the conservation law 

\noindent  (4.1.3) $\qquad \displaystyle 
{{\partial W}\over{\partial t}} \,\, +\,\, {{\partial}\over{\partial x}}F(W) \,\,= \,\,
0 \,$ 

\noindent 
inside the formal expression of the flux (4.1.1), it is easy to see that the finite
difference $\,\, {{1}\over{\Delta x}} (f_{j+1}^{n+1/2} - f_j^{n+1/2}) \,\, $ is
first order accurate~: 

\noindent  (4.1.4) $\qquad \displaystyle  
{{1}\over{\Delta x}} \bigl( f_{j+1}^{n+1/2} - f_j^{n+1/2} \bigr) \,\, = \,\, \biggl( 
{{\partial F(W)}\over{\partial x}}  \biggr)_{j+1/2}^{n+1/2} \,\,\,+ \,\, {\rm O} \bigl( 
\Delta t + \Delta x \bigr) \,. \,$

\noindent 
In a similar way, the use of an explicit scheme in time conducts to

\noindent  (4.1.5) $\qquad \displaystyle  
{{1}\over{\Delta t}} \, ( W_{j+1/2}^{n+1} - W_{j+1/2}^n ) \,\,+\,\, 
{{1}\over{\Delta x}} \bigl( f_{j+1}^{n+1/2} - f_j^{n+1/2} \bigr) \,\, = \,\,0 \,$ 

\noindent
and maintains this first order accuracy for the finite volume scheme. 

\smallskip \noindent   $\bullet \qquad \,\,\, $ 
We  develop in this section the fact that it is possible to improve the method,
{\it i.e.} to define a method  with a relation of the type (4.1.5), and  that 
conduct to a
troncation  error {\bf of second order}~: 

\setbox11=\hbox{ $ \displaystyle 
{{1}\over{\Delta t}} \, ( W_{j+1/2}^{n+1} - W_{j+1/2}^n ) \,\,+\,\, 
{{1}\over{\Delta x}} \bigl( f_{j+1}^{n+1/2} - f_j^{n+1/2} \bigr) \,\, = \,\,  $ }
\setbox12=\hbox {   $ \displaystyle  \qquad \qquad   
=  \,\, \biggl( \, {{\partial W}\over{\partial t}} \,\, +\,\, {{\partial}\over{\partial
x}}F(W) \,\biggr)_{j+1/2}^{n+1/2} \,\,\,+\,\,\, {\rm O}
\bigl(\Delta t^2 + \Delta x^2 \bigr) \,. \,\,$  }
\setbox40= \vbox {\halign{#&# \cr \box11 \cr \box12 \cr}}
\setbox41= \hbox{ $\vcenter {\box40} $}
\setbox44=\hbox{\noindent  (4.1.6) $\qquad \left\{ \box41 \right. \,  $}  
\noindent $ \box44 $

\noindent 
The price to pay is to develop flux formulae  much more complicated than the simple
relation (4.1.1). When the second order precision (4.1.6) is achieved with a stable
scheme, the precision is sufficient to develop predictive computations in acoustics and
aerodynamics, whereas that is not the case with the initial scheme (4.1.1) (4.1.5). 

\bigskip \noindent {\smcaps 4.2 } $ \,  $ { \bf The method of lines.}

\noindent   $\bullet \qquad \,\,\, $ 
The simplest way to extend the first order  finite volume scheme is first to develop a
new  vision of the method with emphasis   more on abstraction. We have presented a
method founded on the integration of the conservation law (4.1.3)  inside the
space-time domain $\,\,  V_{j+1/2}^{n+1/2} =\,\,  ]x_{j} \,,\, x_{j+1}[  \, \times \,
]t^n ,\,t^{n+1} [ \,\,$ as suggested in (2.1.6). With the method of lines, we just
integrate the conservation (4.1.3) in {\bf space} in each control volume  
$\,\,  K_{j+1/2} =  \,\,   ]x_{j} \,,\, x_{j+1}[  . \,\,$ It is straightforward 
to introduce the mean value  $\,\, W_{j+1/2}(t) \,\,$ in this finite element~: 

\noindent  (4.2.1) $\qquad \displaystyle  
 W_{j+1/2}(t) \,\, = \,\, {{1}\over{\mid  K_{j+1/2} \mid}} \,
\int_{\displaystyle x_{j}}^{\displaystyle x_{j\!+\!1}} \, W(x,\,t) \, {\rm d}x \,; \,$ 

\noindent 
then we integrate the conservation law (4.1.3) in space in the cell $\,  K_{j+1/2} \,$
and taking into account the relation $\,\,\, {{\rm d}\over{{\rm d}t}}  
 W_{j+1/2}(t) \,=\, {{1}\over{\Delta x}} 
\int_{\displaystyle x_{j}}^{\displaystyle x_{j\!+\!1}} \, {{\partial
W}\over{\partial t}} (x,\,t) \, {\rm d}x ,\,\, \, $ we get simply

\noindent  (4.2.2) $\qquad \displaystyle  
{{\rm d}\over{{\rm d}t}}    W_{j+1/2}(t) \,\,+ \,\, {{1}\over{\Delta x}} \bigl[ \,
f_{j+1}(t) \,- \, f_{j}(t) \,\bigr] \,\,= \,\, 0 \,$

\noindent
with

\noindent  (4.2.3) $\qquad \displaystyle  
f_{j}(t) \,\,= \,\, F \bigl( W (x_j,\,t) \bigr) .\,$

\smallskip \noindent   $\bullet \qquad \,\,\, $ 
As usual with the finite volume method, a numerical scheme can be obtained from the
relations (4.2.2) (4.2.3) by replacing the relation (4.2.3) by some explicit function
over the set of all discrete variables introduced for the relation (4.2.1). To fix the
ideas, we introduce a {\bf dynamic state vector} $\,\, Z(t) \,\,$ composed by all the
dynamic variables on the finite mesh~:

\noindent  (4.2.4) $\qquad \displaystyle  
Z(t) \,\,= \,\, \bigl( \, W_{1/2}(t),\, \cdots ,\,  W_{j+1/2}(t) ,\, \cdots ,\, 
W_{J-1/2}(t) \, \bigr) \quad  \in \,\,(\R^m)^J \,.\,$

\noindent 
The discretization in space is achieved if we are able to determine the numerical
flux $\,\, f_{j}(t) \,\,$ with the help of both  the  dynamic state vector $\,
Z({\scriptstyle \bullet}) \,$ and the boundary conditions, that is the  input pressure
$\,\, \Pi(t) \,\,$ in the example   considered in the last section. As in
relation (2.1.12), we introduce a {\bf  local numerical flux function} $\,\,
\Psi_j({\scriptstyle \bullet} ,\, {\scriptstyle \bullet} ) \,\,$  relative to the
vertex $\, x_j \,$~: 

\noindent  (4.2.5) $\qquad \displaystyle  
f_{j}(t) \,\,= \,\, \Psi_j \bigl( \Pi(t) ,\, Z(t) \bigr) \,. \,$ 

\noindent
We replace the relation (4.2.3) by the numerical approximation (4.2.5) inside the
equation (4.2.2) of dynamic evolution of the state variable $\,\, 
W_{j+1/2}({\scriptstyle \bullet}) .\,\,$ We obtain the following {\bf ordinary
differential equation} 

\noindent  (4.2.6) $\qquad \displaystyle  
{{\rm d}\over{{\rm d}t}}    W_{j+1/2}(t) \,\,+ \,\, {{1}\over{\Delta x}} \bigl[ \,
\Psi_{j+1} \bigl( \Pi(t) ,\, Z(t) \bigr) \,-\,  \Psi_j \bigl( \Pi(t) ,\, Z(t) \bigr) 
\,\bigr] \,\,= \,\, 0 \,. \,$

\smallskip \noindent   $\bullet \qquad \,\,\, $ 
The {\bf method of lines} is a semi-discrete version of the finite volume method.  It
is obtained by integration in space of the conservation law without integration in
time. The result is not a numerical scheme but {\it just} an ordinary differential
equation for the dynamic state vector $\, Z({\scriptstyle \bullet}) \,$ described
component by component  with the equation (4.2.6). The method is parameterized by the 
local numerical flux functions $\,\,  \Psi_j \bigl( {\scriptstyle \bullet} ,\,
{\scriptstyle \bullet}) \,\,$ and take the general form of a dynamical system
parameterized by the pressure function $\,\, t \longmapsto \Pi(t) \,$~:

\noindent  (4.2.7) $\qquad \displaystyle  
{{\rm d}\over{{\rm d}t}} Z(t) \,\,=\,\, G(Z,\,t) \,.\,$

\noindent 
The {\bf discrete dynamic function} $\,\,(\R^m) ^J \times [0,\,+\infty[ \, \, \ni
\, (Z,\,t)  \mapsto G(Z,\,t)  \, \in (\R^m) ^J  \,\,$   is a vector valued
expression  with $ J $ components~:

\noindent  (4.2.8) $\,   \displaystyle  
G(Z,\,t) \, = \, \bigl( \, G_{1/2}(Z,\,t) ,\, \cdots ,\,  G_{j+1/2}(Z,\,t)  ,\,
\cdots ,\,  G_{J-1/2}(Z,\,t) \, \bigr) \,   \in \,(\R^m) ^J  \,$

\noindent 
and it  is defined from the $\,(J \!+\! 1) \,$ local numerical fluxes $\,\,
(\Psi_j)_{j=0,\, \cdots,\,J} \,\,$  with the very simple algebra relative to  the
finite volume method~: 

\noindent  (4.2.9) $\quad \displaystyle 
 G_{j+1/2}(Z,\,t) \,=\, - {{1}\over{\Delta x}} \bigl( \, \Psi_{j+1} (\Pi(t)
 ,\, Z) \,-\,  \Psi_j  (\Pi(t),\, Z) \,\bigr) \,,\,\,  0 \leq j \leq J-1 \,. \,$

\bigskip  \noindent  {\bf Proposition 4.1.  $\quad$  Explicit Euler scheme.}

\noindent 
With the choice of the first order scheme in space, that is 

\noindent  (4.2.10) $\qquad \displaystyle 
\Psi_j \bigl( \Pi^n ,\, Z^n \bigr) \,\,= \,\, \Phi(W_{j-1/2}^n ,\, W_{j+1/2}^n) \,\, $
\qquad if  $\,\,\, j = 1,\cdots,\, J \!-\! 1 ,\,$

\noindent
and the first order {\bf explicit} forward Euler scheme for the ordinary differential
equation (4.2.7), {\it id est}  

\noindent  (4.2.11) $\qquad \displaystyle 
{{1}\over{\Delta t}} \, (Z^{n+1} - Z^n) \,\,=\,\, G(Z^n,\,t^n) \,\,,\,$

\noindent
we recover the previous first order finite volume scheme 

\setbox21=\hbox {$\displaystyle {{1}\over{\Delta t}} \bigl( W_{j+1/2}^{n+1} - 
W_{j+1/2}^{n} \bigr)  \,\,+ \,\, {{1}\over{\Delta x}} \Bigl( \Phi \bigl( W^n_{j+1/2}
,\, W^n_{j+3/2}  \bigr) \, $}
\setbox22=\hbox {$\displaystyle  \qquad  \quad  
-\,\,  \Phi \bigl(W^n_{j-1/2} ,\, W^n_{j+1/2} \bigr) \Bigr) \,\,= \,\, 0 $ \qquad
\quad  for  $\,\,\,  j = 1,\cdots,\, J \!-\! 2 . \,$}
\setbox30= \vbox {\halign{#&# \cr \box21 \cr \box22    \cr   }}
\setbox31= \hbox{ $\vcenter {\box30} $}
\setbox44=\hbox{\noindent  (4.2.12) $\displaystyle  \quad \left\{ \box31 \right.$}  
\noindent $ \box44 $

\smallskip \noindent   $\bullet \qquad \,\,\, $ 
We write the relation (4.2.10) for the particular control volume $\,   K_{j+1/2}
\,$ and we get~:  $ \displaystyle \qquad  
W_{j+1/2}^{n+1} \,\,= \,\, W_{j+1/2}^{n} \,+\, \Delta t \,\, G(z^n,\, t^n) \, $ \qquad
due to (4.2.11), then 

\noindent $ \displaystyle 
W_{j+1/2}^{n+1} \,\,= \, W_{j+1/2}^{n} \,- \, {{1}\over{\Delta x}} \bigl( \,
\Psi_{j+1} (\Pi^n,\, Z^n) \,-\,  \Psi_j  (\Pi^n,\, Z^n) \,\bigr)   \,\, \hfill $  
due to (4.2.9)

\noindent $ \displaystyle \quad \quad \,
= \,\, W_{j+1/2}^{n} \,- \,\, {{1}\over{\Delta x}} \bigl( \, \Phi(W_{j+1/2}^n ,\,
W_{j+3/2}^n) -  \Phi(W_{j-1/2}^n ,\, W_{j+1/2}^n) \, \bigr) \, $ {\it c.f.}  (4.2.10)

\noindent
and the relation (4.2.12) is established.   $ \hfill \square \kern0.1mm $

\bigskip \noindent {\smcaps 4.3 } $ \,  $ { \bf The  method of Van Leer.}

\noindent   $\bullet \qquad \,\,\, $ 
We turn now to the construction of a second order accurate version of the finite
volume method as proposed initially with the  ``Multidimensional Upwindcentered Scheme
for Conservation Laws''  of   B. Van Leer [VL79]. The fundamental idea of this scheme
is the  reconstruction of a function $\,\, \R \ni x \longmapsto W(x) \in \R \,\,$
from his mean values $\,\, W_{j\!+\!1/2}  \,\,$ in each cell $\, K_{j+1/2} \,.\, $
The reconstructed function is regular inside each control volume $\,\,  K_{j+1/2} \,\,$
and is  discontinuous at the interfaces $\,\, x_l \,\,$ between two control volumes.
The application to the finite volume method  replaces the scheme  (4.1.1) by the same
Roe flux interaction $\,\, \Phi({\scriptstyle \bullet},\, {\scriptstyle \bullet}) \,\,$
  considered for the two extrapolated data $\,\,  W_{j}^- \,\,$ and  $\,\,  W_{j}^+
\,\,$ on each side of the boundary~: 

\noindent  (4.3.1) $\qquad \displaystyle 
f_j \,\,= \,\,  \Phi( W_{j}^- ,\,  W_{j}^+ ) \,.\, $
 
% version mac 
% \smallskip \vskip 4.cm  \smallskip 
% \qquad  \qquad   \qquad   \special{illustration fig.4.1.epsf scaled 500}  
% fin de la version version mac 
% version linux
\smallskip \centerline {  \epsfysize=4,7cm  \epsfbox  {fig.4.1.epsf} }  \smallskip
% fin de la version linux     

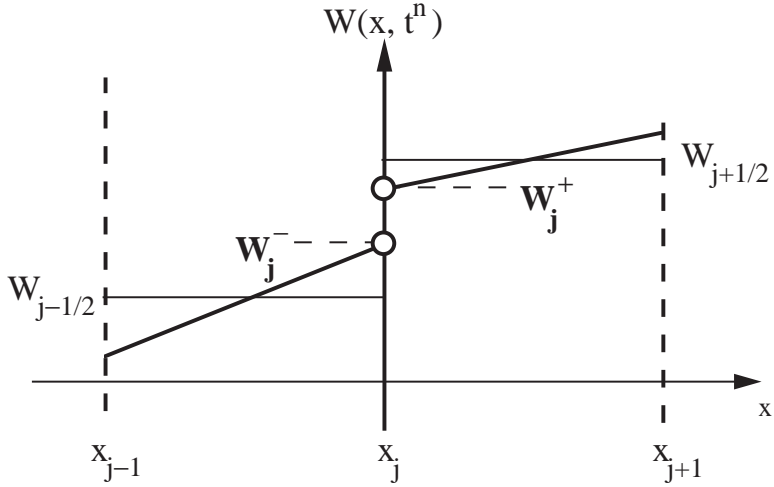
\centerline {\rm  {\bf Figure 4.1.}	\quad   First order and second order interpolation
at the interface $\,\, x_j .\,$ }
 
% version mac 
% \smallskip \vskip 4.8cm  \smallskip 
% \qquad  \qquad   \qquad   \special{illustration fig.4.2.epsf scaled 500}  
% fin de la version version mac 
% version linux
\smallskip \centerline {  \epsfysize=5,2cm  \epsfbox  {fig.4.2.epsf} }  \smallskip
% fin de la version linux     

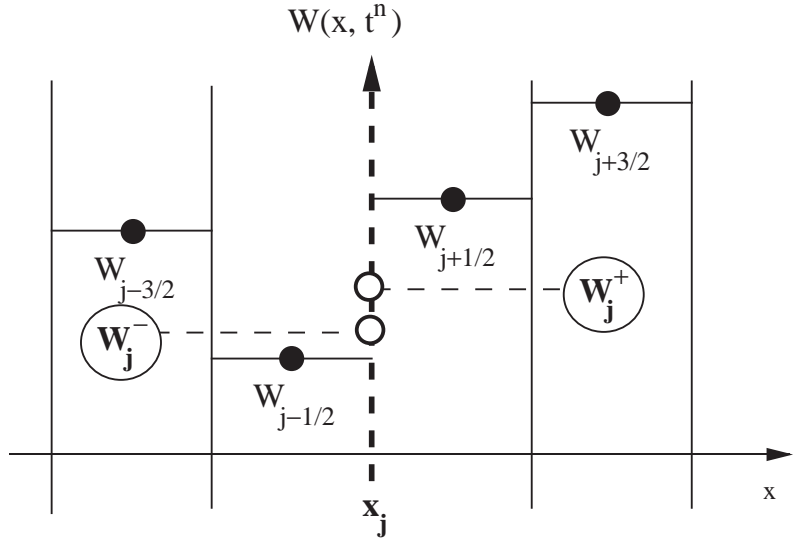
\centerline {\rm  {\bf Figure 4.2.}	\quad   Construction of the nonlinear interpolated
 }

\centerline {\rm  values $\,\,  W_{j}^- \,\,$ and  $\,\,  W_{j}^+ \,\,$  with the Van
Leer method. }
\bigskip 
\noindent   $\bullet \qquad \,\,\, $ 
The simplest case is illustrated on Figure 4.1. It imposes  simply   the
reconstructed function $\, W(x) \,$ to be  {\bf  constant} in each interval~:  

\noindent  (4.3.2) $\qquad \displaystyle 
W(x) \,\, \equiv \,\, W_{j+1/2} \,,\qquad x_j \,\, < \,\, x \,\, < \,\, x_{j+1} \,.$ 

\noindent
The two limit values on each side of the interface $\, x_j \,$ are the following
ones~: $\quad W_{j}^-  \,=\,  W_{j-1/2}  \,,$ $\quad W_{j}^+  \,=\,  W_{j+1/2} 
\,\,$ and the explicit version of this finite volume scheme is the standard first
order numerical flux (4.1.1) as seen at Proposition 4.1. In the context of the method
of lines, we obtain~: 

\noindent  (4.3.3) $\qquad \displaystyle 
f_j \,\,= \,\,  \Phi( W_{j-1/2} ,\,  W_{j+1/2} ) \, . \, $

%  \smallskip 
\vfill \eject   %%% janvier 2011 
\noindent   $\bullet \qquad \,\,\, $ 
The second order accurate Muscl  method consists first in restricting the methodology
to a {\bf scalar} field $\, W({\scriptstyle \bullet})\,\,$  and to construct an {\bf
affine} function in each interval $\,\,  K_{j+1/2}  \,\,$ instead of a constant
function as in (4.3.2). We set 

\noindent  (4.3.4) $\qquad \displaystyle 
W(x) \,\, \equiv \,\, W_{j+1/2} \,+\, p_{j+1/2} \, \bigl(x - x_{j+1/2}\bigr) \,,\qquad
x_j \,\, < \,\, x \,\, < \,\, x_{j+1} \,.$ 

\noindent
The simplest choice for a slope is the one of the {\bf centered} scheme~: 

\noindent  (4.3.5) $\qquad \displaystyle 
p_{j+1/2} \,\,= \,\, {{1}\over{\Delta x}} (W_{j+3/2} - W_{j-1/2} ) \,,\,$  

\noindent 
and due to (4.3.4) and (4.3.5), the extrapolated values $\,\,  W_{j}^- \,\,$ and
$\,\,  W_{j}^+  \,\,$ on each side of the interface located at the position $ \, x_j
\,$  are the following ones~: 

\noindent  (4.3.6) $\qquad \displaystyle 
 W_{j}^- \,\,= \,\, W_{j-1/2} \,+\, {1\over4} \,(W_{j+1/2} - W_{j-3/2} ) \,$

\noindent  (4.3.7) $\qquad \displaystyle 
 W_{j}^+ \,\,= \,\,  W_{j+1/2} \,-\, {1\over4} \,(W_{j+3/2} - W_{j-1/2} ) \, . \,$ 

\smallskip \noindent   $\bullet \qquad \,\,\, $ 
The choice of a numerical flux given according to the relations 

\noindent  (4.3.8) $ \,\,  \displaystyle 
f_j \,=\,  \Phi \bigl(  W_{j-1/2} \,+\, {1\over4} (W_{j+1/2} - W_{j-3/2} )\,,\, 
W_{j+1/2} \,-\, {1\over4} (W_{j+3/2} - W_{j-1/2} ) \bigr) \,\,$ 

\noindent 
lead  to an {\bf unstable} scheme when we consider the particular case of the
advection equation with the first order explicit scheme in time. 

\bigskip  \noindent  {\bf Proposition 4.2.  $\quad$  Linear Muscl scheme is
unstable.}

\noindent 
We apply the linear Muscl approach for the advection equation.  Then  the numerical
scheme obtained by association of  (4.3.8) and the  upwind scheme (2.1.19) 
conducts to the following  explicit first order scheme~: 

\setbox11=\hbox{ $ \displaystyle 
W_{j+1/2}^{n+1} \,-\, W_{j+1/2}^n \,+\, {{a \, \Delta t}\over{\Delta x}} \, \Bigl( \,
\bigl( W_{j+1/2}^n  \,+\, {1\over4} \,(W_{j+3/2}^n - W_{j-1/2}^n \bigr)   $ }
\setbox12=\hbox {   $ \displaystyle  \qquad  \qquad  \qquad  \qquad 
\,- \,  \bigl( W_{j-1/2}^n \,+\, {1\over4} \,(W_{j+1/2}^n - W_{j-3/2}^n \bigr) \,
\Bigr) \,\,= \,\, 0 \,.\, $   }
\setbox40= \vbox {\halign{#&# \cr \box11 \cr \box12 \cr}}
\setbox41= \hbox{ $\vcenter {\box40} $}
\setbox44=\hbox{\noindent  (4.3.9) $\qquad \left\{ \box41 \right. \,  $}  
\noindent $ \box44 $

\noindent
This scheme is unstable for each $\,\, \Delta t > 0 .\,$ 
%  \smallskip }\noindent {\vrule  height \ht12 depth 0pt width 0,05cm} \quad \box12

\smallskip \noindent   $\bullet \qquad \,\,\, $ 
Due to the expression (2.1.19) of the upwind scheme, the discrete first order in time
advection equation can be written~:   

\noindent   $  \displaystyle 
W_{j+1/2}^{n+1} \,-\, W_{j+1/2}^n \,+\, {{a \, \Delta t}\over{\Delta x}} \, \bigl(
W_{j\!+\!1}^{-,\,n} - W_{j}^{-,\,n} \bigr) \,\,=\,\, 0 \,$ 

\noindent
and the expression
(4.3.9) is a consequence of the left extrapolation (4.3.6). For the study of
stability, we introduce a profile of the type $\,\, W_{j+1/2}^n \,= \,  {\rm
e}^{( i \, k \, (j+1/2) \, \Delta x )} \,\,$  with a wave number $\, k .\,$  The
scheme (4.3.9) can be written as 

\noindent   $  \displaystyle 
W_{j+1/2}^{n+1} \,=\, g \bigl( k \, \Delta x
,\,{{a\, \Delta t}\over{\Delta x}} \bigr) \, W_{j+1/2}^{n} \,\,$ 

\noindent
with an amplification
coefficient $\,\,  g ( \xi ,\, \sigma ) \,\,$ $( \xi = k \, \Delta t,\,  \, \sigma = 
{{a\, \Delta t}\over{\Delta x}}) \,\,$  given simply by the expression 

\noindent $ \displaystyle 
g ( \xi ,\, \sigma ) \,\,= \,\,  1 \,- \, \sigma \,\bigl( \,  1 - {\rm e}^{-i \, \xi}
\, \bigr)  \,-\, {{\sigma}\over{4}}\, \bigl( \, {\rm e}^{i \, \xi} - 1 - {\rm e}^{-i \,
\xi} + {\rm e}^{-2 \, i \, \xi} \, \bigr) \,,\, $ \qquad then

\noindent $ \displaystyle 
g ( \xi ,\, \sigma ) \,\,= \,\,  1 \,- \, \sigma \,(\, 1 - {\rm cos} \, \xi \, ) \,-\,
{{\sigma}\over{4}}\, (\, -1 + {\rm cos} \, 2 \xi \,) \,- \, i \, \sigma \, \Bigl( 
 {\rm sin} \, \xi \,+\,  {1\over2}\,  {\rm sin} \, \xi - {1\over4} \, {\rm sin} \, 2
\xi \, \Bigr) \,$

\noindent $ \displaystyle \qquad \quad \,\,\, = \,\, 
1 \,- \, \sigma \,\Bigl( \,  1 - {\rm cos} \, \xi \,  + \, {1\over4}\, ( \,2 \, {\rm
cos}^2 \, \xi - 2 \, ) \, \Bigr)   \,- \, i \, \sigma \, \Bigl( {3\over2}\,  {\rm sin}
\, \xi \,-\,  {1\over2} \, {\rm sin} \,  \xi \, {\rm cos} \, \xi \,  \Bigr) \,\,,\,$

\noindent  (4.3.10) $\qquad \displaystyle 
g ( \xi ,\, \sigma ) \,\,= \,\, 1  \,-\, {{\sigma}\over{2}}\,(\, 1 - {\rm cos} \, \xi
\, )^2 \,-\, i \,  {{\sigma}\over{2}}\, {\rm sin} \, \xi \, ( \, 3 -  {\rm cos} \, \xi
\,) \,. \,$ 

\noindent 
For $\,\, \xi \,\,$ arbitrarily small, we deduce from (4.3.10)~: $\quad \abs{g ( \xi
,\, \sigma ) }^2 \,=\, 1 \,+\,  \sigma^2 \, \xi^2 \,+\, {\rm O}\bigl(\xi^4 \bigr) \, $
which establishes the instability for all $\,\, \sigma \neq 0 .\,$  $ \hfill \square
\kern0.1mm $

\smallskip \smallskip  \noindent   $\bullet \qquad \,\,\, $ 
The above remark motivates the introduction of so-called {\bf slope limiters},
intensively studied during the period 1980-90 after the pioneering work of Van Leer
[VL77]. The idea is to search an interpolation  $\,\, W_{j}^- \,\,$ of the field $\,\,
W({\scriptstyle \bullet}) \,\,$ at the left of the point $\, x_j \,$ from the
neighbouring mean values $\,\, W_{j-3/2} ,\,$ $\,\, W_{j-1/2} \,\,$ and  $\,\,
W_{j+1/2} \,\,$ and by left-right invariance of the procedure, to construct an
interpolated value   $\,\, W_{j}^+ \,\,$ from the first right  neighbours 
$\,\, W_{j-1/2} ,\,$ $\,\, W_{j+1/2} \,\,$ and  $\,\, W_{j+3/2} \,\,$ as suggested on
Figure 4.2. We replace the relations (4.3.6) and (4.3.7) by a {\bf nonlinear
interpolation} parameterized by a    slope limiter  $\,\, \R \ni r 
\longmapsto \varphi(r) \in \R \,\,$~: 

\noindent  (4.3.11) $\qquad \displaystyle 
 W_{j}^- \,\,= \,\, W_{j-1/2} \,+\, {1\over2} \,\varphi \biggl(  {{ W_{j-1/2}
-W_{j-3/2}}\over{W_{j+1/2} - W_{j-1/2}}}  \biggr) \,\, \bigl( W_{j+1/2} - W_{j-1/2}
\bigr) \,$

\noindent  (4.3.12) $\qquad \displaystyle 
 W_{j}^+ \,\,= \,\, W_{j+1/2} \,-\,  {1\over2} \,\varphi \biggl(  {{ W_{j+3/2}
-W_{j+1/2}}\over{W_{j+1/2} - W_{j-1/2}}}  \biggr) \,\, \bigl( W_{j+1/2} - W_{j-1/2}
\bigr) \,. \,$

\noindent
We remark that the relations (4.3.6) and (4.3.7) are a particular case of the general
nonlinear relations (4.3.11) and (4.3.12) with the particular choice $\,\, \varphi(r)
\,=\, {1\over2}\,(1 + r) \,.\,$ The limiter function satisfies very often  the
functional  relation 

\noindent  (4.3.13) $\qquad \displaystyle 
\varphi(r) \,\,\equiv \,\, r \, \varphi \Bigl( {{1}\over{r}} \Bigr) \,\,, \qquad r \, >
\, 0 \,. \,$  

\noindent
Among all the possible choices, we have adopted for fluid mechanics [DM92]
the so-called STS-limiter defined by the relations 

\setbox11=\hbox{ $ \displaystyle 0 \,\,, \qquad   $ }
\setbox21=\hbox{ $ \displaystyle   r \leq 0 \,  $ }
\setbox12=\hbox { $ \displaystyle  {3\over2}\,r  \,\,, \qquad $ } 
\setbox22=\hbox{ $ \displaystyle   0 \, \leq \, r \, \leq \, {1\over2} \,\,$  }
\setbox13=\hbox { $ \displaystyle  {{1 + r}\over{2}} \,,\qquad   $ }
\setbox23=\hbox{ $ \displaystyle    {1\over2}\, \leq \, r \, \leq \,2 \,$ }
\setbox14=\hbox { $ \displaystyle {3\over2} \,\,, \qquad   $ }
\setbox24=\hbox{ $ \displaystyle     r \geq 2 \,, \,  $ }
\setbox40= \vbox {\halign{#&# \cr \box11 & \box21 \cr \box12 & \box22 \cr \box13  &
\box23  \cr \box14  & \box24  \cr}}
\setbox41= \hbox{ $\vcenter {\box40} $}
\setbox44=\hbox{\noindent  (4.3.14) $ \qquad \qquad \quad 
\varphi^{STS}(r) \,\,= \,\,  \left\{ \box41
\right. \,  $}  
\noindent $ \box44 $

\noindent 
and illustrated on Figure 4.3. 
  
% version mac 
% \smallskip \vskip 4.5cm  \smallskip 
% \qquad  \qquad   \qquad  \quad  \special{illustration fig.4.3.epsf scaled 500}  
% fin de la version version mac 
% version linux
\smallskip \centerline {  \epsfysize=5cm  \epsfbox  {fig.4.3.epsf} }  \smallskip
% fin de la version linux     

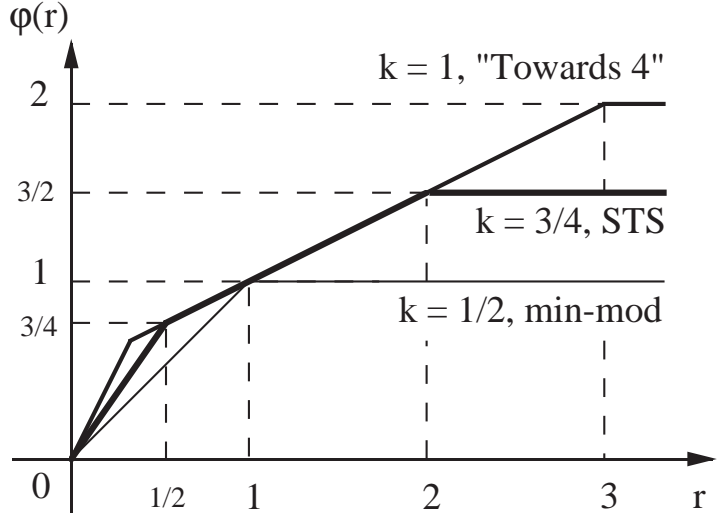
\centerline {\rm  {\bf Figure 4.3.}	\quad   Examples of limiter functions that can be 
easily extended  }

\centerline {\rm  to unstructured meshes. }
\bigskip 
\noindent {\smcaps 4.4 } $ \,  $ { \bf  Second order accurate finite volume
method   for fluid problems.}

\noindent   $\bullet \qquad \,\,\, $ 
We detail in this section a generalization  for unstructured meshes of the  Muscl
scheme proposed by   Van Leer [VL79]. At one space dimension on a uniform mesh, it is
classical to consider a scalar field $\,\, z \,\, $ among the primitive variables,
{\it id est } 

\noindent  (4.4.1) $\qquad \displaystyle
z \,\, \in \,\, \{ \,  \rho \,,\, u \,,\, v \,,\, p \, \} \qquad \, $  
(primitive variables)

\noindent 
and instead of computing the interface  flux with relation (4.1.1), to first construct
two interface states $\,\, W_S^- \,\,$ and $\,\, W_S^+ \,\,$ on each side of the
interface $\, S .\,$ Then the flux is evaluated by the decomposition of the
discontinuity~: 

\noindent  (4.4.2) $\qquad \displaystyle
f_{S} \,\,= \,\, \Phi \bigl( W_S^- \,,\, W_S^+  \bigr) 
\,,\quad \, S \in \bigl\{ \, x_{1},\cdots ,\, x_{J-1} \, \bigr\} \,. \,$ 

\noindent 
This nonlinear interpolation is done with a slope limiter  $\,\, \varphi(
{\scriptstyle \bullet}) \,\,$  that operates on each variable proposed in (4.4.1) and
we have typically when a left-right invariance is assumed [Du91]~:

\noindent  (4.4.3) $\quad \displaystyle 
z_{S}^- \,\,\,=\,\,\, z_{j-1/2}  \,\,+ \,\,{1\over2} \, \varphi \biggl( 
 {{z_{j-1/2}-z_{j-3/2}}\over {z_{j+1/2}-z_{j-1/2}}} \biggr)
\, \bigl( z_{j+1/2}-z_{j-1/2}  \bigr) \,, \quad \, S \,=\, x_{j}\,$ 

\noindent  (4.4.4) $\quad \displaystyle 
z_{S}^+ \,\,\,=\,\,\, z_{j+1/2}  \,\,- \,\,{1\over2} \, \varphi \biggl( 
 {{z_{j+3/2}-z_{j+1/2}} \over {z_{j+1/2}-z_{j-1/2}}} \biggr)
\, \bigl( z_{j+1/2}-z_{j-1/2}  \bigr) \,, \quad \, S \,=\, x_{j} \,. \,$ 

% version mac 
% \smallskip \vskip 4.0cm  \smallskip  
% \qquad  \qquad   \qquad   \qquad  \qquad  \special{illustration fig.4.4.epsf scaled 500}
% fin de la version version mac 
% version linux
\smallskip \centerline {  \epsfysize=5,0cm  \epsfbox  {fig.4.4.epsf} }  \smallskip
% fin de la version linux     

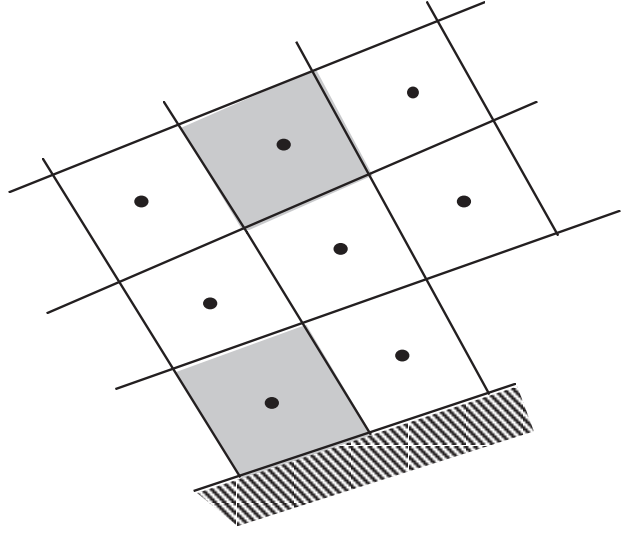
\centerline {\rm  {\bf Figure 4.4.}	\quad   Structured Cartesian mesh. }

\centerline {\rm  The control volumes are exactly the elements of mesh  $\, \cal T
.\,$ }
 
% version mac 
% \smallskip \vskip 3.0cm  \smallskip 
% \qquad  \qquad   \qquad   \qquad   \quad  \special{illustration fig.4.5.epsf scaled 500} 
% fin de la version version mac 
% version linux
\smallskip \centerline {  \epsfysize=4,0cm  \epsfbox  {fig.4.5.epsf} }  \smallskip
% fin de la version linux  

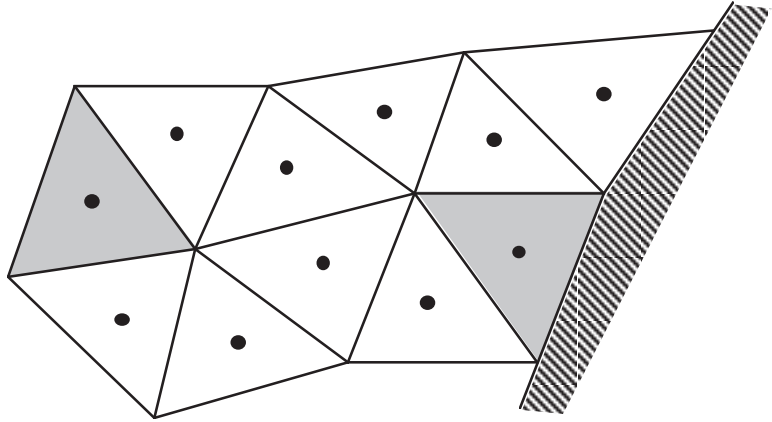
\centerline {\rm  {\bf Figure 4.5.}	\quad   Unstructured mesh composed by triangular 
elements.   }

\centerline {\rm   The control volumes are exactly the elements of the mesh.  }
\bigskip 

\noindent   $\bullet \qquad \,\,\, $ 
We focus now on the use of unstructured meshes  for  the extension to second order  
accuracy of the  finite volume method.  As in the one-dimensional case, the domain of
study is decomposed into finite elements (or control volumes)  $\,\,  K \in {\cal
E}_{\cal T} \,\,$  than can be structured in a Cartesian way (Figure 4.4) or with a
cellular complex as in Figure 4.5. In both cases,  the intersection of two finite
elements define an interface $\, \, f \! \in  \! {\cal F}_{\cal T} .\,\,$ We denote
by   $\, {\bf n}_{f} \,$ the normal at the interface~$\, f \, $ that separates a
left control volume $\, K_l(f) \,$ and a right control volume $ \, K_r(f) .\,$ The
ordinary differential equation (4.2.6) is replaced by a multidimensional version~: 

\noindent  (4.4.5) $\qquad \displaystyle 
\abs{K} \,\, {{{\rm d}W_{K}}\over{{\rm d}t}} \,+\, \sum_{f \subset \partial K}
\,\abs{f} \, \Phi\bigl(  W_{K} \,,\,  {\bf n}_{f} \,,\, W_{ K_r(f)} \bigr)  \,\,=
\,\, 0 \,,\qquad   K \in  {\cal E}_{\cal T} \,. \,$ 

\noindent 
For internal interfaces, the function  $\, \Phi\bigl(   {\scriptstyle \bullet} \,,\, 
{\bf n}_{f} \,,\,  {\scriptstyle \bullet}) \,\,$ is equal {\it e.g.} 
to the Roe flux    between states $\,  W_{ K_l(f)} \, $ and $\,  W_{ K_r(f)} \, $
in the one-dimensional direction along normal $\,  {\bf n}_{f} \,$  in order to take
into account the invariance by rotation of the equations of gas dynamics (see 
[GR96]). 

\smallskip \noindent   $\bullet \qquad \,\,\, $ 
We consider now a finite element $\, K \,$ internal to the domain. The extension to
second order accuracy of the finite volume scheme consists in replacing the arguments 
$\,  W_{ K_l(f)} \, $ and $\,  W_{ K_r(f)} \, $ in relation (4.4.5) by
nonlinear extrapolations  $\,\,  W_{K_l(f),\,f} \,\,$  and $\,\,  W_{K_r(f),\,f}
\,\,$  on each side of the boundary of  state data and  evaluated as described in what 
follows. We first introduce the set $\, {\cal N}(K) \,$ of neighbouring cells of given
finite element $\, K \in {\cal E}_{\cal T} ,\,$ as illustrated on Figure 4.6~: 

\noindent  (4.4.6) $\qquad \displaystyle 
{\cal N}(K) \,\,= \,\, \bigl\{ \, L \in {\cal E}_{\cal T} , \quad \exists \, f \in
{\cal F}_{\cal T} , \quad  f \, \subset \,\, \partial K \cap \partial L \, \bigr\}
\,.\,$ 

\noindent 
For $\, \, L \! \in \! {\cal N}(K) ,\,$ we suppose by convention that the normal  $\, 
{\bf n}_{f} \,$ to the face $\, f  \subset  \partial K \cap \partial L \,$ is external
to the  element $\, K \,$ {\it id est} $ \,\, K_r(f) = K ,\, K_l(f) = L .\,$ We introduce
also the point $\,\, y_{K,\,f}   \,\,$ on the interface $\,\, f \subset \partial K
\,\, $  that links the barycenters $\, x_{K} \,$ and   $\, x_{K_r(f)} \,$~:

\setbox11=\hbox{   $\,\,  y_{K,\,f}   \,\, \equiv \,\,  (1-\theta_{K,\,f} ) \,
x_{K} \,\,+\,\, \theta_{K,\,f}  \,  x_{K_r(f)}  \,\,,\qquad    y_{K,\,f} \in f
\,, \, $ }
\setbox12=\hbox {  $\qquad \qquad  f  \subset  \partial K \, , \qquad \, K \,$ finite 
 element internal to mesh $\,\, {\cal T} \,. \,$  }
\setbox40= \vbox {\halign{#&# \cr \box11 \cr \box12 \cr}}
\setbox41= \hbox{ $\vcenter {\box40} $}
\setbox44=\hbox{\noindent  (4.4.7) $\qquad \left\{ \box41 \right. \,  $}  

\noindent $ \box44 $

\noindent
Then, following Pollet [Po88], for $\, z \,$ equal to one scalar variable of the
 family~: 

\noindent  (4.4.8) $\qquad \,\, \displaystyle
z \,\, \in \,\, \{ \,  \rho \,,\, \rho \,u \,,\, \rho \,v \,,\, p \, \} \qquad \, $  

\noindent
we evaluate a mean value $\, \overline{z_{K,\,f}} \,\, $ on the interface $\,f \,: $ 

\noindent  (4.4.9) $\qquad \,\,  \displaystyle
\overline{z_{K,\,f}} \,\,= \,\,  (1-\theta_{K,\,f} ) \, z_{K} \,\,+\,\,
\theta_{K,\,f}  \,  z_{K_r(f)} \,$

\noindent
and the gradient $\,\, \nabla z(K) \,\,$ of field $\,\, z({\scriptstyle \bullet}) \,
\,$  in volume $\, K \,$ with a Green formula~:

\noindent  (4.4.10) $\quad \displaystyle
\nabla z(K) \,=\, {{1}\over{\mid K \mid}} \, \int_{\partial K} \, \overline{z} \, 
\, {\bf n} \, {\rm d}\gamma \,=\,  {{1}\over{\mid K \mid}} \, \sum_{f \subset
\partial K} \, \abs{f} \,  \overline{z_{K,\,f}}   \,\,  {\bf n}_{f} \,,\,\, 
   K \in  {\cal E}_{\cal T} \,. \,$ 

%%%%%%%%%%%%%%%%%%%%%%%%%%%%%%%%%%%%%%%%%%%%%%%%%%%%%%%%%%%%%%%%%%%%%%%%%%%%%%%%%%%%%%%
% \smallskip 
\vskip  -.2cm     %%%%%%%%%%%%%  janvier 2011 
\centerline {  \epsfysize=4,7cm  \epsfbox  {fig.4.6.epsf} }  %  \smallskip

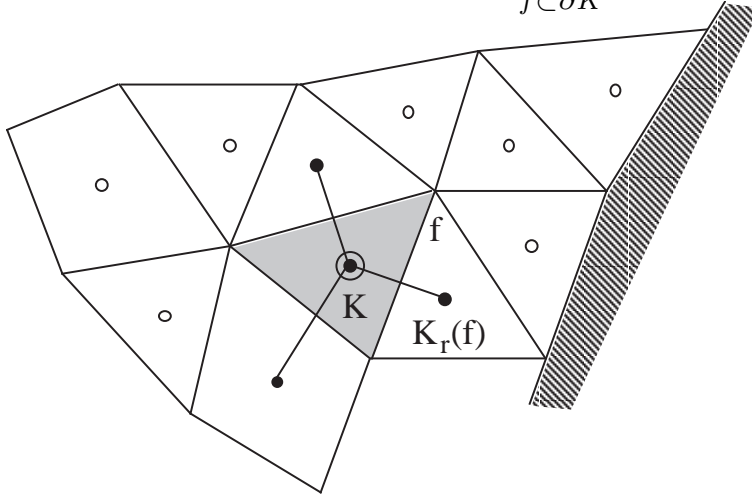
\centerline {\rm  {\bf Figure 4.6.}	\quad   Cellular complex mesh with triangles and
quadrangles. }

\centerline {\rm  Three neighbouring cells are necessary to determine the gradient }

\centerline {\rm  in triangle  $\, K \,$ and to limit eventually its variation. }
%%%%%%%%%%%%%%%%%%%%%%%%%%%%%%%%%%%%%%%%%%%%%%%%%%%%%%%%%%%%%%%%%%%%%%%%%%%%%%%%%%%%%%%

\bigskip 
\noindent   $\bullet \qquad \,\,\, $ 
An ideal extrapolation of field $\,\,  z({\scriptstyle \bullet}) \, \,$ at the
interface $\, f \,$ would be~:

\noindent  (4.4.11) $\qquad \displaystyle
z_{K,\,f} \,\,= \,\, z_{K} \,+\,  \nabla z(K) \, {\scriptstyle \bullet} \, \bigl( 
 y_{K,\,f}  - x_{K} \bigr) \,$
 
\noindent
but the corresponding scheme is unstable as explicited at Proposition 4.2.  When the
variation $\,\,  \nabla z(K) \, {\scriptstyle \bullet} \, \bigl(   y_{K,\,f}  -
x_{K} \bigr) \,\,$  is very  important, it has to be ``limited'' as first suggested by
Van Leer [VL77]. For doing this in a very general way, we introduce the minimum $\, 
m_{K}(z) \,$  and the  maximum $\,  M_{K}(z) \,$  of field  $\,\,  z({\scriptstyle
\bullet}) \, \,$ in the neighbouring cells~: 

\noindent  (4.4.12) $\qquad \displaystyle
 m_{K}(z) \,\,\,=\,\, {\rm min} \, \bigl\{ \, z_{L} \,, \quad L \in {\cal N}(K) 
 \, \bigr\} \, $

\noindent  (4.4.13) $\qquad \displaystyle
M_{K}(z) \,\,=\,\, {\rm max} \, \bigl\{ \, z_{L} \,, \quad L \in {\cal N}(K) 
 \, \bigr\} \, . \, $

\noindent
If the value $\,  z_{K} \,$ is extremum among the neighbouring ones, {\it i.e.} if 
$\,\, z_{K} \, \leq \,  m_{K}(z) \,,\,$ or $\,\, z_{K}\, \geq \,  M_{K}(z) ,\,\,$
we impose that the interpolated value $\,\, z_{K,\,f} \,\,$ is equal to the cell value
$\,\,  z_{K} \,$~: 

\noindent  (4.4.14) $\quad \displaystyle
z_{K,\,f} \,\,= \,\,  z_{K} \qquad {\rm if} \quad  z_{K} \, \leq \,  m_{K}(z)
\quad {\rm or} \quad   z_{K}\, \geq \,  M_{K}(z) \,\,,\quad f \subset \partial K \,.
\,$ 

\noindent 
When on the contrary $\,\,  z_{K} \,\,$ lies inside the interval $\,\, [ m_{K}(z)
,\,  M_{K}(z)] ,\,$ we impose that the variation $\,\, z_{K,\,f} - z_{K} \,\,$ is
limited by some coefficient $\, k \,\,\, ( 0 \leq k \leq 1) \,\,$ multiplied by the
variations $\,\,  z_{K} - m_{K}(z) \,\,$ and $\,\, M_{K}(z) - z_{K} . \,\,$ We
introduce a nonlinear extrapolation of the field $\,\,  z({\scriptstyle \bullet})\,
\,$ between center $\, x_{K} \,$ and boundary face $\,  y_{K,\,f} \, (f \subset
\partial K) $~: 

\noindent  (4.4.15) $\qquad \displaystyle
z_{K,\,f} \,\,= \,\,  z_{K} \,+\, \alpha_{K}(z) \, \nabla z(K) \,{ \scriptstyle
\bullet} \, \bigl( y_{K,\,f}  - x_{K}  \bigr) \,\,,\qquad f \subset \partial K \,$

\noindent 
with a {\bf  limiting coefficient}  $\,\,  \alpha_{K}(z) \,\,$ satisfying the
following conditions~:

\setbox10=\hbox{   $\!\!0 \,\, \leq \,\, \alpha_{K}(z)  \,\, \leq \,\, 1 \,\,, 
\quad  z({\scriptstyle \bullet})
\, \,$ scalar field defined in (4.4.8),$\quad K \in {\cal E}_{\cal T} \,$ }
\setbox11=\hbox{   $\!\!  	k \, (z_{K}- m_{K}(z)) \, \leq \,  \alpha_{K}(z)
\,\nabla z(K) \,{\scriptstyle \bullet} \,  \bigl( y_{K,\,f}  - x_{K}  \bigr) \, \leq
\,   	k \, (M_{K}(z) - z_{K}) \,$ }
\setbox12=\hbox { $ \qquad \qquad \qquad \qquad \qquad \qquad \qquad \qquad \qquad
\forall \, f \subset  \partial K \,, \quad K \in {\cal E}_{\cal T} \,. \,$ }
\setbox40= \vbox {\halign{#&# \cr \box10 \cr \box11 \cr \box12 \cr}}
\setbox41= \hbox{ $\vcenter {\box40} $}
\setbox44=\hbox{\noindent  (4.4.16) $\,\, \displaystyle  \left\{ \box41 \right. \,
 $}  

\noindent $ \box44 $

\noindent
Then $\, \alpha_{K}(z) \,\,$ is chosen as large  as possible and less than  or equal to 1 in
order to satisfy the constraints (4.4.16)~: 

\noindent  (4.4.17) $\qquad \displaystyle
\alpha_{K}(z)\,\,\,=\,\,\, {\rm min} \, \biggl[ \, 1 \,,\, k \,\, {{ {\rm
min} \bigl( \, M_{K}(z) - z_{K} \,,\, z_{K}  - m_{K}(z) \, \bigr)  }\over{ {\rm
max} \, \bigl\{ \, \mid \nabla z(K) \, {\scriptstyle \bullet} \, (y_{K,\,f}
 - x_{K} ) \mid \,,\,f \subset \partial K \, \bigr\}  }} \, \biggr] \,. \,$

\smallskip \noindent   $\bullet \qquad \,\,\, $ 
In the one dimensional case with a regular mesh, it is an exercice to re-write the
extrapolation (4.4.15) under the usual form (4.4.3)  in the context of finite
differences. In this particular case, some limiter functions $\,\, r \longmapsto
\varphi_{k}(r) \,\,$ associated with particular  parameters $\,\, k \,\,$ are shown 
on Figure 4.3.  For $\,\, k=1 ,\,\,$ we recover the initial limiter proposed by Van
Leer in the fourth paper of the family   ``Towards the ultimate finite difference
scheme...''  [VL77]~; for this reason, we have named it the ``Towards~4'' limiter (see
Figure 4.3). When $\,\, k={1\over2}  \,\,$ we obtain the ``min-mod'' limiter proposed
by Harten [Ha83]. The intermediate value $\,\, k={3\over4}  \,\,$ is a good compromise
between the ``nearly  unstable'' choice $\,\, k=1 \,\,$ and the ``too compressive'' 
min-mod choice. We have named it  STS (see also (4.3.14))  and it has been chosen for
our Euler computations in [DM92].

\smallskip \centerline {  \epsfysize=4,7cm  \epsfbox  {fig.4.7.epsf} }  \smallskip

\centerline {\rm  {\bf Figure 4.7.}	\quad  Slope limitation at a fluid boundary. }
\smallskip 
\smallskip \centerline {  \epsfysize=4,7cm  \epsfbox  {fig.4.8.epsf} }  \smallskip

\centerline {\rm  {\bf Figure 4.8.}	\quad   Slope limitation at a solid boundary.}
\bigskip

\noindent   $\bullet \qquad \,\,\, $ 
We explain now the way the preceding scheme is adapted near the boundary. We first 
consider  a {\bf fluid boundary}.  When $\, K \,$ is a finite element with some
face $\, g \subset \partial K \,$  lying on the boundary, we still define the set $\,\,
{\cal N}(K) \,\,$  of   neighbouring cells by the relation (4.4.6) as shown on Figure
4.7. The number of neighbouring cells is just less important in this case. Then points
$\,\, y_{K,\,f}\,\,$ are introduced by relation (4.4.7)  if face $\, f \,$ does
not lie on the boundary and by taking the barycenter of face $\,\, g \,\,$ if it is 
lying on the boundary.  The only difference is the way the values $\,\, 
\overline{z_{K,\,g}} \,\,$ are extrapolated for the face $\,\, g \,\,$ that is on the
boundary~; we set 

\setbox11=\hbox{ $ \displaystyle \,\, 
\overline{z_{K,\,g}} \,\,= \,\,   \, z_{K} \,\,,\quad g \subset \partial K \,\,, $ }
\setbox12=\hbox {  $   \displaystyle  \qquad \qquad   \qquad   \quad   
 g $ face  lying on the boundary of the domain.  }
\setbox40= \vbox {\halign{#&# \cr \box11 \cr \box12 \cr}}
\setbox41= \hbox{ $\vcenter {\box40} $}
\setbox44=\hbox{\noindent  (4.4.18) $\qquad \left\{ \box41 \right. \,  $}  
\noindent $ \box44 $

\noindent 
When values $\,\, \overline{z_{K,\,f}} \,\,$ are determined for all the faces $\, \, 
f \subset \partial K ,\,\,$ the gradient $\, \, \nabla z(K) ,\,\,$  the minimal 
 $\,\,  m_{K}(z) \,\, $ and maximal $\,\,  M_{K}(z) \,\, $ values among the
neighbouring cells  are still determined with the relations (4.4.10), (4.4.12) and
(4.4.13) respectively. The constraints (4.4.16) remain unchanged except that no
limitation process is due to the faces lying on the boundary. In a precise way, we set~: 

\setbox21=\hbox {$\displaystyle \alpha_{K}(z) \,\,=\,\,  \,  $}
\setbox22=\hbox {$\displaystyle  \quad 
{\rm min}  \biggl[  1 ,  { k \quad { {\rm
min} \bigl( \, M_{K}(z) - z_{K} \,,\, z_{K}  - m_{K}(z) \, \bigr)  }\over{ {\rm
max} \, \bigl\{  \abs{ \nabla z(K) {\scriptstyle \bullet} (y_{K,\,f} - x_{K} )
} ,\,f \subset \partial K , \, K_r(f) \in {\cal N}(K) \bigr\}  }}  \biggr] . $}
\setbox30= \vbox {\halign{#&# \cr \box21 \cr \box22    \cr   }}
\setbox31= \hbox{ $\vcenter {\box30} $}
\setbox44=\hbox{\noindent  (4.4.19) $\displaystyle  \,\,   \left\{ \box31 \right.$}  
\noindent $ \box44 $

\noindent
Then the interpolated values $\,\, z_{K,\,f} \,\,$  for all the faces $\,\, f \subset
\partial K \,\,$ are again predicted with the help of relation (4.4.15). 

\noindent   $\bullet \qquad \,\,\, $ 
For a {\bf rigid wall}, the limitation process is a little modified, as presented at
Figure 4.8. We first  introduce the limit face $\,g \,$ inside the set of neighbouring
cells~: 

\setbox21=\hbox {$\displaystyle 
\bigl\{  L \!\in\! {\cal E}_{\cal T} , \, 
\exists f \!\subset \! \partial K \cap \partial L \, \bigr\} \,\, \,  \cup  \,  $}
\setbox22=\hbox {$\displaystyle    \qquad   \qquad  \qquad  
\bigl\{ g \!\in\! {\cal F}_{\cal T},\, g \subset 
\partial K ,\, g $ on the boundary$ \bigr\}  .\, $}
\setbox30= \vbox {\halign{#&# \cr \box21 \cr \box22    \cr   }}
\setbox31= \hbox{ $\vcenter {\box30} $}
\setbox44=\hbox{\noindent  (4.4.20) $\displaystyle  \qquad  
{\cal N}(K) \,\,=\,\, \left. \box31 \right.$}  
\noindent $ \box44 $

\noindent
For the face(s) $\,\,  g \subset \partial K \,\,$ lying on the solid boundary, we
determine preliminary values $\,\,  \overline{z_{K,\,g}} \,\,$ by taking in
consideration    at this level the impenetrability boundary condition $\,\,\, {\bf
u}\, {\scriptstyle \bullet}\, {\bf n}_{g} \,=\, 0 .\,\,$ We introduce the two
components $\,\, n^x_{g} \,$ and $\,\, n^y_{g} \,$ of the normal $\,\, {\bf
n}_{g} \,\,$ at the boundary and we  set, in coherence with variables (4.4.8)~: 

\setbox11=\hbox{   $\!\!  \overline{\rho_{K,\,g}} \,\,= \,\, \rho_{K}\, $} 
\setbox12=\hbox{   $\!\!  \overline{\rho_{K,\,g}} \,\,\,  \overline{u_{K,\,g}}
\,\,= \,\, \rho_{K}\, \bigl( u_{K} - ({\bf u}_{K} \, {\scriptstyle \bullet} \,
{\bf n}_{g}) \,  n^x_{g} \,\bigr) \,  $} 
\setbox13=\hbox{   $\!\!  \overline{\rho_{K,\,g}} \,\,\,  \overline{v_{K,\,g}}
\,\,= \,\, \rho_{K}\, \bigl( v_{K} \,- ({\bf u}_{K} \, {\scriptstyle \bullet} \,
{\bf n}_{g}) \,  n^y_{g} \,\bigr) \,  $} 
\setbox14=\hbox{   $\!\!  \overline{p_{K,\,g}} \,\,= \,\, p_{K}\,.\,  $} 
\setbox40= \vbox {\halign{#&# \cr  \box11 \cr \box12 \cr \box13 \cr \box14 \cr}}
\setbox41= \hbox{ $\vcenter {\box40} $}
\setbox44=\hbox{\noindent  (4.4.21) $\qquad  \displaystyle  \left\{ \box41 \right.\,$}  

\noindent $ \box44 $

\noindent
We consider  also these values for the limitation algorithm. We define ``external
values'' $\,\, z_{L} \,\,$ for $\,\,  L \!=\! g \,\,$ and face $\,\,g \,\, $  lying
on the boundary as equal to the ones defined in relation (4.4.21)~: 

\noindent  (4.4.22) $\quad \displaystyle 
z_{g} \,\equiv \, \overline{z_{K,\,g}} \,,\quad    z({\scriptstyle \bullet})
\, \,$  field defined in (4.4.21),$\,\,\,\,   g \subset  \partial K  $ on the
boundary.

\noindent
Then the extrapolation algorithm that conducts to relation (4.4.15) for extrapolated
values $\,\,  z_{K,\,f} \,\,$ is used as  in the internal case.

\bigskip 
%\vfill \eject
\noindent {\smcaps 4.5 } $ \,  $ { \bf  Explicit Runge-Kutta integration with 
respect to time.}

\noindent   $\bullet \qquad \,\,\, $ 
When all values  $\,\,  z_{K,\,f} \,\,$ are known for all control volumes $\,\, K \in
{\cal E}_{\cal T} ,\,\,$ all faces $\,\, f \!\subset \! K \,\,$ and all fields $\,\,
z({\scriptstyle \bullet}) \,\,$ defined at relation (4.4.8), extrapolated states  
$\,\,  W_{K,\,f} \,\,$ are naturally defined by going back to the conservative
variables. Then we introduce  these states as arguments of the flux function $\,\,
\Phi\bigl(   {\scriptstyle \bullet} \,,\,  {\bf n}_{f} \,,\,  {\scriptstyle \bullet})
\,\,$ and obtain by this way a new system of ordinary differential equations~: 

\noindent  (4.5.1) $\,\,\, \displaystyle 
\abs{K} \,\, {{{\rm d}W_{K}}\over{{\rm d}t}} \,+\, \sum_{f \subset \partial K}
\,\abs{f} \, \Phi\bigl(    W_{K,\,f}  \,,\,  {\bf n}_{f} \,,\,   W_{K_r(f),\,f} 
\bigr)  \,\,= \,\, 0 \,,\qquad   K \in  {\cal E}_{\cal T} \,. \,$ 

\smallskip \noindent   $\bullet \qquad \,\,\, $ 
The numerical integration of such kind of system  is done with a Runge-Kutta
scheme as presented in [CDV92]. We have used with success in [DM92] the Heun scheme of
second order accuracy for discrete integration of (4.5.1) between time steps $\,\,  n\,
\Delta t \,\,$ and $\,\, (n \!+\! 1) \, \Delta t \,\,$~: 

\noindent  (4.5.2) $\,\,\, \displaystyle 
{{\abs{K}}\over{\Delta t}} \,\, \Bigl( \, \widetilde{ W_{K}} -  W_{K}^{n} \, \Bigr)
+ \sum_{f \subset  \partial K} \, \abs{f} \, \Phi \Bigl(    W^{n}_{K,\,f}  \,,\, 
{\bf n}_{f} \,,\,   W^{n}_{K_r(f),\,f}  \Bigr)  \,= \, 0 \,, 
\,   K \in  {\cal E}_{\cal T} \, \,$ 

\noindent  (4.5.3) $\,\,\, \displaystyle 
{{\abs{K}}\over{\Delta t}} \,\, \Bigl( \, \widetilde{ \widetilde{ W_{K}}}  -  
\widetilde{ W_{K}} \, \Bigr)+  \sum_{f \subset  \partial K} \, \abs{f} \,
\Phi \Bigl(   \widetilde{ W}_{K,\,f}  \,,\,  {\bf n}_{f} \,,\,    \widetilde{ W}
_{K_r(f),\,f}  \Bigr) \,= \, 0 \,, 
\,  K \in  {\cal E}_{\cal T} \, \,$ 

\noindent  (4.5.4) $\,\,\, \displaystyle 
W_{K}^{n+1} \,\,= \,\, {1\over2} \,\, \Bigl( \, \widetilde{ \widetilde{ W_{K}}} \, +
\,  W_{K}^{n}  \, \Bigr) \,,\qquad   K \in  {\cal E}_{\cal T} \,. \,$

\bigskip \smallskip
% \bigskip
\noindent  {\smcaps 5) $ \,\,\,\,\,\,\, $   References.} 
\smallskip \noindent

%%%%%%%%%%%%%%%%%%%%%%%%%%%%%%%%%       modif janvier 2011 
%%  \titredroite={\pecaps      References } 

 \hangindent=9mm \hangafter=1 \noindent     [CDV92] $\,\,$
 D. Chargy, F. Dubois, J.P. Vila.
M\'ethodes num\'eriques pour le calcul  d'\'ecoulements compressibles, applications
industrielles, {\it cours de l'Institut pour la Promotion des Sciences de
l'Ing\'enieur}, Paris,  septembre 1992.

\smallskip \hangindent=9mm \hangafter=1 \noindent    [Ci78]  $\,\,$
P.G. Ciarlet. {\it The Finite Element Method
for Elliptic Problems,} North Holland, Amsterdam, 1978. 

\smallskip \hangindent=9mm \hangafter=1 \noindent    [DM92]  $\,\,$
F. Dubois, O. Michaux.    Solution of 
the Euler Equations Around a Double Ellipso\"{\i}dal  Shape Using Unstructured Meshes
and Including Real Gas Effects, {\it Workshop on Hypersonic Flows for Reentry
Problems}, (D\'esid\'eri-Glowinski-P\'eriaux Editors), Springer Verlag, vol. II,
p.~358-373, 1992.

\smallskip \hangindent=9mm \hangafter=1 \noindent    [DM96]  $\,\,$
F. Dubois, G. Mehlman.  A non-parameterized
entropy correction for Roe's approximate Riemann solver, {\it Numerische Mathematik},
vol 73,  p.~169-208, 1996.

\smallskip \hangindent=9mm \hangafter=1 \noindent    [Du91]  $\,\,$
F. Dubois. Nonlinear Interpolation and Total
Variation Diminishing Schemes, {\it  Third International Conference on Hyperbolic
Problems},  (Eng\-quist-Gustafs\-son Editors), Chartwell-Bratt, p.~351-359, 1991.
See also hal-00493555. 

\smallskip \hangindent=9mm \hangafter=1 \noindent    [Du01]  $\,\,$
F. Dubois.  Partial Riemann problem, Boundary
conditions and gas dynamics, in {\it Artificial Boundary Conditions, with Applications
to Computational Fluid Dynamics Problems}, (L. Halpern and L. Tourette Eds), Nova
Science Publishers,  p.~16-77, 2001. See also   hal-00555600. 

\smallskip \hangindent=9mm \hangafter=1 \noindent    [FGH91]  $\,\,$
I. Faille, T. Gallou\"et, R. Herbin. 
Les Math\'ematiciens d\'ecouvrent les Volumes Finis, {\it  Matapli},  
 n$^{\rm o}$~23, p.~37-48,  octobre 1991.  

\smallskip \hangindent=9mm \hangafter=1 \noindent    [GR96]  $\,\,$
 E. Godlewski, P.A. Raviart. {\it
Numerical Approximation of Hyperbolic Systems of Conservation Laws}, Applied
Mathematical Sciences, vol$.\,$118, Springer, New York,  1996.  

\smallskip \hangindent=9mm \hangafter=1 \noindent    [Go59]  $\,\,$
S.K. Godunov. A Difference Method for the
Numerical Computation of Discontinuous Solutions of the Equations of Fluid Dynamics,
{\it Math. Sbornik}, vol. 47, p.~271-290, 1959. 

\smallskip \hangindent=9mm \hangafter=1 \noindent    [GZIKP79]  $\,\,$
 S.K. Godunov, A. Zabrodine, M.
Ivanov, A. Kraiko, G. Prokopov. {\it  R\'esolution num\'erique des probl\`emes
multidimensionnels de la dynamique des gaz,} Editions de Moscou, 1979.  

\smallskip \hangindent=9mm \hangafter=1 \noindent    [Ha83]  $\,\,$
A. Harten.   High Resolution Schemes for
Hyperbolic Conservation Laws, {\it Journal of Computational Physics,} vol. 49, p.
357-393, 1983. 

\smallskip \hangindent=9mm \hangafter=1 \noindent    [HLV83]  $\,\,$
A. Harten, P.D. Lax, B. Van Leer.  
On Upstream Differencing and Godunov-type Schemes for Hyperbolic Conservation Laws,
{\it  SIAM Review,} vol. 25, n$^{\rm o}$~1, p.~35-61, January 1983.  

\smallskip \hangindent=9mm \hangafter=1 \noindent   [HGMW96] $\,\,$
A. Hirschberg, J. Gilbert, R. Msallam, A.P.J. Wijnands. 
Shock waves in trombones, {\it J. Acoust. Soc. Am.}, 
vol.~99,  n$^{\rm o}$3, p.~1754-1758,~1996. 

\smallskip \hangindent=9mm \hangafter=1 \noindent   [Kr70] $\,\,$
H.O. Kreiss. Initial Boundary Value Problems
for Hyperbolic Systems, {\it Comm. Pure Applied Math.}, vol. 23, p.~277-298, 1970.  

\smallskip \hangindent=9mm \hangafter=1 \noindent   [LL54] $\,\,$
L. Landau, E. Lifchitz.  {\it Fluid Mechanics},
Pergamon Press, 1954.

\smallskip \hangindent=9mm \hangafter=1 \noindent  [MD99] $\,\,$
 R. Msallam, F. Dubois. Mathematical
model for coupling a quasi-unidi\-men\-sional perfect flow with an acoustic boundary
layer, {\it Research report  CNAM-IAT} n$^{\rm o}$ 326/99, 1999.
See also hal-00491417. 

\smallskip \hangindent=9mm \hangafter=1 \noindent  [Pa80] $\,\,$
 S.V. Patankar. {\it Numerical Heat Transfer
and Fluid Flow}, Hemisphere publishing, 1980.

\smallskip \hangindent=9mm \hangafter=1 \noindent [Po88] $\,\,$
M. Pollet.  M\'ethodes de calcul relatives aux
interfaces missiles-propul\-seurs, {\it Internal report}, Aerospatiale Les Mureaux,
1988.

\smallskip \hangindent=9mm \hangafter=1 \noindent [Roe81] $\,\,$
P. Roe.   Approximate Riemann Solvers,
Parameter Vectors and Difference Schemes, {\it Journal of Computational Physics,} vol.
43, p.~357-372, 1981.

\smallskip \hangindent=9mm \hangafter=1 \noindent [Roe85] $\,\,$
P. Roe.   Some contributions to the
Modelling of Discontinuous Flows, in {\it Lectures in Applied Mathematics}, vol. 22, 
(Engquist, Osher, Sommerville Eds), AMS, p.~163-193, 1985.

\smallskip \hangindent=9mm \hangafter=1 \noindent [RM67] $\,\,$
  R.D. Richtmyer, K.W. Morton. {\it Difference
Methods for Initial-Value Problems}, Interscience Publishing, J. Wiley \& Sons, New
York, 1967.  

\smallskip \hangindent=9mm \hangafter=1 \noindent [VL77] $\,\,$
B. Van Leer. Towards the Ultimate 
Conservative Difference Scheme IV. A New Approach to Numerical Convection, {\it
Journal of Computational Physics,} vol. 23, p.~276-299, 1977. 

\smallskip \hangindent=9mm \hangafter=1 \noindent [VL79] $\,\,$
 B. Van Leer.   Towards the Ultimate
Conservative Difference Scheme V. A Second Order Sequel to Godunov's Method,
{\it Journal of Computational Physics,} vol 32, n$^0$ 1, p.~101-136, 1979.

%%%%%%%%%%%%%%       mis plus haut, janvier 2011  
%%%      \bigskip  \bigskip \noindent   $\bullet \qquad \,\,\, $ 
%%%      {\it Acknowledgments}. \quad The author thanks Alexandre Gault, listener at the
%%%      spring 2000  ``lectures in computational acoustics'' at the  Conservatoire National des
%%%      Arts et M\'etiers (Paris, France), for providing his personal  manuscript  notes.

\bye